\def\1ton{1,2,\ldots,n}
\def\det{{\rm det}}

\def\div{{\rm div}}
\def\Div{{\rm Div}}

\def\diam{{\rm  diam}}
\def\dist{{\rm dist}}

\def\d{\textnormal d}

\documentclass[11pt]{memo-l}

\usepackage{amssymb}
\usepackage{amsthm}
\usepackage{mathrsfs, amsfonts, amsmath}
\usepackage{graphicx}
\usepackage{psfrag}
\newcommand{\norm}{\,|\!|\,}
\newcommand{\bydef}{\stackrel{{\rm def}}{=\!\!=}}
\newcommand{\rn}{{\mathbb R}^n}
\newcommand{\onto}{\xrightarrow[]{{}_{\!\!\textnormal{onto}\!\!}}}
\newcommand{\into}{\xrightarrow[]{{}_{\!\!\textnormal{into}\!\!}}}
\newcommand{\loc}{\text{loc}}
\newcommand{\Tr}{\text{Tr}}
\newcommand{\A}{\mathbb{A}}
\newcommand{\R}{\mathbb{R}}
\newcommand{\X}{\mathbb{X}}
\newcommand{\rr}{\mathfrak{r}}
\newcommand{\cc}{\mathfrak{c}}
\newcommand{\Y}{\mathbb{Y}}
\newcommand{\W}{\mathscr{W}}

\newcommand{\normav}[1]{\pmb{\boldsymbol{[}} #1 \pmb{\boldsymbol{]}}}
\newcommand{\abs}[1]{\lvert#1\rvert}
\def\Aint{- \nobreak \hskip-10.5pt \nobreak \int}

\newtheorem{theorem}{Theorem}[chapter]
\newtheorem{lemma}[theorem]{Lemma}
\newtheorem{proposition}[theorem]{Proposition}
\newtheorem{corollary}[theorem]{Corollary}

\theoremstyle{definition}
\newtheorem{definition}[theorem]{Definition}
\newtheorem{example}[theorem]{Example}

\theoremstyle{remark}
\newtheorem{remark}[theorem]{Remark}

\def\charfn_#1{{\raise1.2pt\hbox{$\chi_{\kern-1pt\lower3pt\hbox{{$\scriptstyle#1$}}}$}}}

\numberwithin{section}{chapter}
\numberwithin{equation}{chapter}

\begin{document}
\frontmatter

\title{ $n$-Harmonic Mappings Between Annuli \\ \Large {\it The Art of Integrating Free Lagrangians}}

\author{Tadeusz Iwaniec}
\address{Department of Mathematics, Syracuse University, Syracuse,
NY 13244, USA and Department of Mathematics and Statistics,
University of Helsinki, Finland}
\email{tiwaniec@syr.edu}

\author{Jani Onninen}
\address{Department of Mathematics, Syracuse University, Syracuse,
NY 13244, USA}
\email{jkonnine@syr.edu}

\thanks{Iwaniec was supported  by the National Science
Foundation grant DMS-0800416 and the Academy of Finland project 1128331, and Onninen by  the
National Science Foundation grant DMS-1001620. A part of this research was
done when the first author was visiting  the University of Michigan,
as  Gehring visiting professor. He thanks the Department of the
University of Michigan for the support and hospitality.}

\date{\today}



\begin{abstract}
 The central theme of  this paper is the variational analysis of homeomorphisms $h
 \colon  \mathbb X  \onto  \mathbb Y$ between two given domains $\mathbb X , \mathbb Y \subset \mathbb R^n$.  We look for the extremal mappings in the Sobolev space
 $\mathscr W^{1,n}(\mathbb X,\mathbb Y)$ which minimize the energy integral
\[
{\mathscr E}_h=\int_{{\mathbb X}} \norm Dh(x) \norm^n\, \textrm
dx
\]
Because of the  natural connections with quasiconformal mappings this $\,n$-harmonic alternative to
the classical Dirichlet integral (for planar domains)  has drawn the attention of
researchers in Geometric Function Theory. Explicit analysis is made here for a pair of concentric spherical
annuli where many unexpected phenomena about minimal $n$-harmonic mappings are observed. The underlying integration of
nonlinear differential forms, called \textit{free Lagrangians}, becomes truly a work of art.
\end{abstract}

\maketitle{}

\setcounter{page}{4}

\tableofcontents

\chapter*{Preface}

The future developments in modern geometric analysis and its governing partial differential equations (PDEs) will continue to rely on physical and geometric intuition. In recent years, this trend has become more pronounced
and has led to increasing efforts of pure and applied mathematicians, engineers and other
scientists, to share the ideas and problems of compelling interest.  The present paper takes on  concrete
questions about energy minimal deformations of annuli in $\,\mathbb
R^n\,$. We  adopted the
interpretations and central ideas of nonlinear elasticity where the
applied aspects of our results originated. A novelty of our
approach  is that we allow the mappings to slip freely along the boundaries of the domains. It is precisely in this setting that one faces  a real
challenge in establishing the existence, uniqueness  and invertibility
properties of the extremal mappings. The underlying concept  of
 \textit{Free Lagrangians}, is the core of the matter.

Our approach is purely mathematical though the questions
are intimately derived from Nonlinear
Elasticity.  Both the
theoretical and practical aspects of this work  culminate in actual construction of the
mappings with smallest conformal energy.  Special efforts have
been devoted to somewhat subtle  computational details to
present them as simply and clearly as possible.

 We believe the final conclusions  shed considerable new light on the Calculus of Variations, especially for deformations that are free on the boundary. We also feel that some new facts discovered  here have the potential for applications in Geometric Function Theory  as well as  for better  understanding the mathematical models of Nonlinear Elasticity.

\mainmatter

\chapter{Introduction and Overview}
\section{Basic notation}
Let us take a moment to recall a very much needed notation from the calculus of vector fields and matrix fields. Consider
a mapping $h \colon \X \into \Y$ between domains $\X \subset \R^n$ and $\Y \subset \R^m$, $h=(h^1, h^2, \dots, h^m)$, where $h^1, \dots, h^m$ are scalar functions in the Sobolev space $\mathscr W^{1,p}_{\loc} (\X)$. The differential  $Dh(x)$, defined at almost every $x\in \X$, represents a linear transformation of $\R^n$ into $\R^m$, $Dh(x) \colon \R^n \to \R^m$. With the standard choice of the coordinates in $\R^n$ and $\R^m$ we have a matrix field, again denoted by $Dh$,
\begin{equation}
Dh = \left[ \begin{array}{cccc}  h^1_{x_1} &  h^1_{x_2} & \cdots  &  h^1_{x_n}   \\
\vdots & \cdots &   \vdots \\
 h^m_{x_1}  &  h^m_{x_2}  & \cdots  &  h^m_{x_n}
\end{array}\right]\in \mathscr L^p_{\loc} (\X, \R^{m \times n}).
\end{equation}
Hereafter, we abbreviate the notation of distributional partial derivative $\frac{\partial F}{\partial x_i}$ to $F_{x_i}$, $i=1, \dots , n$ for $F \in \mathscr L^1_{\loc} (\mathbb X)$. The differential matrix, also called Jacobian matrix or {\it deformation gradient}, acts on a vector field $V= (V^1, V^2, \dots , V^n) \in \mathscr L^q_{\loc} (\X, \R^n)$ by the rule
\begin{equation}
[Dh]V = \left[ \begin{array}{cccc}  h^1_{x_1} &  h^1_{x_2} & \cdots  &  h^1_{x_n}   \\
\vdots & \cdots &   \vdots \\
 h^m_{x_1}  &  h^m_{x_2}  & \cdots  &  h^m_{x_n}
\end{array}\right]
\left[ \begin{array}{cccc} V^1   \\
\vdots \\
 V^n \end{array}\right] = \left[ \begin{array}{cccc} \langle \nabla h^1 , V \rangle    \\
\vdots \\
 \langle \nabla h^m , V \rangle  \end{array}\right]
\in \mathscr L^1_{\loc} (\X, \R^{m }). \nonumber
\end{equation}
Here $\nabla$ stands for the gradient operator acting on real-valued functions in $\W^{1,p}_{\loc}(\X)$. More generally, consider an arbitrary matrix field
\begin{equation}
M=\left[ \begin{array}{cccc}  M_1^{1}(x) & \cdots  & M_n^{1} (x)  \\
\vdots & \cdots &   \vdots \\
 M_1^{m}(x) & \cdots  & M_n^{m} (x)
\end{array}\right]\in \mathscr L^p_{\loc} (\X, \R^{m \times n})
\end{equation}
and denote its row-vector fields by $\rr^1, \dots , \rr^m \in \mathscr L^p_{\loc} (\X, \R^n)$. Similarly, the column-vector fields will be denoted by $\cc_1, \dots , \cc_n \in \mathscr L^p_{\loc} (\X, \R^m)$. The divergence operator acting on a vector field $\rr=(r_1, \dots , r_n) \in \mathscr L^p_{\loc} (\X, \R^n)$ is a Schwartz distribution defined by
\[\div \,  \rr = \sum_{i=1}^n \frac{\partial r_i}{\partial x_i} \in \mathcal D' (\X, \R).\]
Then the divergence of a matrix field $M\in \mathscr L^p_{\loc} (\X, \R^{m \times n})$ is a distribution, valued in $\R^m$,
\[\Div M = \left[ \begin{array}{cccc} \div \,  \rr^1    \\
\vdots \\
 \div \,  \rr^m  \end{array}\right] \in \mathcal D' (\X, \R^m).   \]
In particular,
\[\Div Dh = \Delta h \in \mathcal D' (\X, \R^m)\]
where $\Delta = \frac{\partial^2}{\partial x_1^2} + \dots + \frac{\partial^2}{\partial x_n^2} $ is the usual Laplacian. The matrix fields $M \in  \mathscr L^p_{\loc} (\X, \R^{m \times n})$ which satisfy the equation $\Div M \equiv 0$ will be called {\it divergence free}, meaning that
\[\int_\X \langle M , D\eta \rangle =0 \quad \mbox{for every test mapping $\eta \in \mathscr C^\infty_\circ (\X, \R^m)$.}\]
Hereafter $\langle A, B \rangle= \Tr (A^\ast B)$ is the inner product of matrices.
We will be typically working with the Hilbert-Schmidt norm of a matrix
\[\norm M \norm^2 = \langle M, M \rangle = \sum_{j=1}^m \sum_{i=1}^n \abs{M_i^j}^2. \]

\section{Mathematical  Model of Hyperelasticity}
Geometric Function Theory (GFT) is currently a field of enormous
activity where  the language  and   general framework  of  Nonlinear
Elasticity  is very helpful. As this
interplay develops, the  $n$-harmonic deformations become  well
acknowledged as a possible generalization of  mappings of finite distortion. We
have also found a place for $n$-harmonic deformations in the theory
of nonlinear hyperelasticity.  J. Ball's fundamental paper \cite{Ba} accounts
for the principles of this theory and sets up mathematical
models.
Historically, the relation between  hyperelasticity and quasiconformal theory, has not been clearly manifested, but  it is indeed very basic and fruitful.

One can roughly describe the hyperelasticity  as a study of weakly differentiable  homeomorphisms  $\,h \colon \mathbb X \overset{\textnormal{\tiny{onto}}}{\longrightarrow} \mathbb Y\;$
between domains in $\,\mathbb R^n\,$ (or $n$-manifolds) that minimize a given energy integral,

\begin{equation}
\mathscr E_h \;=\; \int_\mathbb X \mathbf E(x,h, Dh)\, \textrm{d}x  \; <\infty\quad\quad\quad Dh \colon \mathbb X \rightarrow \mathbb R^{n\times n},
\end{equation}
The condition on the injectivity of $\,h\,$ is imposed in order to avoid  interpenetration of matter.  The Jacobian matrix $\,Dh(x) \in \mathbb R^{n\times n}\,$, defined at almost every point $\,x\in\mathbb X\,$, is referred to as the \textit{deformation gradient.}  In this model the so-called \textit{stored energy function} $\, \mathbf E\colon \mathbb X \times\mathbb Y \times \mathbb R^{n\times n}\rightarrow \mathbb R\,$ is given; it characterizes  mechanical properties of the elastic material in $\,\mathbb X\,$ and $\mathbb Y\,$.

Motivated by GFT we will be essentially concerned with the $n$-harmonic energy, also called {\it conformal energy}
\begin{equation}
\mathscr E_h = \int_\X \norm Dh(x) \norm^n \, \d x.
\end{equation}
Another energy integral of interest in GFT is
\begin{equation}
\mathscr F_h =  \int_\X \frac{ \norm Dh(x)\norm^n }{\abs{h(x)}^n}  \, \d x.
\end{equation}
We devote Chapter~\ref{cha7} to this latter integral.

\section{Variational Integrals in GFT}
 In another direction, we recall Geometric Function Theory in
$\,\mathbb R^n$ and  its governing   variational integrals.
Let as begin with a conformal mapping $\,h: {\mathbb X}
\stackrel{\textrm{\tiny{onto}}}{\longrightarrow}  {\mathbb Y}\, $. Thus  at every $x\in \mathbb X\,$  we have the relation between the norm of the Jacobian matrix  and its determinant    $\norm Dh(x)\norm^n = \,n^{\frac{n}{2}}J(x,h)$. This can be expressed in the form of a  \textit{nonlinear Cauchy-Riemann system }of PDEs;
\begin{equation}
       D^\ast h \cdot Dh  =  J(x,h)^{\frac{2}{n}} \;\textbf{I}\;\;
\end{equation}
It is evident   that the $\,n$ -harmonic energy of  $\,h: \mathbb X
\stackrel{\textrm{\tiny{onto}}}{\longrightarrow }\mathbb Y$ depends only on the deformed configuration. Indeed, we have
\begin{equation}
{\mathscr E}_h=\int_{{\mathbb X}} \norm Dh(x) \norm^n\, \textrm dx
\;=\;n^{\frac{n}{2}}\int_{{\mathbb X}} J(x,h)\, \textrm dx
\;=\,\,n^{\frac{n}{2}}|\mathbb Y|
\end{equation}
For other homeomorphisms $\,g: {\mathbb X}
\stackrel{\textrm{\tiny{onto}}}{\longrightarrow}  {\mathbb Y}$, in the Sobolev space  $\,  \mathscr W^{1,n}(\mathbb X,
\mathbb Y)\,$, we  only have a lower bound, due to Hadamard's inequality for determinants:
$${\mathscr E}_g = \int_\mathbb X \norm Dg (x)\norm^n \, \d x \geqslant  n^{\frac{n}{2}}\int_\mathbb X J(x,g)\, \textrm dx = n^\frac{n}{2}|{\mathbb Y}| $$
Thus
conformal deformations $\,h: {\mathbb X}
\stackrel{\textrm{\tiny{onto}}}{\longrightarrow}  {\mathbb Y}\, $ are none other than those having
  the $\,n$-harmonic energy equal to $n^\frac{n}{2}|\mathbb Y|$, the smallest possible. It is for this reason that conformal mappings are frequently characterized as  {\it absolute minimizers} of the $n$-harmonic integral. However, it is rare
in higher dimensions that two topologically equivalent domains are
conformally equivalent, because of Liouville's rigidity theorem. Even in the plane, multiply connected domains like annuli are of various conformal type.  From this point of concerns Quasiconformal
Theory \cite{AIMb,  IMb} offers significantly larger class of mappings.
\begin{definition}
 A homeomorphism   $\,h: {\mathbb X}
 \to   \mathbb R^n\, $ of   Sobolev space  $ \mathscr W^{1,1}_{\textrm{loc}}(\mathbb X, \mathbb R^n)\,$
  is
said to have {\it finite outer distortion}  if $$\norm Dh(x) \norm^n \leqslant n^\frac{n}{2} \, K(x) \, J(x,h) \; \; \;  $$
for some measurable function $1 \leqslant K(x) < \infty$. The smallest such $\,K(x)\,$  is called the outer distortion, denoted by $\, \mathbb K_{_O}(x,h)\,$. Then $\,h\,$ is  $\,K$ -\textit{quasiconformal} if  $\mathbb K_{_O}(x,h) \leqslant K$ for some constant $K$.
\end{definition}

A concept somewhat dual to outer distortion  is the inner distortion. For this,  consider the cofactor matrix $D^\sharp h$ (it represents infinitesimal deformations of $\,(n-1)$-dimensional entities)  defined for invertible $\,Dh\,$ via  Cramer's rule \[ D^\sharp h \cdot D^\ast h = J(x,h)\, \mathbf I\]

  For a map $h\in \mathscr W^{1,1}_{\textrm{loc}}(\mathbb X, \mathbb R^n)$, not necessarily a homeomorphism, but with nonnegative Jacobian,  we introduce  the \textit{ inner distortion function}
\begin{equation} \mathbb K_{_I}(x,h)= \begin{cases}   \frac{\| D^\sharp h(x)\|^n }{ n^{n/2}\;[J(x,h)\, ]^{n-1}} & \; \textnormal{ if } \; J(x,h)>0 \\
1 &  \; \textnormal{ if } \; J(x,h)=0 \end{cases}\end{equation}

\begin{remark}
Any map of finite outer distortion has finite inner distortion and $\mathbb K_{_I}(x) \leqslant \mathbb K_{_O}(x)^{n-1}$\,, but not vice versa.
\end{remark}
 It is again interesting to find the place for
such mappings in continuum mechanics. The latter deals with the positive definite matrix $
\textbf{C}(x) = D^\ast h(x)\, Dh(x)$   as the
\textit{right Cauchy-Green} deformation tensor. While on the other hand,
there is a fundamental interplay between mappings of finite distortion and  the \textit{Beltrami equation}
\begin{equation}\label{3}
  D^\ast h(x)\, Dh(x) \;=\; J(x,h)^{\frac{2}{n}}\; \mathbf G(x)\;,\;\;\;\;\;\mathrm{det}\;\mathbf G(x) \equiv 1
\end{equation}
  Thus $\mathbf G=\mathbf G(x)$, called the \textit{distortion
  tensor} of $h$, is none other than the    Cauchy-Green tensor renormalized
   so as to have  determinant  identically equal to one.
The symmetric positive definite matrix function $\mathbf G =\mathbf G(x) = [G_{ij}(x)] \in
  \mathbb R^{n\times n}\,$  can be viewed as a Riemann metric tensor on $\mathbb X$.
  In this way  $h$ becomes conformal with respect
  to this, usually only measurable, metric structure on $\mathbb X$. Thus ${\bf G}(x)$ is uniformly elliptic in case of $K$-quasiconformal mappings.

 It is in this  Riemannian manifold framework   that
  variational interpretations  of quasiconformal mappings really crystalize.   For example, the solutions to the Beltrami equation (\ref{3})
   are none other then the {\it absolute  minimizers} of their own
  energy integrals.   Indeed, a homeomorphism
$h: \mathbb X \rightarrow \mathbb Y$ of  Sobolev class ${\mathscr
W}^{1,n} (\mathbb X, \mathbb Y)$ solves the Beltrami equation
(\ref{3})  if and only if
\begin{equation}\label{xx}
{\mathscr E}_h \bydef \int_{{\mathbb X}} \mathbf E(x,  Dh)
\,\textrm dx \;=\;n^{\frac{n}{2}}\int_{{\mathbb X}} J(x,h)\, \textrm
dx \;=\,\,n^{\frac{n}{2}}|\mathbb Y|
\end{equation}
where  the integrand is defined on  $ \mathbb X\!\times\!\mathbb R^{n\times n}$ by the rule
$$
\mathbf E(x,  \xi) =  \;\left(\,  \textrm {Tr}\,[\, \xi \;\textbf{G}^{-1}(x)\;
\xi^{\ast}\, ]\, \right)^\frac{n}{2}\;,\quad  \xi \in \mathbb R^{n\times n}  $$
As in the conformal case, for all other homeomorphisms $\,g \colon  \mathbb X \stackrel{\textrm{\tiny{onto}}}{\longrightarrow}\,\mathbb Y$, in the Sobolev space  $\mathscr W^{1,n} (\mathbb X , \mathbb Y)$, we have the lower bound $\mathscr E_g \geqslant n^\frac{n}{2}[\mathbb Y]$.
  The most appealing conclusion  is a connection between the $n$ -harmonic energy of
  $\, h : \,\mathbb X \rightarrow\,\mathbb Y\,$
  and  the inner distortion function
 of the inverse mapping  $\, f = h^{-1} : \,\mathbb Y \rightarrow\,\mathbb X\,$.

\begin{proposition}\label{9x}$\textnormal{\textsc{(Transition to the inverse map)}}$ \\ Let $f\in\mathscr W^{1,n-1}_{\textnormal loc}(\mathbb Y,\mathbb X)$ be a
homeomorphism of finite outer distortion between bounded  domains,   with $\mathbb K_{_I}(y,f) \in \mathscr L^1(\mathbb Y)$. Then the inverse map
$\,h = f^{-1}: \mathbb X \stackrel{\textrm{\tiny{onto}}}{\longrightarrow}\,\mathbb Y \,$ belongs to the
 Sobolev class $\,\mathscr W^{1,n}(\mathbb X,\mathbb Y)\,$ and we have
\begin{equation}\label{11}
   n^{\frac{n}{2}}\int_\mathbb Y  \mathbb K_{\!_I}(y, f) \textrm dy \;\;\;=\;\;\int _\mathbb X \| Dh(x)\| ^n \textrm dx
\end{equation}
 \end{proposition}
This identity gains in significance if we realize that the polyconvex variational integrand in the left hand side turns into a convex one, a rarity that one can exploit when studying  quasiconformal mappings of smallest $\,\mathscr L ^1$ -mean distortion.
From yet another perspective, it is worth mentioning the   classical Teichm\"{u}ller theory which
is concerned, broadly speaking, with extremal mappings between
Riemann surfaces.  The extremal Teichm\"{u}ller mappings are exactly the ones
whose distortion function has the smallest possible $\,\mathscr L^\infty$-norm. The existence and uniqueness of such an extremal quasiconformal map within a given homotopy class of quasiconformal mappings is the heart  of Teichm\"uller's theory.
Now, in view of the identity (\ref{11}),   minimizing the
$\,\mathscr L^1$-norm of the inner distortion  function offers a
study of $n$-harmonic mappings. Is there any better motivation?

\section{Conformal Energy}

For $\mathbb X\,$ and $\,\mathbb Y\,$  open regions in $\rn$, we
shall consider mappings
\begin{equation}
h=\left(h^1,..., h^n\right) : \; \mathbb X  \longrightarrow \mathbb
Y
\end{equation}
in the Sobolev class  ${\mathscr W}^{1,n} (\mathbb X, \mathbb Y)$.
Thus the Jacobian matrix of $h$ and its determinant are well defined
at almost every point $x\in \mathbb X$. We recall the notation
\begin{equation}
Dh\, = \, \left[\frac{\partial h^i}{\partial x_j}\right]\in
{\mathscr L^n(\mathbb X,\mathbb R}^{n \times n})\, , \hskip0.5cm
J(x,h)=\mbox{det}\, Dh\;\in\mathscr L^1(\mathbb X)
\end{equation}
Here, as usual, $ {\mathbb R}^{n \times n}$ is supplied with the inner product and the
Hilbert-Schmidt norm:
\begin{equation}
\langle A, B \rangle \, = \, \mbox{Tr}\, \left(A^\star B\right) =
\sum_{ij=1}^n A_j^iB_j^i\;\;\;\;\;\;\;\;\;\;\;\norm A \norm =
\langle A, A \rangle^\frac{1}{2}
\end{equation}
  At the initial stage of our undertaking the  \textit{$n$-harmonic integral} will be subjected to the orientation preserving homeomorphisms
   $\,h
 : \; \mathbb X  \stackrel{\textrm{\tiny{onto}}}{\longrightarrow} \mathbb Y\,$, so that
\begin{equation}\label{xxxx}
{\mathscr E}_h=\int_\mathbb X \norm Dh(x) \norm^n\,
dx\;\geqslant\;n^{\frac{n}{2}}\int_\mathbb X J(x,h)\,
dx\;=\;n^{\frac{n}{2}} |\mathbb Y|
\end{equation}
In dimension $n \geqslant 3$, it may well be that no homeomorphism  $h: {\mathbb X}
\stackrel{\textrm{onto}}{\longrightarrow} \mathbb Y$ of finite
$n$-harmonic energy exists, as the following result~\cite{IOj} shows.

\begin{theorem}\label{blahthm}
Let $\,\mathbb X \subset\mathbb R^n\,$ be a ball with a
$\,k$-dimensional closed disk removed, and let $\,\mathbb Y \subset
\mathbb R^n\,$ be a ball with a $\,(k+1)$-dimensional closed disk
removed, $1\leqslant k< n-1$. Then every homeomorphism $h:
{\mathbb X} \stackrel{\textrm{onto}}{\longrightarrow} \mathbb Y$ has
infinite $\,n$-harmonic energy.
\end{theorem}
Note that both $\X$ and $\Y$ are {\it topological annuli}; that is, homeomorphic images of a spherical annulus $\A=\{x \colon r< \abs{x}<R\}$. Let us  view the disks removed from the balls as cracks. It can be
easily  shown, by means of an example,  that mappings of finite conformal energy may
outstretch an $\,(n-1)$-dimensional crack into an $\,n$-dimensional
hole. However, Theorem~\ref{blahthm}  ensures us  that, in principle,
mappings of finite energy cannot increase the dimension of lower
dimensional cracks, a fact highly nontrivial to observe and prove
\cite{IOj}. From now on we assume without explicit mention that the domains $\,\mathbb X\,$
and $\,\mathbb Y\,$ admit at least one homeomorphism $\,h
 : \; \mathbb X  \stackrel{\textrm{\tiny{onto}}}{\longrightarrow} \mathbb Y\,$ in the Sobolev space
 $\,\mathscr W^{1,n}(\mathbb X,\mathbb Y)$
\section{Weak limits of homeomorphisms}

But the true challenge is to find a deformation $\,h^\circ: {\mathbb
X} \stackrel{\textrm{onto}}{\longrightarrow} \mathbb Y$ with the
smallest possible energy. In general, when passing to the weak limit
of the minimizing sequence of homeomorphisms, the injectivity of the extremal map will
be lost. Nevertheless, from the point of view of the elasticity
theory~\cite{An, Ba, Ci}, such limits are still
legitimate deformations to consider. For, if this is the case, they create no new cracks or holes
in $\;\mathbb Y\,$. Let $\,\mathscr P(\mathbb X,\mathbb Y)\,$ denote
the class of weak limits of  homeomorphisms $\,h
 : \; \mathbb X  \stackrel{\textrm{\tiny{onto}}}{\longrightarrow} \mathbb Y\,$ in the Sobolev space
 $\,\mathscr W^{1,n}(\mathbb X,\mathbb Y)$. We refer to such limits as \textit{permissible
 deformations}.\footnote{Homeomorphisms converging weakly in $\,\mathscr W^{1,n}(\mathbb X,\mathbb Y)$ also converge
$c$-uniformly, so
 their limits are still continuous, taking $\,\mathbb X\,$ into $\,\overline{\mathbb Y}\,$.}
 Differentiability and  geometric features of permissible mappings are  not as clear as one may have expected.
In Theorems~\ref{thme14},~\ref{thme12} and~\ref{ThMono} we assume that $\X$ and $\Y$ are bounded domains of the same topological type, like spherical annuli, having at least two though finitely many boundary components. We consider a sequence $h_j \colon \mathbb X \onto \mathbb Y$ of homeomorphisms converging weakly in $\W^{1,n} (\mathbb X , \mathbb Y)$ to a mapping $h \colon \mathbb X \to \overline{\mathbb Y}$. In general,  homeomorphisms $h_j: {\mathbb X} \onto {\mathbb Y}$ do not
 extend as  continuous maps to  the closure of $\,\mathbb
X\,$, but the distance functions $x \to \dist \big(h_j (x), \partial \mathbb Y\big)$ do extend. This is also true for the limit mapping $h$. The precise result to which we are referring is the
following:
 \begin{theorem}\label{thme14}  \cite[Theorem 1.1]{IOt}  For the above-mentioned pair of domains $\X$ and $\mathbb Y$, there  exists a nonnegative continuous function $\eta=\eta (x)$ defined on $\overline{\mathbb X}$ such that
\begin{equation}\label{distest}
\dist \big(h_j(x), \partial \mathbb Y \big) \leqslant \eta (x)  \norm Dh_j \norm_{\mathscr L^n (\mathbb X)}, \qquad \eta \equiv 0 \mbox{ on } \partial \X.
\end{equation}
In particular,
\begin{equation}
\dist \big(h(x), \partial \mathbb Y \big) \leqslant \eta (x) \sup_{j \geqslant 1} \norm Dh_j \norm_{\mathscr L^n (\mathbb X)}.
\end{equation}
\end{theorem}
The weak limit $h$ actually covers the target domain, but  this may  fail if $\mathbb X$ and $\mathbb Y$ have only one boundary component.
\begin{theorem}\label{thme12} \cite[Theorem 1.4]{IOt}
The mapping $h$ is continuous and $\mathbb Y \subset h(\mathbb X) \subset \overline{\mathbb Y}$. Furthermore, there exists a measurable mapping $f \colon \mathbb Y \to \mathbb X$, such that
\[h \circ f = \textnormal{id} \,: \mathbb Y\rightarrow \mathbb Y, \]
everywhere on $\mathbb Y$. This right inverse mapping has bounded variation,
\[ \norm f\norm_{_{\textnormal{BV}(\mathbb
  Y)}}
  \;\leqslant \;\int_\mathbb X \norm Dh(x)\norm ^{n-1} \;\textnormal dx .
\]
\end{theorem}
As noted in~\cite[Remark 9.1]{IOt} the weak limit $h$ is monotone in the  sense of C.B.
Morrey \cite{Mor}.
\begin{theorem}\label{ThMono}
The mapping $h$ is   monotone, meaning that for every continuum $\mathbb K
\subset \mathbb Y\,$ its preimage $\;h^{-1}(\mathbb K) \subset
\mathbb X\,$ is also a continuum; that is compact and connected.
 \end{theorem}
The proof of this  theorem is presented in Section~\ref{secprmon}.

\section{Annuli}

The first nontrivial case is that of doubly connected domains. Thus we consider mappings $\, h:\mathbb A \rightarrow \mathbb
A^{\!\ast}\,$ between concentric spherical annuli in $\R^n$.
\begin{equation}
{\mathbb A}= \mathbb A(r,R) =\left\{x\in \mathbb R^n\; ; \; \; \; r
< |x| < R \right\} \;,\;\;\;\;\;\;\; 0\leqslant r < R < \infty
\nonumber\end{equation}
 \begin{equation} {\mathbb A^{\!\ast} }= \mathbb
A(r_{_{\!\!\ast}}\,,R_{_{\ast}}) =\left\{y \in \mathbb R^n \; ; \;
\; \; r_{_{\!\!\ast}} < |y| < R_{_{\ast}}\right\}\;,\;\;\;
0\leqslant r_{_{\! \! \ast}} < R_{_{\ast}} < \infty \nonumber
\end{equation}
Such domains   are of different conformal type
unless the ratio of the two radii is the same for both annuli. As
for the domains of higher connectivity in dimension $\,n=2\,$, the
conformal type of a domain of connectivity $\,\ell > 2\,$ is
determined by $\,3\ell - 6\,$ parameters, called Riemann moduli of
the domain.\footnote{In this context the mappings are orientation
preserving.} This means that two $\,\ell$-connected domains are
conformally equivalent if and only if they agree in all
$\,3\ell-6\,$ moduli. But we shall have considerably more
freedom in deforming $\,\mathbb X\,$ onto $\,\mathbb Y$, simply by means of
mappings of finite energy. An obvious question to ask is whether
minimization of the $\;n$-harmonic integral is possible within
homeomorphisms between domains of different conformal type.

Concerning uniqueness, we note that the energy  ${\mathscr
E}_h$ is invariant under conformal change of
the variable $x\in {\mathbb A}$. Such a change of variable is
realized by {\it conformal automorphism} of the form
\begin{equation}
x^\prime = \left(\frac{rR}{|x|^2}\right)^k\, Tx
\end{equation}
where $k=0,1$ and $T$ is an orthogonal matrix.
\section{Hammering  a part of an annulus into a circle, $n=2$}

Let us caution the reader that a minimizer $h^\circ \colon \mathbb A
\rightarrow \mathbb A^{\!\ast}$, among all permissible deformations,
does not necessarily satisfy the Laplace equation. A loss of harmonicity  occurs exactly at the points where $\,h^\circ\,$ fails
to be injective. This is the case when the target annulus ${\mathbb
A}^\ast$ is too thin as compared with ${\mathbb A}$; precisely, if
\begin{equation}\label{XX}
\frac{R_{\,{\!\ast}}}{r_{{\!\ast}}} \;< \;\frac{1}{2}\left(
\frac{R}{r} +\frac{r}{R} \right) \, , \;\;\;\;\;\textit{-annuli below the Nitsche bound}
\end{equation}
By way of illustration, consider the so-called \textit{critical Nitsche map}
 \begin{equation}\label{XXXX}
  \aleph(z) \;=\; \frac{1}{2} \left( z\,+\,\frac{1}{\overline z}\right)\;,\;\;\;\;\;0 <|z|<\infty
 \end{equation}
This harmonic mapping takes an annulus $\,\mathbb A(1,R)\,$  univalently  onto
$\,\mathbb A^{\!\ast} = \mathbb
 A(1,R_{\,{\!\ast}})\,$, where $\,R_{\,{\!\ast}} =  \frac{1}{2}\left( R + \frac{1}{R}\right)$. We have equality at
 (\ref{XX}),  and $\aleph$ is  the energy minimizer.   Note the symmetry $\aleph \left(\frac{1}{\overline{z}}\right)= \aleph (z)$. Thus  the same Nitsche map takes reflected annulus $\,\mathbb A(R^{-1} , 1)\,$
 univalently onto $\,\mathbb A^{\!\ast}\,$.  Let us  paste these two annuli along their common boundary
 $$
 \mathbb A \bydef \mathbb A(r,R) \;=\mathbb A(r,1] \cup \mathbb A[1,R)\;,\;\;\;\; r= \frac{1}{R}
 $$
 Now the same harmonic map $\,\aleph :
 \mathbb A  \stackrel{\textrm{\tiny{onto}}}{\longrightarrow} \mathbb A\, [1,R_\ast )$ is a double
cover. Its Jacobian determinant vanishes along the unit circle, the
branch set of $\aleph$. Therefore, $\aleph$ is
  not permissible (it is not a weak $\W^{1,2}$-limit of homeomorphisms). An extension of $\,\aleph :
  \mathbb A(1,R)  \stackrel{\textrm{\tiny{onto}}}{\longrightarrow} \mathbb A^{\!\ast}\,$ inside the unit disk to a permissible mapping of $\mathbb A (r,R)$ onto $\mathbb A^\ast$ can  be nicely
 facilitated by squeezing    $\,\mathbb A(r,1)\,$  onto the unit circle. This procedure will hereafter be referred to as hammering the inner portion of the domain annulus onto the inner boundary of the target.
 Precisely, the map we are referring to takes the form

$$
h^\circ(z) =H(|z|)\, \frac{z}{|z|} \bydef  \left\{\begin{array}{lll} \frac{z}{|z|} \; \;\;&\;\;\;\frac{1}{R} <|z|\leqslant 1\;,\;\;\;&\;\;\;\textrm{hammering part}\\
\\
 \frac{1}{2} \left( z\,+\,\frac{1}{\overline z}\right)\; & \;\;\;\;1\leqslant |z|\leqslant R \;,\;\;\;&\;\;\;\textrm{harmonic part} \end{array} \right.$$

\begin{figure}[!h]
\begin{center}
\psfrag{H}{\small {${|h^\circ |}$}} \psfrag{t}{\small ${t}$}
\psfrag{rh}{\small ${1}$} \psfrag{r}{\small ${r}$}
\psfrag{r'}{\small ${1}$} \psfrag{R}{\small ${R}$}
\psfrag{R'}{\small ${R_\ast}$} \psfrag{h}{$h^\circ$}
\psfrag{A}{${\mathbb A}$} \psfrag{A'}{${\mathbb A}^\ast$}
\includegraphics*[height=1.2in]{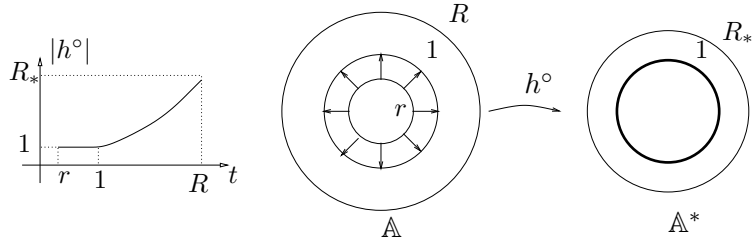}
\caption{Hammering the inner ring into the unit circle.}\label{figur1}
\end{center}
\end{figure}

It is true, though  somewhat less obvious,  that: $\,h^{\circ}\,$ is a $\mathscr W^{1,2}$-limit of homeomorphisms $h \colon \mathbb A \onto \mathbb A^\ast$ and its energy is smaller than that of any homeomorphism from $\mathbb A$ onto $\mathbb A^\ast$, see Figure~\ref{figur1}.

We shall actually prove the following theorem.
\begin{theorem}\label{thn=2}
Let $\mathbb A= \mathbb A(r,R)$ and $\mathbb A^\ast = \mathbb A(r_\ast, R_\ast)$ be planar annuli, $0<r<R<\infty$ and $0<r_\ast<R_\ast<\infty$. We have:\\
\underline{Case 1.}  (Within the Nitsche bound) If
\[\frac{R_\ast}{r_\ast} \geqslant \frac{1}{2} \left(\frac{R}{r}+ \frac{r}{R}\right)\]
then the harmonic homeomorphism
\[h^\circ (z)= \frac{r_\ast}{2} \left(\frac{z}{r}+ \frac{r}{\bar z}\right), \quad h^\circ \colon \mathbb A \onto \mathbb A^\ast\]
attains the smallest energy among all homeomorphisms $h \colon \mathbb A \onto \mathbb A^\ast$, and as such is unique up to a conformal automorphism of $\mathbb A$.\\
\underline{Case 2.} (Below the Nitsche bound) If
\[\frac{R_\ast}{r_\ast} < \frac{1}{2} \left(\frac{R}{r}+ \frac{r}{R}\right)\]
then the infimum energy among all homeomorphisms $h \colon \mathbb A \onto \mathbb A^\ast$ is not attained. Let a radius $r< \sigma < R$ be determined by the equation
\[\frac{R_\ast}{r_\ast}= \frac{1}{2} \left(\frac{R}{\sigma}+ \frac{\sigma}{R}\right) \quad \mbox{-critial Nitsche configuration.}\]
Then the following mapping
\[h^\circ(z) =  \left\{\begin{array}{lll} r_\ast \frac{z}{|z|} \; \;\;&\;\;\; r <|z|\leqslant \sigma \\ \\
 \frac{r_\ast}{2} \left(\frac{z}{\sigma} + \frac{\sigma}{\bar z}  \right)\; & \;\;\;\;\sigma \leqslant |z|<R  \end{array} \right.\]
is a $\W^{1,2}$-limit of homeomorphisms $h_j \colon \mathbb A \onto \mathbb A^\ast$, and its energy is smaller than that of any homeomorphism $h \colon \mathbb A \onto \mathbb A^\ast$.
\end{theorem}
The proof of this theorem was first given in~\cite{AIM}; here in Section~\ref{secn2} we present another one, based on free Lagrangians.

 Let us emphasize that this result does not rule out the existence of
univalent harmonic mappings from $\,\mathbb A\,$ onto $\,\mathbb
A^{\!\ast}\,$, simply because harmonic homeomorphisms need not be the ones that minimize the energy.  Nonexistence of harmonic homeomorphisms between such annuli $\mathbb A$ and $\mathbb A^\ast$ was conjectured by J. C. C. Nitsche~\cite{Nitsche}. After several partial results were obtained in~\cite{Ha, Ka, L, W}, the conjecture  was proved in~\cite{IKOni1, IKOni}. The connection between the
Nitsche conjecture and minimal surfaces is further explored in~\cite{IKOge}. Hammering also occurs in free boundary problems for minimal graphs, where it is called {\it edge-creeping}~\cite{Ch, HN, Tu}.  Similar  hammering phenomena  will be
observed in higher dimensions as well.

\section{Principal $n$-harmonics}

In studying the extremal deformations between spherical annuli it is
natural to look for the radially symmetric solutions of the
$n$-harmonic equation
\begin{equation}\label{Ist}
\mbox{div}\, \left(\, \norm Dh \norm^{n-2}Dh\right)=0\, , \; \; \;
\; \; h(x)=H\big(|x|\big)\, \frac{x}{|x|}
\end{equation}
There is a nice reduction of this problem to the first order
(nonlinear) differential equation for the strain function $H=H(t)$,
called the  {\it characteristic equation} for $H$,
\begin{equation}\label{Idost}
{\mathcal L} H= \left[H^2+\frac{t^2
\dot{H}^2}{n-1}\right]^\frac{n-2}{2}\, \left(H^2-t^2
\dot{H}^2\right) \equiv \mbox{const.}
\end{equation}
Although this equation provides a very convenient tool for studying
properties of radial $n$-harmonics, a little caution is needed
because ${\mathscr C}^1$-solutions to (\ref{Idost}) may fail to
satisfy the original equation (\ref{Ist}). We shall distinguish four so-called
{\it principal $n$-harmonics} $\aleph_\circ\, , \, \aleph^\circ \, , \,
\aleph_+$ and $\aleph_-$. The first two are conformal mappings
\begin{equation}
\aleph_\circ (x)=x \; \; \; \mbox{ and }\; \; \; \aleph^\circ
(x)=\frac{x}{|x|^2}
\end{equation}
The other two are more involved, but still radially symmetric
\begin{equation}
\aleph_+ =H_+\big(|x|\big)\frac{x}{|x|} \; \; \; \mbox{ and }\; \;
\; \aleph_- =H_-\big(|x|\big)\frac{x}{|x|}
\end{equation}
They cannot be described in any elementary way, except for the case
$n=2$.
\begin{equation}
\aleph_\circ (z)=z\, , \; \aleph^\circ (z)= \frac{1}{\overline{z}}
\end{equation}
\begin{equation}\label{kaava22}
 \aleph_+(z)= \frac{1}{2}\left(z+
\frac{1}{\overline{z}}\right)\; \mbox{ and }\; \aleph_-(z)=
\frac{1}{2}\left(z- \frac{1}{\overline{z}}\right)
\end{equation}
The principal $n$-harmonics $\aleph_+$ and $\aleph_-$ are determined
 by solving the following Cauchy problems for their strain functions
\begin{eqnarray}
{\mathcal L}H_+ &\equiv & 1 \hskip1cm H_+(1)=1\\
{\mathcal L}H_- &\equiv & -1 \hskip1cm H_-(1)=0
\end{eqnarray}
It has to be emphasized that all radial solutions to the
$n$-Laplace equation take the form
\begin{equation}\label{Itrist}
h(x)=\lambda\, \aleph(kx)\, , \; \; \; \; \; \; \lambda \in {\mathbb
R}\, ,\;  k>0
\end{equation}
where $\aleph$ is one of the four principal solutions. We then see
that every radial solution (originally defined in an annulus)
 extends $n$-harmonically to the entire punctured space
${\mathbb R}^n_\circ=\rn \setminus\{0\}$, and is evidently
${\mathscr C}^\infty$-smooth. Their connections with~\eqref{kaava22} motivate our calling
$\aleph_\circ\, , \, \aleph^\circ \, , \, \aleph_+$ and $\aleph_-$
the
Nitsche maps in $\rn$.
\section{Elasticity of stretching}

There is an important entity associated with the radial mappings;
namely the elasticity of stretching
\begin{equation}\label{IE17}
\eta_{_H}(t)= \frac{t\, \dot{H}(t)}{H(t)}\, , \; \; \; \mbox{
provided } H^2 + \dot{H}^2 \not=0
\end{equation}
The elasticity function of a power stretching $H(t)=t^\alpha$ is a constant,
equal to $ \alpha$. Two particular cases $\eta \equiv 1$ and
$\eta \equiv -1$ correspond to the conformal  maps
$\aleph_\circ$ and $\aleph^\circ$. In both cases ${\mathcal L}H
\equiv 0$. There are exactly three types of radial $n$-harmonics:
\begin{itemize}
\item[(a)]{{\it Conformally Expanding.} These are the mappings with $\left|\eta_{_H}\right|>1$, everywhere. Equivalently, ${\mathcal L}H \equiv \mbox{const}>0$. Geometrically it means that $h$ exhibits  greater change in the radial direction
than in spherical directions.}
\item[(b)]{{\it Conformally Contracting.}  These are the mappings with $\left|\eta_{_H}\right|<1$, everywhere. Equivalently, ${\mathcal L}H \equiv \mbox{const}<0$.}
\item[(c)]{{\it Conformally Balanced.}  These are conformal  mappings with $\left|\eta_{_H}\right|\equiv 1$, or,  equivalently ${\mathcal L}H \equiv 0$.}
\end{itemize}
For every permissible radial map $h(x)=H\big(|x|\big)\frac{x}{|x|}\, :\;
{\mathbb A} \rightarrow {\mathbb A}^\ast$ (weak $\mathscr W^{1,n}$-limit of homeomorphisms) the elasticity function
does not change sign. Indeed, this follows from the identity $J(x,h)= \dot{H}\big(\abs{x}\big) \cdot \left(\frac{H\big(\abs{x}\big)}{\abs{x}}\right)^{n-1}$. The following ${\mathscr L}^1$-mean (with respect to the
conformal density on ${\mathbb A}$) is  equal to the ratio of the
moduli of the annuli, see (\ref{1321})  for the definition of the modulus,

\begin{equation}
\Aint_{\mathbb A} \left|\eta_{_H}(x)\right|\, \d \mu (x) =
\frac{\mbox{Mod}\, {\mathbb A}^\ast}{\mbox{Mod}\, {\mathbb A}}\, ,
\; \; \; \; \d  \mu = \frac{\d x}{|x|^n}\; \; \;  \footnote{The
integral mean notation $\underset{\mathbb E}{- \nobreak \hskip-7pt
\nobreak \int} f(x)\, \d \mu (x)$ stands for the ratio
$\int_{\mathbb E}f\, \d \mu \big/ \int_{\mathbb E} \d \mu$.}
\end{equation}
In this way we are led to three types of pairs of the
annuli ${\mathbb A}$ and ${\mathbb A}^\ast$.
\begin{itemize}
\item[(1)]{{\it Conformally expanding}; it  pertains to a pair of annuli such that $$\mbox{Mod}\, {\mathbb A}^\ast > \mbox{Mod}\, {\mathbb A}$$}
\item[(2)]{{\it Conformally contracting}; it pertains to a pair of annuli such that $$\mbox{Mod}\, {\mathbb A}^\ast < \mbox{Mod}\, {\mathbb A}$$}
\item[(3)]{{\it Conformally equivalent}; these are the annuli having the same modulus $$\mbox{Mod}\, {\mathbb A}^\ast = \mbox{Mod}\, {\mathbb A}$$}
\end{itemize}
These three cases will be treated by using somewhat different
estimates.
\section{Conformally expanding pair}

\begin{figure}[!h]
\begin{center}
\psfrag{t}{\small${ t}$} \psfrag{tn}{\small${ \;}$}
\psfrag{N}{\small $\aleph_\dag(t)$} \psfrag{Ga}{\small
$\aleph^\dag(t)$} \psfrag{45}{\small ${\tiny 45^\circ}$}
\includegraphics*[height=1.6in]{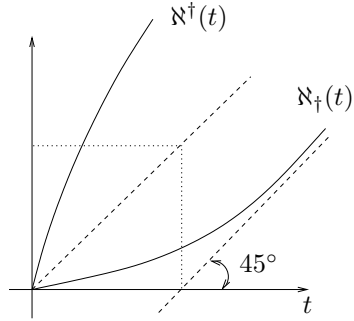}
\caption{The Nitsche functions.}\label{figu2}
\end{center}
\end{figure}

We shall show that the principal solution
$\aleph_-(x)=H_-\big(|x|\big)\frac{x}{|x|}$ generates all
minimizers of the conformal energy. It is rather easy to show that
for a given expanding pair ${\mathbb A}$, ${\mathbb A}^\ast$ there exist
unique $k>0$ and $\lambda
>0$ such that the $n$-harmonic map $h^\circ (x)=\lambda \,
\aleph_-(kx)$ takes ${\mathbb A}$ homeomorphically onto ${\mathbb
A}^\ast$. However, the answer to the question  whether  this map minimizes the $n$-harmonic energy among all homeomorphisms is not obvious. When $n=2$
or $n=3$, the answer is "yes".
\begin{theorem}\label{ThNe1}
Let $\mbox{Mod}\, \mathbb A^\ast > \mbox{Mod}\, \mathbb A $. Then for  $n=2,3$, the $n$-harmonic radial map $h^\circ= \lambda \aleph_- (kx)$ assumes the minimum  conformal energy within all homeomorphisms. Such a minimizer is unique up to a conformal automorphism of $\mathbb A$.
\end{theorem}
Surprisingly, for $n \geqslant 4$ the answer will depend on how wide is the target annulus
${\mathbb A}^\ast$, relatively to $\mathbb A$.
\begin{theorem}\label{ThNe2}
For dimensions $n \geqslant 4$, there exists a function ${\mathcal
N}^\dag = {\mathcal N}^\dag (t)$, $t< {\mathcal N}^\dag (t) < \infty$ for $t>0$, see Figure~\ref{figu2}, such that: if
\begin{equation}\label{n4}
\textnormal{Mod}\, {\mathbb A} \leqslant \textnormal{Mod}\, {\mathbb A}^\ast  \leqslant  {\mathcal
N}^\dag(\textnormal{Mod}\, {\mathbb A}) \quad \mbox{-upper Nitsche bound for  $n \geqslant 4$,}
\end{equation}
then the  map $h^\circ :
{\mathbb A} \rightarrow {\mathbb A}^\ast$ is a unique (up to an automorphism of $\mathbb A$)   minimizer of the
conformal energy among all homeomorphisms.
\end{theorem}
\begin{theorem}\label{n4no}
In dimensions $n \geqslant 4$, there are annuli $\mathbb A$ and $\mathbb A^\ast$ such  that no radial stretching  from $\mathbb A$ onto $\mathbb A^\ast$ minimizes the conformal energy.
\end{theorem}


\section{Conformally contracting pair}

In this case we  obtain the minimizers from the
principal solution $\aleph_+=H_+\big(|x|\big)\frac{x}{|x|}$. As
before, we observe that the mappings
\begin{equation}
h^\circ (x)= \lambda \, H_+(kx)\, , \; \; \; k>0\, , \; \;
\lambda >0
\end{equation}
are radial $n$-harmonics. Recall that ${\mathbb A}^\ast$ is
conformally thinner than ${\mathbb A}$. But it is not enough.  In contrast to the
previous case, such $n$-harmonic mappings  take
the annulus ${\mathbb A}$ homeomorphically onto ${\mathbb A}^\ast$
only when ${\mathbb A}^\ast$ is not too thin. The precise necessary
condition reads as
\begin{equation}\label{LoNi}
\mbox{the  lower Nitsche bound;} \qquad  {\mathcal N}_\dag(\mbox{Mod}\, {\mathbb A}) \leqslant \mbox{Mod}\,
{\mathbb A}^\ast \leqslant \mbox{Mod}\, {\mathbb A}.
\end{equation}
Numerically, the lower Nitsche function ${\mathcal N}_\dag$ is given by \begin{equation}\label{equ34p}0<\aleph_\dag (t) = \omega_{n-1}\, \log
H_+\left(\textnormal{exp} \frac{t}{\omega_{n-1}}\right)<t \quad  \mbox{ for }\;   0<t<\infty. \end{equation}

\begin{theorem}\label{ThNe3}
Under the condition at (\ref{LoNi}) there exist unique $k>0$ and
$\lambda
>0$ such that $h^\circ (x)=\lambda H_+ (kx)$ takes ${\mathbb
A}$ homeomorphically onto ${\mathbb A}^\ast$. This map is a unique
 (up to a conformal automorphism of ${\mathbb A}$) minimizer   of the
conformal energy among all  homeomorphisms of $\mathbb A$ onto $\mathbb A^\ast$.
\end{theorem}

Rather unexpectedly, the injectivity of the weak limit of a minimizing sequence of homeomorphisms $h_j \colon \mathbb A \onto \mathbb A^\ast$ will fail once the lower bound in~\eqref{LoNi} is violated.  Let us look more closely at the
critical   configuration of annuli; that is,
\begin{equation}
{\mathcal N}_\dag(\mbox{Mod}\, {\mathbb A}^\prime ) = \mbox{Mod}\,
{\mathbb A}^\ast\, , \; \; \; \; \; \;
 \left\{\begin{array}{ll}& \textnormal{where } {\mathbb A}^\prime ={\mathbb A}(1,R) \textnormal{ and}\\
& {\mathbb A}^\ast ={\mathbb A}(1,R_\ast)\, , \; \; R_\ast =H_+(R)
\end{array}\right.
\end{equation}
Thus the critical Nitsche map $\aleph_+ (x)= H_+(|x|)\,
\frac{x}{|x|}$, defined for $1\leqslant |x| \leqslant R$ takes
${\mathbb A}^\prime$ homeomorphically onto ${\mathbb A}^\ast$. The Jacobian determinant of $\aleph$ vanishes for $|x|=1$. Now
let us built a pair of annuli below the lower bound at (\ref{LoNi}), simply  by pasting an
additional spherical ring to ${\mathbb A}^\prime$ along the unit sphere. Thus,  we consider a slightly fatter
annulus
$${\mathbb A}={\mathbb A}\, (r,R)={\mathbb A}(r,1]\cup {\mathbb A}[1,R)$$
As in dimension $n=2$, we have
\begin{theorem}\label{ThNe4}
The following deformation
\begin{equation}
h^\circ(x) \;\bydef\; \left\{\begin{array}{lll} \frac{x}{|x|} \; \; \; \; & r <|x|\leqslant 1\; , \;\;\; & \textrm{ hammering part}\\
 \aleph_+ (x) & 1\leqslant |x|\leqslant R\; , & \textrm{ $n$-harmonic part} \end{array} \right.
\end{equation}
is a $\W^{1,n}$-limit of homeomorphisms $h \colon \mathbb A \onto \mathbb A^\ast$ and its energy is smaller than that of any homeomorphism from $\mathbb A$ onto $\mathbb A^\ast$. Among all such mappings $h^\circ$ is unique up to conformal automorphisms of $\mathbb A$.
\end{theorem}
\section{The Conformal Case $\textnormal{Mod}\, {\mathbb
A} = \textnormal{Mod}\, {\mathbb A}^\ast$}

Obviously the minimizers are conformal mappings. Precisely, they
take the form
\begin{equation}
h^\circ (x)= \frac{\sqrt{r_\ast
R_\ast}}{|x|}\left(\frac{\sqrt{rR}}{|x|}\right)^{\pm 1}Tx
\end{equation}
where $T$ is an orthogonal transformation. This case receives
additional consideration in Section \ref{SecQc}.

\section{The energy function $\mathscr F_h$}

So far we considered the domains ${\mathbb X}$ and ${\mathbb Y}$ equipped
with the Euclidean metric. However, one may ask analogous questions
for different metrics on ${\mathbb X}$ and ${\mathbb Y}$. We do not
 enter into a general framework here, but instead  illustrate this possibility by introducing a conformal density  on the target space,
which we continue to  assume to be a spherical annulus
\begin{equation}
{\mathbb Y}={\mathbb A}^\ast = {\mathbb A}(r_\ast \, , \, R_\ast)\,
, \; \; \mbox{ equipped with the measure }\;  \d \mu (y) = \frac{\d
y}{|y|^n}
\end{equation}
The modulus of ${\mathbb A}^\ast$ is none other than
its "conformal" volume
\begin{equation}\label{1321}
\mbox{Mod}\, {\mathbb A}^\ast = \int_{\mathbb A^\ast} \, \d \mu (y)=\int_{{\mathbb A}^\ast}\frac{\d
y}{|y|^n}=\omega_{n-1}\log \frac{R_\ast}{r_\ast}
\end{equation}
More generally, if  $h
 : \; \mathbb X  \stackrel{\textrm{\tiny{onto}}}{\longrightarrow} {\mathbb
 A}^\ast$ is conformal then the modulus of ${\mathbb X}$ is defined by pulling back the measure $ \d
 \mu$ to ${\mathbb X}$.
 \begin{equation}
 \mbox{Mod}\, {\mathbb X} = \int_{\mathbb X} \frac{J(x,h)\, \d
 x}{|h(x)|^n}=n^{-\frac{n}{2}}\int_{\mathbb X} \frac{\norm Dh(x)\norm^n \, \d x}{|h(x)|^n}
 \end{equation}
Obviously, this definition is free from the choice of the conformal map $\,h\,$. For $n=2$ the modulus of a doubly connected domain is the
only conformal invariant; that is,  the Riemann moduli space is one
dimensional. The situation is much more rigid for $n \geqslant 3$.
If now  $h
 : \; \mathbb X  \stackrel{\textrm{\tiny{onto}}}{\longrightarrow} {\mathbb
 A}^\ast$ is any permissible deformation, then
\begin{equation}
{\mathscr F}_h \bydef \int_{\mathbb X} \frac{\norm Dh(x) \norm^n \,
\d x}{|h(x)|^n} \geqslant n^\frac{n}{2} \int_{\mathbb X} \frac{J(x,h)\, \d
x}{|h(x)|^n}=n^\frac{n}{2}\mbox{Mod}\, {\mathbb A}^\ast
\end{equation}
Naturally, this integral tells us  about how much $h$ differs   from a conformal mapping  in an average sense. On the other hand, Quasiconformal
Theory deals with point-wise distortions. Among them is the {\it
outer distortion function}
\begin{equation}
{\mathcal K}_{_O} (x,h) \;\bydef\; \left\{\begin{array}{ll}
\frac{\norm Dh(x) \norm^n}{n^\frac{n}{2}J(x,h)}\, , \; \; \; \; & \textnormal{ if }\; J(x,h)>0\\
  1\, , & \textnormal{ otherwise }\end{array} \right.
\end{equation}
Let us push forward ${\mathcal K}_{_O}$ to the target space via the
mapping $h$ itself, so as to obtain a function
\begin{equation}
{\mathcal K}_h(y)\; \bydef \; {\mathcal K}_{_O} (x,h)\, , \; \; \;
\; \; \mbox{ where }\; x=h^{-1}(y)
\end{equation}
We note, without proof, that
\begin{equation}
{\mathscr F}_h=n^\frac{n}{2}\int_{\mathbb Y} {\mathcal K}_h(y) \, \d
\mu (y)
\end{equation}
This discussion leads us to the minimization of the
${\mathscr L}^1(\mathbb Y, \d \mu)$-norm of the outer distortion. One might suspect that the minimum will be
attained when the distortion function is constant. Indeed, when the domain
is also an annulus, say ${\mathbb X}={\mathbb A}= {\mathbb A}(r,R)$,
then the power stretching
\begin{equation}
h^\alpha (x)=\lambda |x|^{\alpha-1}x\, , \; \; \; \; \; \mbox{ where
}\; \; \alpha = \frac{\mbox{Mod}\, {\mathbb A}^\ast}{\mbox{Mod}\,
{\mathbb A}} \mbox{ and } \lambda =r_\ast r^{-\alpha}=R_\ast
R^{-\alpha}
\end{equation}
has  constant outer distortion ${\mathcal K}_{_O}(x,
h^\alpha)=\alpha^{-1}n^{-\frac{n}{2}}\, (\alpha^2+n-1)^\frac{n}{2}$.
Nevertheless, it takes some effort to show that in dimensions $n=2,3$,  $h^\alpha$  is actually  a
minimizer of $\mathscr F_h$ among all homeomorphisms.
\begin{theorem}\label{th3}
Let ${\mathbb A}$ and ${\mathbb A}^\ast$ be spherical annuli in
${\mathbb R}^n$, $n=2,3$. Then for every homeomorphism $h\colon \mathbb A \onto \mathbb A^\ast$  we have
\begin{equation}\label{IE9}
{\mathscr F}_h \bydef  \int_{{\mathbb A}} \frac{\norm Dh \norm^n}{|h|^n}
\geqslant  \left(n-1 + \alpha^2\right)^\frac{n}{2}
\textnormal{Mod}\, {\mathbb A}\, , \; \; \; \;  \textnormal{ where }
\alpha = \frac{ \textnormal{Mod}\, {\mathbb A}^\ast}{
\textnormal{Mod}\, {\mathbb A}}
\end{equation}
Equality holds for the power stretching $h (x) =r_\ast \,
r^{-\alpha} |x|^{\alpha -1} x$, uniquely up to conformal
automorphisms of ${\mathbb A}$.
\end{theorem}
It is somewhat surprising that in dimensions $n \geqslant 4$,  this
feature of $h^\alpha$ no longer holds when the target
annulus is conformally too fat. We have the following result.
\begin{theorem}\label{th3half}
For each $n\geqslant 4$, there exists $\alpha_n > 1$ such that
(\ref{IE9}) holds whenever
\begin{equation}\label{eqy146}
\alpha \bydef \frac{ \textnormal{Mod}\, {\mathbb A}^\ast}{
\textnormal{Mod}\, {\mathbb A}} \,<\, \alpha_n
\end{equation}
The power stretching $h (x) =r_\ast r^{-\alpha} |x|^{\alpha -1} x$
is the only minimizer of $\mathscr F_h$  modulo conformal automorphism of ${\mathbb A}$.
\end{theorem}
Examples will be given to show that the extremals are no longer
power stretchings if
\begin{equation}
\frac{\mbox{Mod}\, {\mathbb A}^\ast}{\mbox{Mod}\, {\mathbb A}}
\geqslant \sqrt{\frac{n-1}{n-3}}
\end{equation}
In other words, the upper bound for $\alpha$ in~\eqref{eqy146} lies in the interval
$1< \alpha_n < \sqrt{\frac{n-1}{n-3}}$. See Chapter \ref{Sec33} for
more precise estimates of $\alpha_n$. Moreover, if $\mbox{Mod}\, \mathbb A^\ast$ is too large relative to $\mbox{Mod}\, \mathbb A$, then
the extremals cannot be found even  within general radial mappings, see Chapter \ref{SecCro}.
\section{Free Lagrangians}

In 1977 a novel approach towards minimization of polyconvex energy
functionals was developed and published by J. Ball \cite{Ba}. The
underlying idea was to view the integrand as convex function of null
Lagrangians. The term null Lagrangian pertains to a nonlinear
differential expression whose integral over any open region depends
only on the boundary values of the map, like integrals of an exact
differential form. The interested reader is referred to~\cite{BCO, Ed, Iw}. But we are concerned with mappings $\,h :
\mathbb X \rightarrow \mathbb Y $ that are free on the boundary. The
only condition we impose on $h$ is that it is a weak ${\mathscr
W}^{1,n}$-limit of homeomorphisms from $\mathbb X$ onto $\mathbb Y$.
There still exist some nonlinear differential forms, associated with a given pair of domains $\mathbb X$ and $\mathbb Y$,  whose integral means over $\mathbb X$ remain the same within a given class of deformations $h \colon \mathbb X \onto \mathbb Y$, regardless of their boundary values. These
are rather special null Lagrangians. The simplest example is
furnished by the Jacobian determinant of an orientation preserving  homeomorphism $h\in \W^{1,n}(\X, \mathbb Y)$
\begin{equation}
\int_{\mathbb X}J(x,h)\, \d x = |{\mathbb Y}|.
\end{equation}
One might ask which expressions of the type $E(x,h,Dh)$ enjoy
identities such as this? We call them {\it free Lagrangians}. Such a
notion lies fairly deep in the topology  of the mappings $h:
{\mathbb X} \rightarrow {\mathbb Y}$. For example, using the
topological degree, we find that the differential expression
\begin{equation}
E(x,h,Dh)\, \d x= \sum_{i=1}^n \frac{h^i\, dh^1 \wedge ... \wedge
dh^{i-1}\wedge d|x|\wedge dh^{i+1}\wedge ... \wedge d h^n}{|x|\,
|h|^n}
\end{equation}
is a free Lagrangian within the class of orientation preserving homeomorphisms $h:
{\mathbb A} \rightarrow {\mathbb A}^\ast$ between annuli. Indeed, we
have the desired identity
\begin{equation}\label{Inint}
\int_{\mathbb A}E(x,h,Dh)\, \d x \; = \; \mbox{Mod}\, {\mathbb A} \;
\; \; \; \mbox{ for every }\; h\in {\mathscr P}({\mathbb X},
{\mathbb Y})
\end{equation}
A peculiarity of this example is further emphasized by the fact that
the target annulus ${\mathbb A}^\ast$ does not even enter into this
identity.

Like in the theory of  polyconvex energies, the minimization of an energy integral whose  integrand is a convex function of a number of free Lagrangians   poses no challenge; Jensen's
inequality usually gives the desired sharp lower bounds. Unfortunately, our
integrand $\norm Dh \norm^n$ cannot be expressed as a convex
function of free Lagrangians; though it is a convex function of the
usual null Lagrangians.
\section{Uniqueness}

To reach the uniqueness conclusions, we first  show that any extremal
mapping $h : {\mathbb A} \rightarrow {\mathbb A}^\ast$ satisfies the
following system of nonlinear PDEs
\begin{equation}
D^\ast h(x)\, Dh(x)={\bf G}(x,h)
\end{equation}
where the  Cauchy-Green tensor {\bf G} actually depends  only on  two variables $|x|$ and $|h|$.
As this system is overdetermined it comes as no
surprise that any two solutions $h^\circ (x)$ and $h(x)$  for which  $|h^\circ (a)|=|h(a)|$, at some  $a\in {\mathbb A}$, must be equal
modulo  an orthogonal transformation $T: {\rn} \rightarrow {\rn}$,
namely $h(x)=Th^\circ (x)$. Our proof exploits the classical
computation of curvature of ${\bf G}$ in terms of its Christoffel
symbols. As regards the existence of such point $a\in {\mathbb A}$,
we shall
again rely on estimates of free Lagrangians.
\section{The ${\mathscr L}^1$-theory of inner distortion}

The conformal energy can naturally  be turned around so as to yield, the ${\mathscr L}^1$-integrability of a
distortion function of the inverse mapping \cite{AIMO}.

We recall that a homeomorphism $f:{\mathbb Y} \rightarrow {\mathbb X} $ of Sobolev
class ${\mathscr W}^{1,1}_{\loc} ({\mathbb Y}, {\mathbb X})$ is said
to have finite distortion if there exists a measurable function
$1\leqslant K(y) < \infty$ such that
\begin{equation}\label{E47}
|Df(y)|^n \leqslant K(y)\, J(y,f)
\end{equation}
for almost every $y\in {\mathbb Y}$. Here $\abs{Df}$ stands for the operator norm of $Df$. Using the Hilbert-Schmidt norm
of matrices the outer distortion
function of $f$ takes the form
\begin{equation}\label{E48}
{\mathbb K}_{_O}   (y,f)=\frac{\norm
Df(y)\norm^n}{n^\frac{n}{2}J(y,f)}
\end{equation}
if $ J(y,f)>0 $ and we set ${\mathbb K}_{_O}(y , f) =1$, otherwise.
There are many more distortion functions of great importance in
Geometric Function Theory.
Among them are the inner distortion functions.
\begin{equation}
{K}_{_I}   (y,f) = \frac{| D^\sharp f(y)|^n}{\det D^\sharp f(y)}\;
\; \footnote{Here we use the operator norm of the cofactor matrix $D^\sharp
f$.}
\end{equation}
and
\begin{equation}
{\mathbb K}_{_I}   (y,f) = \frac{\norm D^\sharp
f(y)\norm^n}{n^\frac{n}{2}\det  D^\sharp f(y)}
\end{equation}
at the points where $J(y,f)>0$. Otherwise, we set
\begin{equation}
{\mathbb K}_{_I}   (y,f) = {K}_{_I}   (y,f) =1
\end{equation}
In recent years there has been substantial interest in mappings with integrable distortion \cite{AIM, AIMO,  HKO, IKO1, IS}. Suppose now that
$f\in {\mathscr W}^{1,n-1}_{\loc}({\mathbb Y}, \rn)$. The following
identities connect the $n$-harmonic integrals with the theory of mappings with integrable inner distortion
\begin{equation}\label{E48}
\int_{{\mathbb Y}}{K}_{_I}   (y,f)\, \d y = \int_{\mathbb X} | D
h(x)|^n\, \d x
\end{equation}
and
\begin{equation}\label{E49}
n^\frac{n}{2}\int_{\mathbb Y}{\mathbb K}_{_I}   (y,f)\, \d y = \int_{\mathbb X}
\norm D h(x)\norm^n\, \d x
\end{equation}
where $h$ denotes the inverse of $f$, \cite{CHM, HK, HKO}. These identities imply that the inverse of a mapping of
integrable distortion always lies in the Sobolev space $h\in
{\mathscr W}^{1,n} ({\mathbb X}, {\mathbb Y})$.  Similarly, the  $\mathscr L^1 (\mathbb A^\ast, \d \mu)$-integral means of ${K}_{_I} (y,f)$ and ${\mathbb
K}_{_I} \, (y,f)$ with respect to the dimensionless weight $\d \mu
=|y|^{-n}\, \d y$ are the energy functionals $\mathscr F_h$:
\begin{equation}\label{E50}
\int_{{\mathbb A}^\ast}\frac{{K}_{_I} (y,f)}{|y|^n}\, \d y =
\int_{\mathbb A} \frac{| D h(x)|^n}{|h(x)|^n}\, dx
\end{equation}
and
\begin{equation}\label{E51}
n^\frac{n}{2}\int_{{\mathbb A}^\ast}\frac{{\mathbb K}_{_I} (y,f)}{|y|^n}\, \d y =
\int_{\mathbb A} \frac{\norm D h(x)\norm^n}{|h(x)|^n}\, dx
\end{equation}
Returning to~\eqref{E49}, we infer from minimal $n$-harmonic mappings that

\begin{theorem}\label{th5}
Under the Nitsche bounds
\begin{equation}
\aleph_\dag (\textnormal{Mod}\, {\mathbb A}) \leqslant
\textnormal{Mod}\, {\mathbb A}^\ast \leqslant \aleph^\dag
(\textnormal{Mod}\, {\mathbb A})
\end{equation}
the $\mathscr L^1 (\mathbb Y)$-norm of the inner distortion  ${\mathbb K}_{_I} (y,f)$
assumes its minimum value on a mapping $f \colon \mathbb A^\ast \onto \mathbb A$ whose  is  inverse is  a radial $n$-harmonic
mapping $h^\circ \colon \mathbb A \onto \mathbb A^\ast$. Such an  extremal mapping $f$ is unique up to a conformal automorphism  of ${\mathbb A}$.
\end{theorem}
\begin{theorem}\label{th6}
If the domain annulus $\mathbb A^\ast$ is too thin relative to the target annulus $\mathbb A$; precisely, under the condition
\begin{equation}
\textnormal{Mod}\, {\mathbb A}^\ast \; < \;  \aleph_\dag
(\textnormal{Mod}\, {\mathbb A}) \qquad \mbox{-below the Nitsche bound.}
\end{equation}
then the infimum of the $\mathscr L^1(\mathbb A^\ast)$-norm  of   ${\mathbb K}_{_I} (y,f)$ is
not attained among homeomorphisms $f \colon \mathbb A^\ast \onto \mathbb A$.
\end{theorem}
Nevertheless, we were able to find the infimum of the $\,\mathscr L^1$-norms  and
identify the minimizing sequences. The weak BV-limits of such sequences and the underlying concept of their distortion (to be defined) are worth carrying out.

\section*{Conclusion}

At the first glance the problems we  study here may appear to be merely technical. However, their
solutions require truly innovative approaches with surprising outcomes. For instance, the case $n \geqslant 4$ is different than one might
a priori expect; upper bounds for the modulus of ${\mathbb A}^\ast$
are necessary in order to ensure  radial symmetry of the minimizers. The underlying technique of  integration  of various nonlinear
differential forms is interesting in its own right. The {\it
free Lagrangians} play a pivotal role, like  null
Lagrangians did play  in the  polyconvex calculus of
variations.  The  entire subject grew out of fundamental questions
of  Quasiconformal Geometry about mappings of integrable
distortion. The paper also embraces a number of important aspects of the Calculus of Variation.

We present this little theory  here in two parts.

\part{Principal Radial   $n$-Harmonics}

\chapter{Nonexistence of $n$-harmonic homeomorphisms}
That nice smooth domains $\X$ and $\mathbb Y$, such as annuli, may not admit a homeomorphism $h \colon \X \onto \mathbb Y$ of smallest conformal energy is a sequel of even more general observation.
Let $\mathcal A$ be a topological annulus in $\R^n$, we consider all possible $n$-harmonic homeomorphisms $h \colon \mathcal A \to \R^n$. Unconcerned about the energy (finite or infinite) we address the  question; which topological annuli $\mathcal A^\ast$ can be obtained as images of $\mathcal A$ under such mappings? In general, this problem lies beyond our methods.  Even in dimension $n=2$ and when $\mathcal A^\ast$ is a circular annulus (the Nitsche conjecture) the answer to this question required rather sophisticated ideas~\cite{IKOni}. Nevertheless, one may roughly  say that $\mathcal A^\ast$ cannot be very thin. A specific instance  is as follows. Let $\mathbb S$ be a convex $(n-1)$-hypersurface in $\R^n$, given by
\[\mathbb S= \{y\in \R^n \colon F(y)=0\},\]
where $F$ is a $\mathscr C^2$-smooth function in a neighborhood of $\mathbb S$, such that $\nabla F \ne 0$ on $\mathbb S$. We assume that the Hessian matrix
\[\nabla^2 F = \left[\frac{\partial^{\,2} F}{\partial y_i \,\partial y_{\!j}}\right] \quad \mbox{is positive definite on $\mathbb S$.}\]
Consider the $\epsilon$-vicinity of $\mathbb S$; that is,  $\mathbb V_\epsilon = \{y \colon \abs{F(y)} < \epsilon\}$.
\begin{proposition}
If  $\epsilon$ is small enough  then there is no $n$-harmonic homeomorphism $h \colon \mathcal A \into  \mathbb V_\epsilon$ such that $\mathbb S \subset h(\mathcal A )$.
\end{proposition}
\begin{proof}
We shall not give any explicit bound for $\epsilon$, but instead we argue by contradiction. Suppose that for every positive integer $\ell$, thus $\epsilon = 1/\ell$, there exists an $n$-harmonic homeomorphism
\[h_\ell= (h^1_\ell, \dots, h_\ell^n) \colon \mathcal A \onto h_\ell (\mathcal A ), \quad \mbox{such that } \quad \mathbb S \subset h_\ell (\mathcal A )  \subset \mathbb V_\epsilon, \quad \epsilon = \frac{1}{\ell}.\]
We  appeal to the $\mathscr C^{1, \alpha}$-regularity theory of $n$-harmonic mappings~\cite{DM, Mib, Uh, Ur, Wa}. Accordingly, $h_\ell$ are uniformly bounded in $\mathscr C^{1, \alpha}$-norm on every compact subset of $\mathcal A$. Moreover, there is a subsequence, again denoted by $\{h_\ell\}$, which converges  together with the first order derivatives to a mapping $h \colon \mathcal A \to \mathbb S$  uniformly on compact subsets. The limit mapping $h=(h^1, \dots, h^n)$ still satisfies the $n$-harmonic equation. We aim to show that $Dh \equiv 0$ in $\mathcal A$, meaning that $h$ is constant. Recall that $Dh \colon \mathcal A \to \R^{n \times n}$ is (H\"older) continuous. The  computation below is certainly valid in the open region where $Dh \ne 0$, because $h$ is $\mathscr C^\infty$-smooth in such region. We have
\[F\big( h(x) \big)=0 \qquad \mbox{for all } x \in \mathcal A.\]
Applying partial differentiation $\frac{\partial}{\partial x_\nu}$ and chain rule yields
\begin{equation}\label{equ1111}
\sum_{i=1}^n \frac{\partial F}{\partial y_i} \frac{\partial h^i}{\partial x_\nu}=0, \qquad \mbox{for } \nu =1,2, \dots , n.
\end{equation}
On the other hand, the $n$-harmonic equation takes the form
\begin{equation}\label{equ2222}
\sum_{\nu =1}^n \frac{\partial}{\partial x_\nu} \left[\lambda^{n-2} \frac{\partial h^i}{\partial x_\nu}\right]=0 \quad \mbox{for } i=1,2, \dots , n, \quad \mbox{where } \lambda = \norm {Dh(x)} \norm
\end{equation}
or, equivalently
\begin{equation}
\sum_{\nu =1}^n \left[(n-2) \lambda^{n-3} \frac{\partial \lambda}{\partial x_\nu} \frac{\partial h^i}{\partial x_\nu} + \lambda^{n-2} \frac{\partial^2 h^i}{\partial x_\nu \partial x_\nu}\right]=0 \quad i=1,\dots , n.
\end{equation}
We multiply these equations by $\frac{\partial F}{\partial y_i}$, sum them up with respect to $i$, and use the identity~\eqref{equ1111} to obtain
\begin{equation}\label{equ3333}
\lambda^{n-2} \sum_{i=1}^n \sum_{\nu=1}^n  \frac{\partial F}{\partial y_i} \frac{\partial^2 h^i}{\partial x_\nu \partial x_\nu}=0.
\end{equation}
Next we differentiate~\eqref{equ1111} with respect to $\nu$ and sum up the equations,
\begin{equation}\label{equ4444}
\sum_{i=1}^n \sum_{\nu=1}^n \left(\frac{\partial F}{\partial y_i} \frac{\partial^2 h^i}{\partial x_\nu \partial x_\nu} + \sum_{j=1}^n \frac{\partial^2 F}{\partial y_i \partial y_j} \frac{\partial h^j}{\partial x_\nu}   \frac{\partial h^i}{\partial x_\nu}\right)=0.
\end{equation}
Here the double sum for the first term in~\eqref{equ4444} vanishes due to~\eqref{equ3333}, so we have
\[\sum_{\nu =1}^n \left(\sum_{i, j} \frac{\partial^2 F}{\partial y_i \partial y_j} \;\frac{\partial h^i}{\partial x_\nu} \frac{\partial h^j}{\partial x_\nu}  \right)=0.\]
Since the Hessian matrix of $F$ is positive definite, this equation yields $\frac{\partial h^i}{ \partial x_\nu}=0$ for all  $i$, $\nu =1,2, \dots , n$, as desired.

To reach a contradiction we look more closely at the homeomorphisms $h_\ell \colon \mathcal A \onto h_\ell(\mathcal A )$. Choose and fix an $(n-1)$-dimensional hypersurface $\Sigma \subset \mathcal A$ (topological sphere) which separates the boundary components of $\mathcal A$. Its images $ h_\ell (\Sigma) \subset h_\ell(\mathcal A )$ separate the boundary  components of $ h_\ell(\mathcal A )$. Thus, in particular,
\[0< \liminf_{\ell \to \infty} \;[\diam \, h_\ell (\Sigma)] =  \;\diam \, h (\Sigma)\;=0\]
a clear contradiction.
\end{proof}
It is not so clear, however,  how thick $\mathcal A^\ast$ should be to ensure existence of $n$-harmonic homeomorphisms $h \colon \mathcal A \onto \mathcal A^\ast$. In dimension $n=2$ the condition $\mbox{Mod}\,  \mathcal A^\ast \geqslant \mbox{Mod}\,  \mathcal A$ is sufficient~\cite{IKKO} and is not far from being necessary~\cite{IKOaf}. Questions of existence of harmonic diffeomorphisms between surfaces are treated in~\cite{JS}.

\chapter{Generalized $n$-harmonic mappings}
In the classical Dirichlet problem one asks for the energy minimal mapping $h \colon \X \to \R^n$ of the Sobolev class $h\in h_\circ + \W^{1,n}_\circ (\X , \R^n)$ whose boundary values are explicitly prescribed by means of a given mapping $h_\circ \in  \W^{1,n} (\X , \R^n)$.  We refer to the works of C. B. Morrey~\cite{Mo, Mo1},  where a systematic study of variational methods in the theory of harmonic integrals in vectorial case originated, see also~\cite{Mob}.   The   variation  $h \leadsto h\,+ \,\epsilon \eta $,  in which $\eta \in \mathscr C^\infty_\circ (\X , \R^n)$ and $\epsilon \to 0$, leads to the integral form of the familiar $n$-harmonic system of equations
\begin{equation}\label{equa1}
\int_\X \langle \norm Dh\norm^{n-2}Dh , \, D\eta \rangle =0, \quad \mbox{ for every } \eta \in \mathscr C^\infty_\circ (\X , \R^n).
\end{equation}
Equivalently
\begin{equation}\label{equa2}
\Delta_n h = \Div \big( \norm Dh \norm^{n-2}Dh\big)=0, \quad \mbox{in the  sense of distributions}
\end{equation}
or, entry-wise, for $h=(h^1, \dots, h^n)$
\begin{equation}\label{equa6}
\sum_{i=1}^n \big(\norm Dh\norm^{n-2} h^\alpha_{x_i}\big)_{x_i}=0, \quad \alpha =1,2, \dots, n.
\end{equation}
In nonlinear elastostatics the matrix field
\[\mathcal S =\mathcal S h= \norm Dh(x)\norm^{n-2}Dh(x)\]
is known as {\it Piola-Kirchoff} tensor for the energy density function $W(x)=\norm Dh(x)\norm^n$. This tensor represents the stress induced by $h \colon \X \to \R^n$. Mappings of Sobolev class $\W^{1,n}_{\loc}(\X, \R^n)$ that satisfy the equation $\Div \mathcal S h=0$ are called the {\it equilibrium solutions}. Equilibrium solutions in a given class $h_\circ + \W^{1,n}_\circ (\X, \R^n)$  represent unique minimizers within this class. The situation is dramatically different if  we allow $h$ to slip freely along the boundaries. The {\it inner variation} works well in this case. This is simply a change of the independent variable; $h_\epsilon=h \circ \eta_\epsilon $, where  $\eta_\epsilon \colon \X \onto \X$ are $\mathscr C^\infty$-smooth automorphisms of $\X$ onto itself, depending smoothly on a parameter $\epsilon \approx 0$ where  $\eta_\circ = id \colon \X \onto \X$. Let us take on, as a initial step, the inner variation of the form
\begin{equation}\label{equa7}
\eta_\epsilon (x)= x + \epsilon \, \eta (x), \qquad \eta \in \mathscr C_\circ^\infty (\X, \R^n).
\end{equation}
We compute
\[
\begin{split}
Dh_\epsilon (x) &= Dh (x+\epsilon \eta) (I+ \epsilon D\eta)\\
\norm Dh_\epsilon (x)\norm^n & = \norm Dh\norm^n
+ n \epsilon\,  \langle \norm{Dh}\norm ^{n-2} D^\ast h\cdot  Dh\, ,\,  D \eta \rangle + o(\epsilon).
\end{split}
\]
Be aware that in this equation $Dh$ is evaluated at the point $y=x+\epsilon \, \eta (x) \in \X$. Integration with respect to $x$-variable yields a formula for the energy of $h_\epsilon$,
\[\mathscr E_{h_\epsilon} = \int_\X \left[ \norm Dh \norm^n +  n \epsilon \langle \norm{Dh}\norm ^{n-2} D^\ast h\cdot  Dh\, ,\,  D \eta \rangle \right]\, \d x + o(\epsilon) .\]
We now make the substitution $y=x + \epsilon \, \eta (x)$ for which the following transportation rules apply:  $x= y- \epsilon \, \eta (y)+ o(\epsilon)$, $\eta (x)= \eta (y)+o(1)$ and the change of volume element $\d x = [1-\epsilon \, \Tr \,D \eta (y) ]\, \d y + o(\epsilon) $.  The equilibrium equation for the inner variation is obtained from $\frac{\d}{\d \epsilon} \mathscr E_{h_\epsilon}\,=\,0\,$ at $\epsilon =0$,
\begin{equation}\label{intstar}\int_\X \langle \norm Dh \norm^{n-2} D^\ast h \cdot Dh - \frac{1}{n} \norm Dh \norm^n I \, , \, D \eta \rangle \, \d y=0 \end{equation}
or, by means of distributions
\begin{equation}\label{enhe}
\Div \left(\norm Dh \norm^{n-2} D^\ast h \cdot Dh - \frac{1}{n} \norm Dh \norm^n I  \right)=0.
\end{equation}
Now we introduce the divergence free tensor
\[\Lambda = \norm Dh \norm^{n-2} D^\ast h \cdot Dh - \frac{1}{n} \norm Dh \norm^n I =  \left(C - \frac{1}{n} \Tr \,  C \right) \Tr^\frac{n-2}{n}C \]
where we recall the right Cauchy-Green tensor $C=C(x)= D^\ast h \cdot Dh$.
The name {\it generalized $n$-harmonic equation} will be given to~\eqref{enhe} because of the following:
\begin{lemma}
Every $n$-harmonic mapping $h\in \W^{1,n}_{\loc} (\X , \R^n)$ solves the generalized $n$-harmonic equation~\eqref{enhe}.
\end{lemma}
\begin{proof}
We consider again the perturbed mappings $h_\epsilon (x)= h(x+\epsilon \,\eta (x))$ which coincide with $h$ outside a subdomain $\mathbb U$ compactly contained in $\X$. Applying the integral form~\eqref{equa1} of the $n$-harmonic equation,  but with the test function $h-h_\epsilon \in \W^{1,n}_\circ (\mathbb U, \R^n)$ in place of $\eta$, \footnote{This is justified because $\mathscr C^\infty_\circ (\mathbb U , \R^n)$ is dense in $\W^{1,n}_\circ (\mathbb U , \R^n)$.} we estimate the energy of $h$ over $\mathbb U$ as follows
\[
\begin{split}
\mathscr E_h & = \int_{\mathbb U} \norm Dh \norm^n = \int_{\mathbb U} \norm Dh \norm^{n-2} \langle Dh \, , \, Dh \rangle \\
& =  \int_{\mathbb U} \norm Dh \norm^{n-2} \langle Dh \, , \, Dh_\epsilon \rangle \leqslant \left( \int_{\mathbb U} \norm Dh \norm^n \right)^\frac{n-1}{n} \left(\int_{\mathbb U} \norm Dh_\epsilon \norm^n \right)^\frac{1}{n}.
\end{split}
\]
Hence
\[\mathscr E_h \leqslant \mathscr  E_{h_\epsilon}, \quad \mbox{with equality at $\epsilon =0$.} \]
This means $\frac{\d}{\d \epsilon} \mathcal E_{h_\epsilon}=0$ at $\epsilon =0$. Equivalently, $h$ satisfies the generalized $n$-harmonic equation.
\end{proof}
In dimension $n=2$, the generalized harmonic equation reduces to
\begin{equation}\label{e2he}
\Div \left(D^\ast h \, Dh - \frac{1}{2} \norm Dh \norm^2 I  \right)=0.
\end{equation}
This equation is known as Hopf-Laplace equation~\cite{CIKO}. In complex notation it takes the form
\begin{equation}\label{eqv}
\frac{\partial}{\partial \bar z} \left(h_z \overline{h_{\bar z}}\right)=0, \qquad z= x_1+ix_2.
\end{equation}
If, by some reason, $h\in \mathscr C^2 (\X, \mathbb C)$ then~\eqref{eqv} reads as
\begin{equation}\label{eqvi}
J(z,h)\Delta h =0, \qquad \mbox{where } \Delta = 4 \frac{\partial^2}{\partial z \partial \bar z}.
\end{equation}
There exist diverse non-harmonic solutions to~\eqref{eqvi}, many of them  of great interest in the theory of minimal surfaces and some  with potential applications to nonlinear elasticity (elastic plates), see~\cite{CIKO}.

Equally in higher dimensions, the $n$-harmonic mappings of Sobolev class $\W^{1,n}_{\loc} (\X , \R^n)$ are only particular solutions to the generalized $n$-Laplacian.
\section[Non $n$-harmonic solutions]{Solutions to the generalized $n$-harmonic equation that are not $n$-harmonic}

Let us take on stage a radial mapping
\[h(x)= H\big(\abs{x}\big) \frac{x}{\abs{x}}, \qquad \mbox{where $H=H(t)$ is absolutely continuous}\]
see Chapter~\ref{cha4} for extensive treatment of such mappings. We find that
\begin{equation}\label{eqvii}
\begin{split}
\Lambda & = \norm Dh \norm^{n-2} \left(D^\ast h \cdot Dh - \frac{1}{n} \norm Dh \norm^2 I\right)\\
& = (n-1)^\frac{n-2}{n} \left(H^2 + \frac{\abs{x}^2 \dot{H}^2}{n-1}  \right)^\frac{n-2}{2} \left(H^2 - \abs{x}^2 \dot{H}^2\right) \frac{1}{\abs{x}^n} \left(\frac{x \otimes x}{\abs{x}^2}- \frac{1}{n}I\right).
\end{split}
\end{equation}
It is shown in Chapter~\ref{secrad}, that if $h$ is a $\mathscr C^2$-smooth $n$-harmonic mapping then $H=H(t)$ must satisfy the  {\it characteristic equation}
\begin{equation}\label{eqviii}
 \left(H^2 + \frac{\abs{x}^2 \dot{H}^2}{n-1}  \right)^\frac{n-2}{2}\cdot\left(H^2 - |x|^2 \dot{H}^2\right)\equiv \mbox{const.}
\end{equation}
This also confirms that $h$ satisfies the generalized $n$-harmonic equation  because of the identity
\[\Div \,  \frac{1}{\abs{x}^n} \left( \frac{x \otimes x}{\abs{x}^2} - \frac{1}{n}I \right)=0.\]
However, the {\it hammering mapping} $h(x)= \frac{x}{\abs{x}}$, corresponding to $H(t) \equiv 1$, solves equation~\eqref{enhe} but is not $n$-harmonic. If we paste this hammering mapping with a compatible  smooth radial $n$-harmonic solution  (the critical Nitsche mapping) there will emerge a $\mathscr C^{1,1}$-smooth minimal deformation between annuli. There are many more non $n$-harmonic solutions to~\eqref{enhe}, some seem to be unsatisfactory; for example those with Jacobian changing sign are forbidden in elasticity theory because of the principle of non-interpenetration of matter.

\section{Slipping along the boundaries}
Let us return to the inner variation in~\eqref{equa7}, $h_\epsilon (x)=h(x+ \epsilon \, \eta (x))$, but this time with $\eta$ not necessarily having compact support. We assume that $\X$ is a $\mathscr C^1$-smooth  domain and $\eta \in \mathscr C^1 (\overline{\X})$. The fact that $h$ is allowed to freely slip along the boundaries amounts to saying that the vector field $\eta = \eta (x)$ is tangent to $\partial \X$ at every point $x\in \partial \X$. The integral form of the resulting variational equation is the same as in~\eqref{intstar}. Integration by parts (Green's formula) will produce no integral over $\X$, because of~\eqref{enhe}; there will remain only boundary integrals. Precisely, a general formula we are referring to is:
\[0= \int_\X \langle \Lambda, D\eta \rangle \, \d x = - \int_\X \langle \Div \Lambda, \eta \rangle \, \d x + \int_{\partial \X} \langle \Lambda (x) \mathfrak n (x), \eta (x) \rangle \, \d \sigma (x)\]
where $\mathfrak n (x)$ is the outer unit vector field to $\partial \X$ and $\d \sigma (x)$ stands for the surface measure. This is justified under an appropriate assumption on the degree of integrability of $\Lambda$, $\Div \Lambda$, $\eta$ and $D \eta$. Nevertheless in our case, we obtain
\begin{equation}\label{eqsli1}
\int_{\partial \X} \langle ( \norm Dh \norm^{n-2} D^\ast h \cdot Dh - \frac{1}{n} \norm Dh \norm^n I )\mathfrak n(x) \, , \, \eta (x) \rangle  \, \d \sigma (x) =0
\end{equation}
 Since $\eta (x)$ can be any tangent field the equation~\eqref{eqsli1} is possible if and only if the vector field
\[\left(\norm Dh \norm^{n-2} D^\ast h \cdot Dh - \frac{1}{n} \norm Dh \norm^n I  \right) \mathfrak n (x)\]
is orthogonal to  $\partial \X$. In other words, at each $x\in \partial \X$ the linear mapping $\norm Dh \norm^{n-2} D^\ast h \cdot Dh - \frac{1}{n} \norm Dh \norm^n I $ takes the (one-dimensional) space of normal vectors into itself. Of course, the same holds for the mapping $D^\ast h \cdot Dh$. Since $D^\ast h \cdot Dh$ is symmetric one can say, equivalently, that $D^\ast h \cdot Dh$ preserves the tangent space. We just proved the following.
\begin{proposition}
The equilibrium solution for mappings that are slipping along the boundaries satisfies, in addition to~\eqref{enhe}, the following  condition
\begin{equation}\label{bddcond}
D^\ast h \cdot Dh \colon {\bf T}_x \partial \X \to {\bf T}_x \partial \X
\end{equation}
equivalently,
\begin{equation}
D^\ast h \cdot Dh \colon {\bf N}_x \partial \X \to {\bf N}_x \partial \X  .
\end{equation}
where ${\bf T}_x \partial \X$ and  ${\bf N}_x \partial \X$ designate the tangent and normal spaces at $x\in \partial \X$.
\end{proposition}
In dimension $n=2$ this amounts to saying that the Hopf quadratic differential $h_z \overline{h_{\bar z}}\, \d z^2$ is real along the boundary components of $\X$. In an annulus $\mathbb A = \{z \colon r<\abs{z}<R\}$ every Hopf differential which is real on $\partial \mathbb A$ takes the form
\[h_z \overline{h_{\bar z}} \equiv \frac{c}{z^2}, \qquad c \mbox{ is  a real constant}\]
see~\cite{CIKO}. One of the solutions is the mapping we are already encountered in~\eqref{XXXX}
\[
h(z)=  \left\{\begin{array}{lll} \frac{z}{|z|} \; \;\;&\;\;\;\frac{1}{R} <|z|\leqslant 1\;,\;\;\;&\;\;\;\textrm{hammering part}\\
\\
 \frac{1}{2} \left( z\,+\,\frac{1}{\overline z}\right)\; & \;\;\;\;1\leqslant |z|\leqslant R \;,\;\;\;&\;\;\;\textrm{harmonic part} \end{array} \right.
\]
In higher dimensions the radial mapping $h(x)= H \big(\abs{x}\big) \frac{x}{\abs{x}}$ in the annulus $\mathbb A = \{x \colon r<\abs{x}<R\}$ also complies with the boundary condition~\eqref{bddcond}. Indeed, the Cauchy-Green tensor takes the form
\[D^\ast h \cdot Dh = \frac{H^2}{\abs{x}^2}I + \left(\dot{H}^2 - \frac{H^2}{\abs{x}}\right) \frac{x \otimes x}{\abs{x}^2}.\]
The normal vector field at $\partial \mathbb A$ is $\mathfrak n= \mathfrak n (x)= x$, and we see that
\[[D^\ast h(x) \cdot Dh(x)]\, \mathfrak n = \lambda \mathfrak n, \qquad \lambda= [\dot{H}(\abs{x})]^2 .\]
For any tangent vector $\mathfrak t = \mathfrak t (x)$, we have $(x \otimes x) \mathfrak t =0$ so
\[\left[D^\ast h (x) \cdot Dh(x)  \right] \mathfrak t = \frac{H^2}{\abs{x}^2} \mathfrak t.\]
The corresponding singular values of the right Cauchy-Green tensor are  $\dot{H}(x)$ and $\frac{H\big(\abs{x}\big)}{\abs{x}}$; {\it principal stretches}.
\vskip0.3cm
Before moving to a study of $n$-harmonic mappings between spherical annuli, let us fulfill the promise of  proving  that weak $\W^{1,n}$-limits of homeomorphisms are indeed monotone.

\section[Proof of Theorem 1.7]{Proof of Theorem~\ref{ThMono}}\label{secprmon}
We are dealing with a sequence of homeomorphisms $h_j \colon \mathbb X \onto \mathbb Y$, $j=1,2, \dots$ converging $c$-uniformly to a mapping $h \colon \X \to \overline{\mathbb Y}$ such that
\[
\dist \big(h_j(x), \partial \mathbb Y \big) \leqslant \eta (x)  \norm Dh_j \norm_{\mathscr L^n (\mathbb X)}.
\]
see~\eqref{distest}. Passing to the limit we also have
\[
\dist \big(h(x), \partial \mathbb Y \big) \leqslant \eta  (x)  \liminf_{j \to \infty} \norm Dh_j \norm_{\mathscr L^n (\mathbb X)}.
\]
Therefore, for every $\epsilon >0$ there is $\delta >0$ such that
\begin{equation}\label{equys}
\dist \big(h_k(x), \partial \mathbb Y \big) < \epsilon \qquad \mbox{whenever } \dist (x, \partial \X)< \delta
\end{equation}
for every $k=1,2, \dots$ and
\[\dist \big(h(x), \partial \mathbb Y \big) < \epsilon \qquad \mbox{whenever } \dist (x, \partial \X)< \delta .
\]
This shows that  $h^{-1} (\mathbb K) \subset \X$ is compact whenever $\mathbb K \subset \mathbb Y$ is compact. To show that $h$ is also monotone we argue as follows. Consider a continuum $\mathbb K \subset \mathbb Y$ and assume, to the contrary, that $h^{-1} (\mathbb K) = \mathbb A \cup \mathbb B$, where $\mathbb A$ and $\mathbb B$ are nonempty disjoint compact subsets of $\mathbb X$. For every integer $\ell$ such that
\[\frac{1}{\ell} < \epsilon \bydef \frac{1}{2} \dist \big( \mathbb K , \partial \mathbb Y \big)\]
one may consider a neighborhood of $\mathbb K$,
\[\mathbb K_\ell = \{y\in \mathbb Y \colon \dist (y, \mathbb K)< 1/\ell \}.\]
This is an open connected set in $\mathbb Y$. The preimage $h^{-1}(\mathbb K_\ell) \subset \X$ is a neighborhood of the compact set $\mathbb A \cup \mathbb B \subset \X$. Since $h_j \to h$ $c$-uniformly then for sufficiently large $j$, say $j=j_\ell$, we have
\[h^{-1}_{j_\ell} (\mathbb K_\ell) \supset \mathbb A \cup \mathbb B, \qquad \ell \geqslant \frac{1}{\epsilon}. \]
The sets $h^{-1}_{j_\ell} (\mathbb K_\ell)$ are connected. Therefore to every $\ell \geqslant \frac{1}{\epsilon}$ there corresponds a point $x_\ell \in h^{-1}_{j_\ell} (\mathbb K_\ell) \subset \X$
such that
\begin{equation}\label{equyss}
\dist (x_\ell , \mathbb A) = \dist (x_\ell , \mathbb B).
\end{equation}
We also note that
\[\dist \big(h_{j_\ell} (x_\ell), \partial \mathbb Y  \big) \geqslant \dist \big(\mathbb K_\ell , \partial \mathbb Y\big) \geqslant 2 \epsilon - \frac{1}{\ell}> \epsilon. \]
Hence, in view of condition~\eqref{equys}, we have
\[\dist(x_\ell , \partial \X) \geqslant \delta, \qquad \ell =1,2, \dots\]
We then choose a subsequence, again denoted by $x_\ell$, converging to some point $x\in \X$. Passing to the limit in~\eqref{equyss} yields
\[\dist (x, \mathbb A) = \dist (x, \mathbb B), \qquad \mbox{thus } \; x \not\in \mathbb A \cup \mathbb B.\]
On the other hand, since $h_j$ converge to $h$ uniformly on compact subsets, we may pass to the limit with the sequence $h_{j_\ell} (x_\ell) \to h(x)$. Finally, since $h_{j_\ell} (x_\ell)  \in \mathbb K_\ell$ and $\cap \mathbb K_\ell = \mathbb K$, we conclude that $h(x) \in \mathbb K$, meaning that $ x\in h^{-1} (\mathbb K)= \mathbb A \cup \mathbb B$. This contradiction proves Theorem~\ref{ThMono}.

\chapter{Notation}\label{cha4}
This chapter is designed  to describe basic geometric objects and
concepts that will be used throughout  this paper.  The primary
domains here are  annuli in $\rn$, also called spherical rings. They
are subsets of the punctured Euclidean space,
$$
   \mathbb R^n_\circ\;= \mathbb R^n \setminus \{0\}\;=\left\{x=(x_1,...,x_n)\,; \;\,|x| =
   \sqrt{x_1^2 +\cdots+x_n^2} \neq 0\;\right\}\;,\;\;\;\; n\geqslant 2
$$
\section{Annuli and their modulus}
There are four types of annuli, generally denoted by
$\;{\mathbb A} \subset \rn_\circ\;$
\begin{equation}
{\mathbb A}=\left\{\begin{array}{ll} \mathbb A(r,R) = \left\{x\in \rn \; ;  \; \; r < |x| < R \right\}
\;,
\; \;\;\;\; \mbox{where }\; 0\leqslant r<R\leqslant\infty\\
\mathbb A[r,R) = \left\{x\in \rn \; ;  \; \; r \leqslant |x| <R \right\} \;,
\; \;\;\;\; \mbox{where }\; 0< r< R\leqslant\infty\\
\mathbb A(r,R] = \left\{x\in \rn \; ;  \; \; r < |x| \leqslant R \right\} \;,
\; \;\;\;\; \mbox{where }\; 0\leqslant r<R<\infty\\
\mathbb A[r,R] = \left\{x\in \rn \; ;  \; \; r \leqslant |x| \leqslant R \right\} \;, \; \;
\;\;\;\mbox{where }\; 0< r\leqslant R<\infty
\end{array}\right.
\end{equation}
The boundary of each annulus consists of two components (the inner
and the outer sphere) except for the degenerate case $\,\mathbb
A[r,r]$ which reduces to the $(n-1)$-dimensional sphere denoted by
${\mathbb S}^{n-1}_r$.

The  conformal modulus of an annulus is defined by
\begin{equation}
\mbox{Mod}\, {\mathbb A}=\omega_{n-1} \log \frac{R}{r}\;=\;\int_\mathbb A \frac{\textrm
dx}{\,|x|^n}\;\;\;\in[0,\infty]
\end{equation}
Hereafter, $\omega_{n-1}$ denotes the surface area of the unit
sphere $\;{\mathbb S}^{n-1}\subset \rn$.
\begin{equation}
    \omega_{n-1} = \left\{\begin{array}{ll}\frac{2\pi^k}{1\cdot 2\cdots(k-1)} \;,&\;\;\textrm{for}\;\; n=2k \\
            &\\ \frac{2^{k+1}\pi^k}{1\cdot 3\cdots(2k-1)} \;,&\;\;\textrm{for}\;\; n=2k+1\end{array}\right.
\end{equation}
\section{Polar coordinates in $\,\mathbb R^n_\circ\,$}
Polar coordinates are natural to use when working with annuli. We
associate with any point $\,x\in \mathbb R^n_\circ\,$ a pair of
polar coordinates
\begin{equation}
   (r,\omega )\,\in \mathbb R_+ \times \mathbb S^{n-1} \;\sim \;
    \mathbb R^n_\circ
\end{equation}
where $r = |x|\,$ is referred to as the radial distance and
$\,\omega = \frac{x}{|x|} \,$ as the
 spherical coordinate of $\,x\,$. Obviously $\,x =r\,  \omega\,$ and the volume element in
 polar coordinates reads as $\textrm dx\,=\, r^{n-1}\, \d r\,  \d \omega$
\section{Spherical coordinates, latitude and longitude}\label{seclalo}
A position of a point $\omega \in \mathbb S^{n-1}$ is determined
 by its latitude and longitude. These are convenient coordinates for describing deformations of the unit sphere.
We refer to the points ${\omega}_+=(0,...,0,1)$
 and ${\omega}_-=(0,...,0,-1)$ as north and south poles, respectively. The unit sphere with poles removed will be
denoted by ${\mathbb S}^{n-1}_\pm$. The equatorial sphere
$\,{\mathbb S}^{n-2} \subset {\mathbb S}^{n-1}$ is given by the
equations $x_1^2+...+x_{n-1}^2=1$ and $x_n=0$. To every point
$\,\omega \in {\mathbb S}^{n-1}_\pm$ there correspond meridian
coordinates
\begin{equation}
   (\theta ,{\mathfrak s} )\,\in   \left(0, \pi\right) \times \mathbb S^{n-2}  \;\sim \;
    \mathbb S^{n-1}_\pm
\end{equation}
Here ${\mathfrak s}\,$
lies in the equatorial sphere while $\theta$
is the distance  south of the north pole measured in degrees. The
rectangular coordinates of $\,\omega = (x_1, x_2, ..., x_n) $ are
uniquely recovered from ${\mathfrak s}$ and $\theta$ by the rule
\begin{equation}
\omega = (\cos \theta\, , \, {\mathfrak s}\cdot\sin \theta)\, , \;
\; \; \; \; \; \; 0< \theta< \pi
\end{equation}
  The $(n-1)$-surface area of ${\mathbb S}^{n-1}$ is expressed in
  terms of the
meridian angle as
\begin{equation}
\textrm d{\omega} \;=\;(\sin\theta)^{n-2}\,\, \textrm d\theta\,
\textrm d {\mathfrak s}\;,\;\;\;\;\;\;\;\int_{\mathbb S^{n-1}}
\textrm d \omega \,=\,\omega_{n-1}.
\end{equation}
where $\textrm d {\mathfrak s}$ stands for the $(n-2)$-surface area
of the equatorial sphere. Therefore,
\begin{equation}
  \int_0^\pi \sin^{n-2}\theta \;\textrm d \theta = \frac{\omega_{n-1}}{\omega_{n-2}}
\end{equation} We turn now to the basic examples of mappings
between annuli.

\section{Radial stretching}
These are the mappings which change the radial distance but leave
the spherical coordinate fixed. Let $\mathscr A\mathcal{C}\,
[r,R],\; 0\leqslant r \leqslant R<\infty$, denote the space of
absolutely continuous functions  on a closed interval $[r,R]$.
Associated to every $\,H\in\mathcal A\mathcal{C}\, [r,R]$ is a
mapping $h:\mathbb A \rightarrow \mathbb A^\ast$ of the annulus
${\mathbb A}= \left\{x\; ; \; \; \; r \leqslant |x| \leqslant R
\right\}$ \;
  onto an annulus  $\; {\mathbb A^{\! \ast} }= \left\{y\; ; \; \; \; r_{\! \ast} \leqslant |y| \leqslant
  R_{ \ast}
 \right\}$, defined by the rule
\begin{equation}
 h(x)  =  H(r)\,\omega \, , \;\; \;\;\textrm{ where }\;  r=|x| \; \mbox{ and }\textrm{\footnotesize{$\omega=\frac{x}{|x|}$}} \in\; \mathbb S^{n-1}
\end{equation}
Here we have
\begin{equation}
r_{\! \ast} = \underset{r\leqslant t\leqslant R}{\min}
 \;|H(t)| \;\leqslant \; \underset{r\leqslant t\leqslant R}{\max}
 \;|H(t)| = R_{ \ast}
\end{equation}
We call $h$ a radial stretching and refer to $ H : [r,R] \rightarrow
[r_\ast, R_\ast] $ as its normal  strain. Note that a composition of
two such mappings is again a radial stretching; absolute continuity
is not lost.

\section{Spherical mappings}
A generalization of radial stretching is the so-called
{\it spherical mapping}

\begin{equation}
 h(x)  =  H(|x|)\,\Phi \textrm{\large(\scriptsize{$\frac{x}{|x|}$}\large)}
 \;\;,\;\;\;\;\textrm{where}\;\;\;\;\;\Phi : \mathbb S^{n-1} \rightarrow  \mathbb S^{n-1}
\end{equation}
Here, as in the case of radial stretchings, the normal strain
function is assumed to be absolutely continuous. The spherical part
$\Phi : \mathbb S^{n-1} \rightarrow \mathbb S^{n-1}$, called {\it tangential tension,} is continuous and weakly differentiable. In
what follows we will be actually concerned with homeomorphisms $\Phi
: \mathbb S^{n-1} \rightarrow \mathbb S^{n-1}$ in the Sobolev class
$\mathscr W^{1,n}(\mathbb S^{n-1}, \mathbb S^{n-1})$. Composition of
two spherical mappings results in the  composition of   their normal strain
functions and tangential tensions, respectively.
\chapter{Radial $n$-harmonics}\label{secrad}
In this chapter  we identify all solutions to the $n$-harmonic equation
\begin{equation}
\mbox{div}\, \norm Dh \norm^{n-2}Dh =0\; \mbox{ of the form }\;
h(x)=H(|x|)\, \frac{x}{|x|}
\end{equation}
Such solutions, called {\it radial $n$-harmonics}, will be
originally defined in the annulus
\begin{equation}
{\mathbb A}={\mathbb A}(a,b) = \left\{x\, ; \; \; \; a< |x| < b
\right\}\, , \; \; \; \; \; 0 < a < b < \infty
\end{equation}
Concerning regularity, we assume that the {\it strain function} $H$
belongs to the Sobolev space
\begin{equation}
{\mathscr W}^{1,n}_{\loc} (a,b) \subset {\mathscr C}^{\alpha}_{\loc}
(a,b)\, , \; \; \; \; \; \alpha=\frac{n-1}{n}
\end{equation}
We are going to  give a clear account of how to solve the $n$-harmonic
equation for the radial $n$-harmonics. But before, let us state in advance that the radial $n$-harmonics will actually
extend as $n$-harmonics to the entire punctured space $\rn_\circ$.
Even more, they will extend continuously to the M\"obius space
$\hat{{\mathbb R}}^n =\rn \cup \{\infty \}$, one point compactification of $\R^n$.

Let ${\mathcal H}_n$ denote the class of all $n$-harmonics in
$\hat{{\mathbb R}}^n$, the following elementary transformations of variables $x$ and $h$
preserve this class.
\begin{itemize}
\item{{\it Rescaling}; $$h(x)\in {\mathcal H}_n  \mbox{ implies }  \lambda \, h(k x)\in {\mathcal H}_n \mbox{ for every } k >0 \mbox{ and } \lambda \in {\mathbb R}$$}
\item{{\it Reflection};  $$h(x)\in {\mathcal H}_n \mbox{ implies }  h\left(x|x|^{-2}\right)\in {\mathcal H}_n$$}
\end{itemize}
The radial $n$-harmonics in ${\mathbb A} = \left\{x\, ; \; \; a< |x|
< b \right\}$, $0<a<b<\infty$ are none other than the unique solutions of
the Dirichlet problem
\begin{equation}\label{E21}
\begin{cases}\mbox{div}\, \norm Dh\norm^{n-2}Dh=0 \; \; & \mbox{ for } a<|x|<b \\
h(x)=\alpha x &  \mbox{ for } |x|=a\\
h(x)=\beta x &  \mbox{ for } |x|=b
\end{cases}
\end{equation}
where $\alpha$ and $\beta$ can be any real numbers, so $H(a)= \alpha
\, a$ and $H(b)= \beta \, b$. It is understood here that $h$ extends
continuously to the closed annulus ${\mathbb A}\, [a,b]$.  We shall
distinguish four so-called {\it principal $n$-harmonics} in $\hat{{\mathbb
R}}^n$ and use them to generate all radial $n$-harmonics via
rescaling of the variables $x$ and $h$. For this we need  some
computation.
\section{The $n$-Laplacian for the strain function}
Suppose $H \in  {\mathscr C}^{2} (a,b)$, $0 \leqslant a < b
\leqslant \infty$. Thus the derivatives of the radial stretching
\begin{equation}
h(x)=H(|x|)\, \frac{x}{|x|}
\end{equation}
exist up to order 2. The differential matrix of $h$ can be computed as
\begin{equation}
Dh(x)=\frac{H(|x|)}{|x|}{\bf I} + \frac{|x|\dot{H}(|x|)-H(|x|)}{|x|}\cdot \frac{x\otimes x}{|x|^2}
\end{equation}
Hereafter $\dot{H}$ stands for the derivative of $H=H(t)$, $a<t<b$.
Hence the Hilbert-Schmidt norm of the differential matrix is expressed by
means of $H$ as
\begin{equation}
\norm Dh(x) \norm^2 = \mbox{Tr}\left( D^\ast  h\cdot  Dh\right)=
\dot{H}^2+ (n-1)\frac{H^2}{|x|^2}
\end{equation}
A lengthy  computation leads to the following explicit formula for
the $n$-Laplacian
\begin{eqnarray}\label{E211}
 \mbox{div}\, \norm Dh\norm^{n-2}Dh &=& \frac{(n-1)|Dh|^{n-4}}{|x|^4}\cdot \bigg\{t^2 (t^2\dot{H}^2+H^2)\ddot{H} +\\
+ t^3 \dot{H}^3&+&(n-3)t^2 \dot{H}^2 H + t \dot{H}H^2
-(n-1)H^3 \bigg\} \frac{x}{|x|}\nonumber
\end{eqnarray}
where $t=\abs{x}$. Our subsequent analysis will be based on  the {\it
elasticity function} already mentioned at (\ref{IE17}),
\begin{equation}
\eta=\eta_H=\frac{t\dot{H}(t)}{H(t)}.
\end{equation}
Let us make certain here that we shall never deal with   the case in which
both $\dot{H}$ and $H$ vanish simultaneously, so $\eta$ will be well
defined even when $H$ vanishes. Now $h$ is $n$-harmonic if and only
if
\begin{equation}\label{243}
\left(1+\eta^2 \right)t^2 \ddot{H}= \left(1-\eta \right)
\left[\eta^2 +(n-2)\eta +n-1\right]\cdot H
\end{equation}
where the case $\eta = \infty$ should read as $H=0$. In dimension
$n=2$ the $n$-harmonic equation reduces to the familiar linear
Cauchy-Euler equation for $H$
\begin{equation}\label{E24}
t^2 \ddot{H}+ t\dot{H}-H =0
\end{equation}
It has two linearly independent solutions, denoted by
\begin{equation}\label{E216}
H_\circ (t)=t \; \;  \; \mbox{ and } \; \; \; H_\infty
(t)=\frac{1}{t}
\end{equation}
They generate all solutions by means of linear
combinations. We  pick out two
particular solutions
\begin{equation}
H_+ (t)=\frac{1}{2}\left(t+\frac{1}{t}\right) \; \;  \; \mbox{ and } \; \; \; H_- (t)=\frac{1}{2}\left(t-\frac{1}{t}\right)
\end{equation}
and refer to $H_\circ, H_\infty, H_+$ and $H_-$ as principal
solutions. The reason for preferring these principal solutions to the ones in
(\ref{E216}) is simply that we will be able to generate all
radial $n$-harmonics, simply  be rescaling the principal ones. Such a set of principal solutions will be
particularly useful in higher dimensions where the $n$-Laplacian is
nonlinear. For $n=2$ the corresponding complex  principal harmonics are
\begin{equation}
h_\circ (z)=z\, , \; \;  \;  \; \; \; h_\infty (z)=\frac{1}{z}\, ,
\end{equation}
\begin{equation}
h_+ (z)=\frac{1}{2}\left(z+\frac{1}{\bar{z}}\right) \; \;  \; \mbox{ and } \; \; \; h_- (z)=\frac{1}{2}\left(z-\frac{1}{\bar{z}}\right)
\end{equation}
The reader may wish to observe that $h_\circ$ and $h_\infty$ map
$\hat{{\mathbb R}}^2$ univalently onto itself, whereas $h_+ , h_- :
\hat{{\mathbb R}}^2 \rightarrow \hat{{\mathbb R}}^2$ have a branch
set (folding) along the unit circle. Precisely, $h_+$  covers the exterior of
the unit disk twice while $h_-$ is a double cover of the entire
space $\hat{{\mathbb R}}^2$. In spite of nonlinearity the analogous
ordinary differential equations in higher dimensions are still effectively solvable, with some
computational efforts. The key is the following identity which follows from
(\ref{E211}).
\begin{lemma}\label{L21}
Let $h(x)=H(|x|)\, \frac{x}{|x|}$ be a radial stretching of class $
{\mathscr C}^{2} ({\mathbb A}, \rn)$, where ${\mathbb A} =
\left\{x\, ; \; \; \; a<|x| < b \right\}$. Then, with the notation $\,|x| = t\,$, we have
\begin{eqnarray}\label{Ediv}
&\, & \frac{n \,\dot{H}}{n-1} \, \textnormal{div}\, \norm Dh\norm^{n-2}Dh \nonumber \; =\; \\
&\, & - \;\frac{x}{t^{n+1}} \frac{d}{dt}\Big \{\left[(n-1)H^2 + t^2
\dot{H}^2\right]^\frac{n-2}{2}\left(H^2-t^2\dot{H}^2\right)\Big\}\,
\end{eqnarray}
\end{lemma}
\section{The principal solutions}
The main idea is to eliminate  one derivative in the second order $n$-harmonic
equation. This will lead us to a problem of solving a first order
ODE with one free parameter. For nonlinear equations such a
procedure is rather tricky. Fortunately, there is a very
satisfactory realization of this idea due to the identity
(\ref{Ediv}). We wish to express the $n$-harmonic equation  for
$h(x)=H(|x|)\, \frac{x}{|x|}$ in the form
\begin{equation}
\mathcal{L}H \bydef F(t, H, \dot{H})= \mbox{\it const.}
\end{equation}
By virtue of Lemma \ref{L21}, the nonlinear differential operator
${\mathcal L}$ takes the following explicit form
\begin{equation}\label{E222}
{\mathcal L}H = \left[H^2 + \frac{t^2
\dot{H}^2}{n-1}\right]^\frac{n-2}{2} \left(H^2-t^2\dot{H}^2\right)= \mbox{{\it const.}}\;
\; \; \footnote{In general, guessing the identity such as
(\ref{E222}) takes some efforts; though verifying it poses no
challenge.}
\end{equation}
Certainly, all radial $n$-harmonics can be found by solving the
so-called {\it characteristic equation}
\begin{equation}
{\mathcal L}H = \mbox{\it const.}
\end{equation}
However, for the converse statement, caution must be exercised
because this equation only implies that $\dot{H} \, \mbox{div}\,
\norm  Dh \norm^{n-2} Dh =0$. For example, the constant function
$H(t)\equiv 1$ satisfies ${\mathcal L}H\equiv 1$, but the
corresponding radial mapping $h(x)=\frac{x}{|x|}$ is not
$n$-harmonic. \footnote{This noninjective solution will  play an
important role in the sequel.} More generally, there are many more
solutions
 to the equation ${\mathcal L}H\equiv
1$ which are constant  on a subinterval. For example,  one can paste $\frac{x}{\abs{x}}$ with a suitable $n$-harmonic map on the rest of the interval. Such solutions, however, cannot be $\mathscr C^2$-smooth;
${\mathscr C}^{1,1}$-regular at best. They, nevertheless solve the generalized $n$-harmonic equation~\eqref{enhe}.

\begin{figure}[!h]
\begin{center}
\psfrag{2}{\tiny${ 2}$}
\psfrag{4}{\tiny${ 4}$}
\psfrag{1}{\tiny${ 1}$}
\psfrag{H}{${\tiny H }$}
\psfrag{f1}{${\tiny {\, }_{2\left(t+\frac{1}{t}\right)}}$}
\psfrag{f2}{${\tiny {\, }_{t+\frac{4}{t} }}$}
\psfrag{t}{${\tiny t}$}
\includegraphics*[height=1.6in]{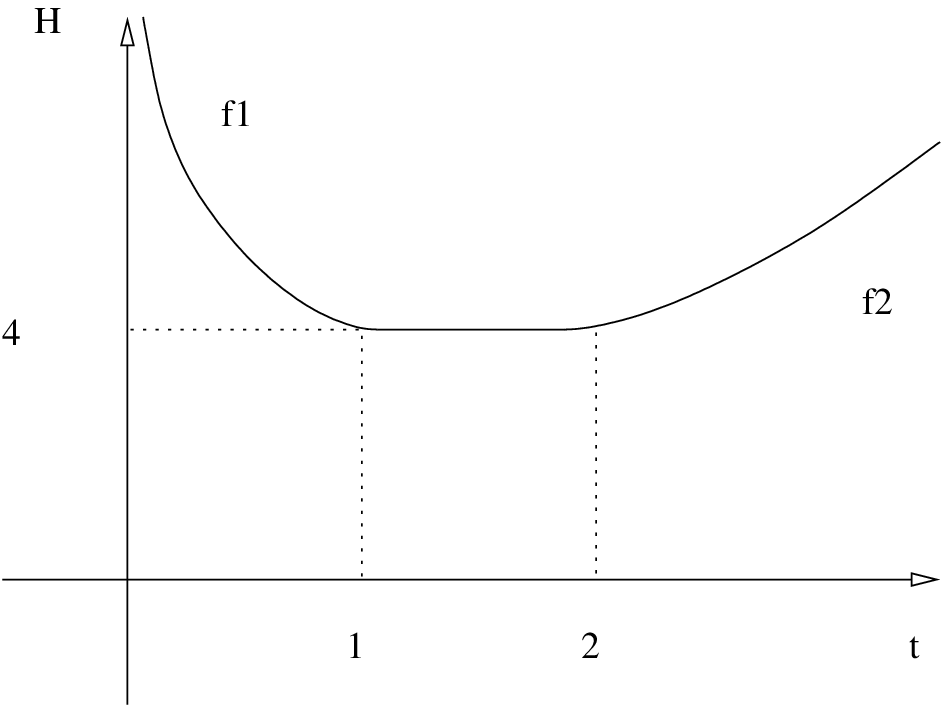}
\caption{A solution which is constant on an interval.}
\end{center}
\end{figure}

Solutions $H=H(t)$ to the characteristic equation   ${\mathcal L}H \equiv c $, that are constant in a proper subinterval like in the above figure, cannot be even ${\mathscr C}^{1,1}$-regular if  $c \leqslant 0$.

Let us summarize our findings as:
\begin{definition}
The term principal solution pertains to each of the following four
functions of class  ${\mathscr C}^2 (0, \infty)$ which solve the
equation ${\mathcal L}H= \mbox{\it constant}$;
\begin{equation}
 H_\circ (t)= t\, , \; \; \; \; \; \; \; \;   {\mathcal L}H_\circ \equiv 0
\end{equation}
\begin{equation}
 H_\infty (t)= \frac{1}{t}\, , \; \; \; \; \; \; \; \;   {\mathcal L}H_\infty \equiv 0
\end{equation}
\begin{equation}
{\mathcal L}H_+ \equiv 1\, , \; \; \; \; \; \; \; \; H_+ (1)= 1
\end{equation}
\begin{equation}
{\mathcal L}H_-\equiv -1\, , \; \; \; \; \; \; \; \; H_- (1)= 0
\end{equation}
\end{definition}
We shall solve these  Cauchy problems and show that both
$H_+$ and $H_-$ are actually ${\mathscr C}^\infty$-smooth.

\begin{figure}[!h]
\begin{center}
\psfrag{1}{\tiny ${ 1}$}
\psfrag{H0}{\small${ H_\circ = t}$}
\psfrag{Hinf}{\small${ H_\infty = \frac{1}{t}}$}
\psfrag{H+}{\small ${ H_+ = \frac{1}{2}\left(t+\frac{1}{t}\right)}$}
\psfrag{H-}{\small ${ H_- = \frac{1}{2}\left(t-\frac{1}{t}\right)}$}
\includegraphics*[height=3.0in]{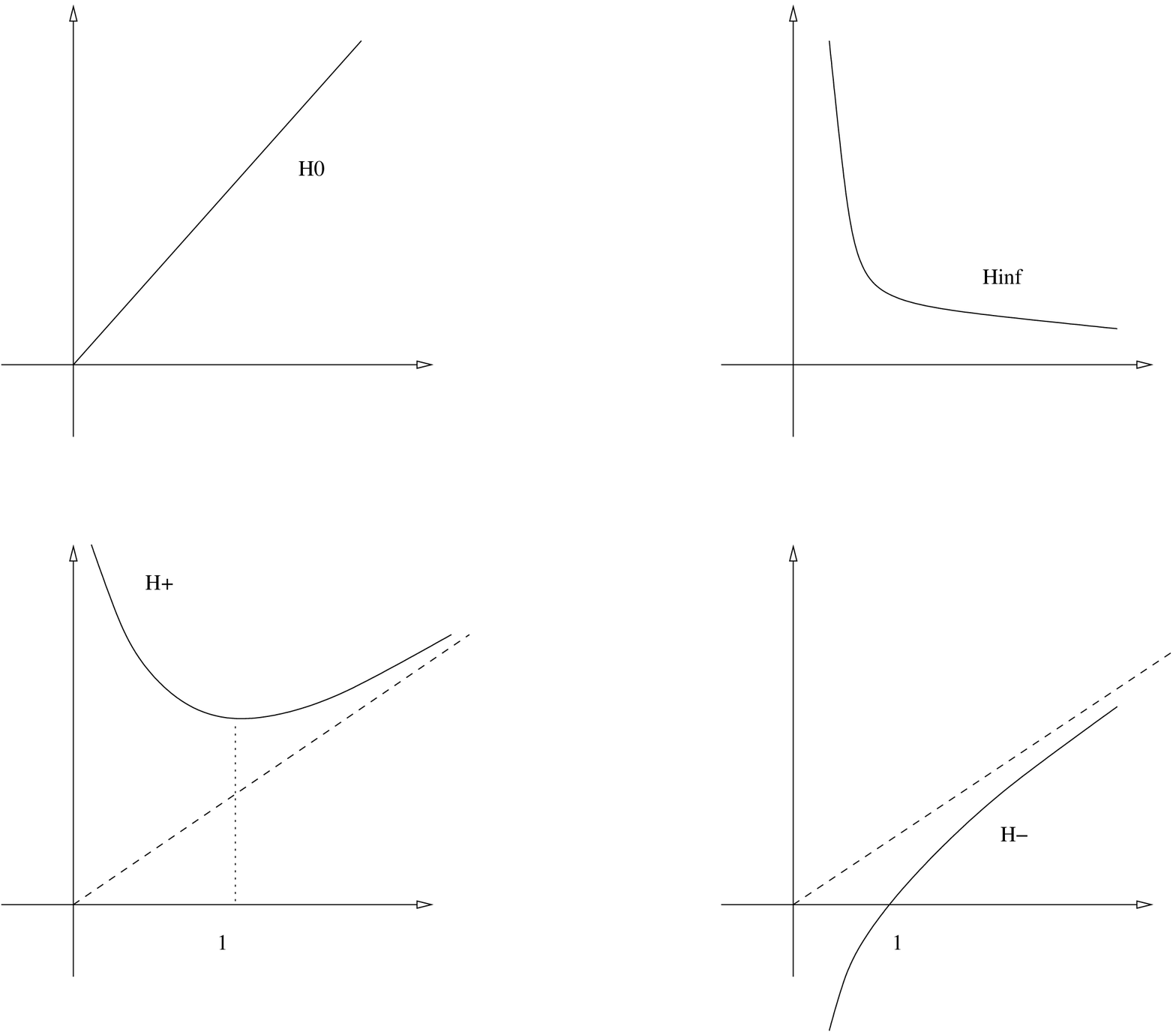}
\caption{The principal solutions for $n=2$.  Similar ones  for $n
\geqslant 3$ are illustrated in Figures~\ref{Fig10} and~\ref{figur10}.}\label{Fig7}
\end{center}
\end{figure}

The name ``principal solutions'' is given to $H_\circ $, $H_\infty
$,  $H_+$ and  $H_- $ because they generate all radial
$n$-harmonics. We denote the corresponding radial mapping by $h_\circ$,
$h_\infty$, $h_+$ and $h_-$, and call them {\it principal $n$-harmonics}.
Precisely, we have
\begin{proposition}
Every radial $n$-harmonic mapping in the annulus ${\mathbb A}= \{x\,
; \; \; \; a < |x|<b\}$ takes the form
\begin{equation}
g(x)=\lambda \cdot h(k x)\, , \; \; \; \; \; \; \; \; \; \; \; \;
\lambda \in {\mathbb R}\, , \; \; k >0
\end{equation}
where $h \in {\mathscr C}^\infty (\rn_{\, \circ})$,  is one of the
four principal $n$-harmonics.
\end{proposition}
It should be observed, as a corollary, that radial $n$-harmonics are
${\mathscr C}^\infty$-smooth in the entire space $\rn_\circ$.
\section{The elasticity function}
We shall distinguish four classes of the radial $n$-harmonics
$h(x)=H(|x|)\frac{x}{|x|}$. The concept of so-called conformal
elasticity underlines this distinction. In what follows it will
never be the case that both $H$ and $\dot{H}$ vanish simultaneously.
Consider any function $H\in {\mathscr C}^1(r,R) $, $0<
r<R < \infty$, such that $H^2 + \dot{H}^2>0$. The elasticity
function
\begin{equation}\label{E232}
\eta (t)=\eta_{_H} (t) = \frac{t \dot{H}(t)}{H(t)}
\end{equation}
is  continuous  on $(r,R)$. In general, it assumes values in the
extended real line $\hat{\mathbb R} = {\mathbb R} \cup \{-\infty ,
\infty\}$. We recall that  $h=H\big(\abs{x}\big) \frac{x}{\abs{x}}$ is:
\begin{itemize}
\item{{\it Conformally contracting}, if $|\eta (t)| < 1$}
\item{{\it Conformally expanding}, if $|\eta (t)| > 1$}
\item{{\it Conformal}, if $|\eta (t)| =1$}
\end{itemize}
Observe that the equation $\eta (t)\equiv 1$ yields
$H(t)=\lambda \cdot t$, while $\eta (t)=-1$ gives
$H(t)=\frac{\lambda}{t}$. These are special cases of the power
functions $H(t)=\lambda\, t^\alpha$, having constant elasticity
$\eta(t) \equiv \alpha$. The above three classes  are invariant
under rescaling and inversion. Precisely, we have the
formulas:
\begin{equation}\label{E21}
\begin{cases}\eta_{_H}(kt) = \eta_{_F}(t)  &  \mbox{ where } F(t)=H(kt) \\
\eta_{_H}(1/t) =-\eta_{_F}(t) \; \;  &  \mbox{ where } F(t)=H(1/t)
\end{cases}
\end{equation}
The  elasticity function tells us something about the infinitesimal relative rate of change of the modulus of an annulus under the deformation $h=H\big(\abs{x}\big) \frac{x}{\abs{x}}$.  Precisely, we have
for $r<t<R$
\begin{equation}
\eta_{_H}(t)= \lim\limits_{\epsilon \rightarrow 0} \,
\frac{\mbox{Mod}\, A^\ast(t+\epsilon )- \mbox{Mod}\,
A^\ast(t)}{\mbox{Mod}\, A(t+\epsilon )- \mbox{Mod}\, A(t)}
\end{equation}
where $A(t)={\mathbb A}\, (r,t)  \mbox{ and }
A^\ast(t)=h[A(t)]$.

 We shall be concerned with ${\mathscr
C}^1$-solutions to the characteristic equation
\begin{equation}\label{E231}
{\mathcal L}H=\left[H^2 + \frac{t^2 \dot{H}^2}{n-1}\right]^\frac{n-2}{2} \left(H^2-t^2\dot{H}^2\right)\equiv c
\end{equation}
This equation is invariant under rescaling
and inversion of the independent variable $t$. Accordingly, all
${\mathscr C}^1$-solutions  fall into three
categories:
\begin{itemize}
\item{$H$ is contracting if $c>0$, equivalently  $|\eta (t)|<1$}
\item{$H$ is expanding if $c<0$, equivalently  $|\eta (t)|>1$}
\item{$H$ is conformal if $c=0$, equivalently  $|\eta (t)|=1$,\\ Furthermore, in the conformal case   $$\left\{\begin{array}{ll} \mbox{$h$ is preserving the order of boundary components if } & \eta (t)\equiv 1\, ,\;  H(t)= \lambda t\\
  \mbox{$h$ is  reversing the order of boundary components if } & \eta (t)\equiv -1, H(t)= \frac{\lambda}{t} \end{array} \right.$$   }
\end{itemize}
This is easily seen by writing (\ref{E231}) as
\begin{equation}\label{equ123}
\left( 1 + \frac{\eta^2}{n-1}\right)^\frac{n-2}{2} \left(1-\eta^2\right)=\frac{c}{|H|^n}
\end{equation}
We then distinguish four classes of the radial $n$-harmonics:
\begin{equation}
{\mathcal H}_n = {\mathcal H}_+ \cup {\mathcal H}_1 \cup {\mathcal H}_\circ \cup {\mathcal H}_\infty
\end{equation}
where
\begin{eqnarray*}
{\mathcal H}_+ &=& \Big\{ h=H(|x|)\, \frac{x}{|x|}\, ; \; \; \; \; {\mathcal L}H\equiv c >0 \Big\} \\
{\mathcal H}_- &=& \Big\{ h=H(|x|)\, \frac{x}{|x|}\, ; \; \; \; \; {\mathcal L}H\equiv c <0 \Big\} \\
{\mathcal H}_\circ &=& \Big\{ h \, ; \; \; \; \; h(x) = \lambda x \, , \; \; \; \lambda \in {\mathbb R}_\circ \Big\} \\
{\mathcal H}_\infty &=& \Big\{ h \, ; \; \; \; \; h(x) =
\frac{\lambda \cdot x}{|x|^2}\, , \; \; \; \lambda \in {\mathbb
R}_\circ \Big\}
\end{eqnarray*}
We strongly emphasize again that not every solution to the equation ${\mathcal L}H \equiv
c$ defines a radial $n$-harmonic map, due to lack of regularity. There are Lipschitz solutions to ${\mathcal L}H=0$ lacking ${\mathscr C}^1$-regularity

\begin{figure}[!h]
\begin{center}
\psfrag{1}{\tiny ${ 1}$}
\psfrag{t}{$t$}
\psfrag{1/t}{$\frac{1}{t}$}
\includegraphics*[height=1.4in]{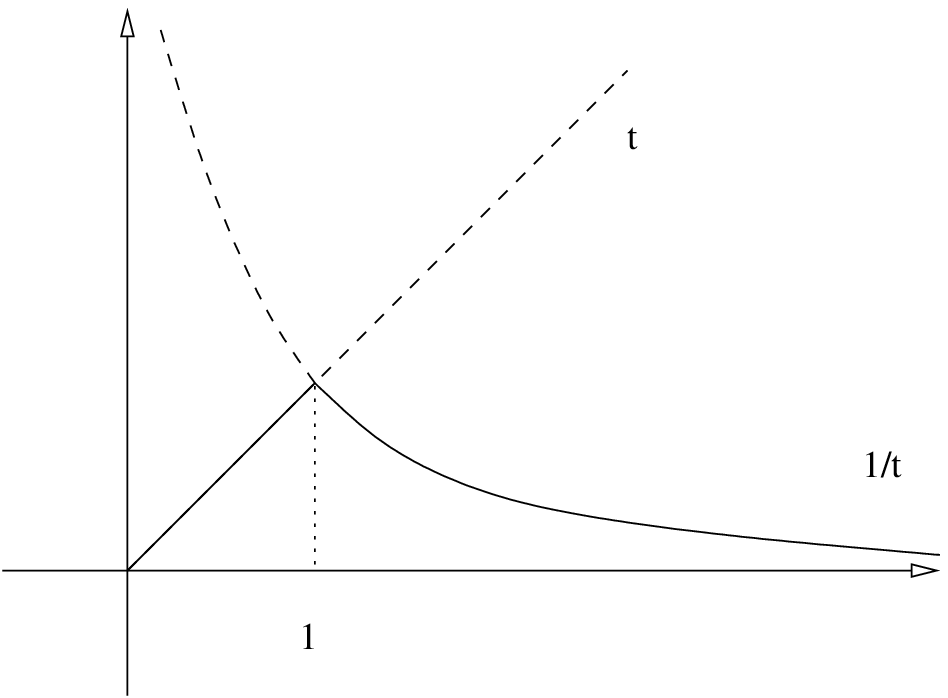}
\caption{Lipschitz solution to ${\mathcal L}H=0$.}
\end{center}
\end{figure}

Such solutions are not particularly desirable. Even ${\mathscr
C}^1$-solutions to the equations ${\mathcal L}H=1$ need not give
radial $n$-harmonics. For example, if we paste $H_+(t)=\frac{1}{2}
\left(t + \frac{1}{t}\right)$ $1\leqslant t < \infty$ with function identically equal to 1 for $ 0<t
\leqslant 1$ then the resulting function will become a ${\mathscr
C}^{1,1}$-solution to the equation ${\mathcal L}H =1$ on $(0,
\infty)$.

\begin{figure}[!h]
\begin{center}
\psfrag{1}{\tiny ${ 1}$}
\includegraphics*[height=1.4in]{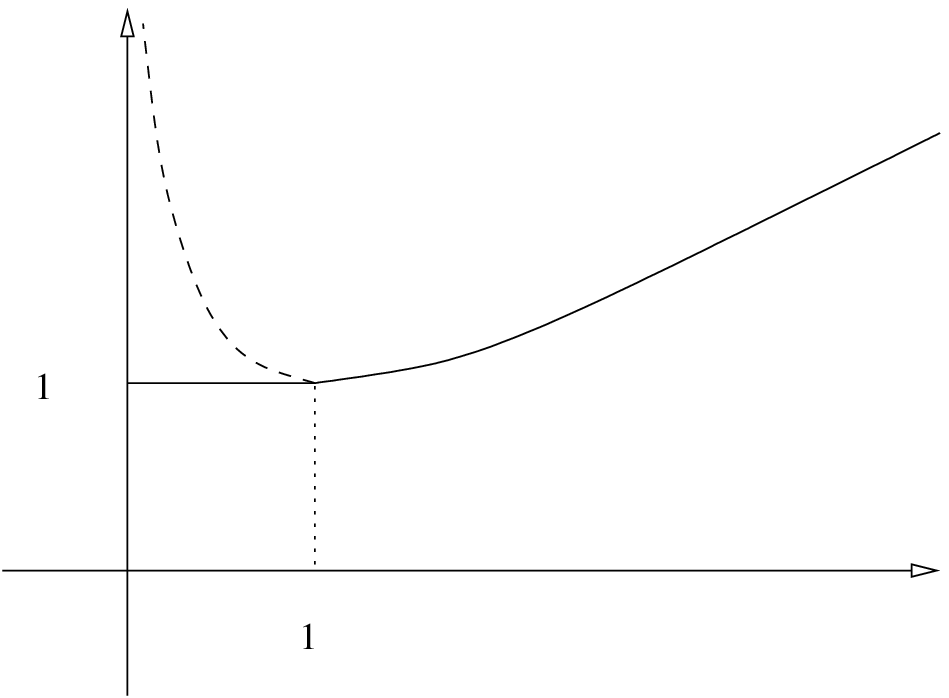}
\caption{A ${\mathscr C}^1$-solution to ${\mathcal L}H=1$.}
\end{center}
\end{figure}

The radial stretching obtained in this way  fails  to be
$n$-harmonic on the punctured disk  ${\mathbb A}\, (0,1)$. Thus we
see various degenerations of ${\mathscr C}^{1}$-solutions,
suggesting that we should restrict ourselves to the solutions of
class ${\mathscr C}^{2}(0, \infty )$.\footnote{In the expanding
case ${\mathscr C}^{1}$ solutions are always smooth.}

The two principal solutions $H_+$ and $H_-$ are particularly  interesting and
important to the forthcoming results, so we devote next two
sections for a detailed treatment.
\section{The principal solution ${H}_+$ (conformal contraction)}
We find  ${H}_+$  as a  ${\mathscr
C}^2$-solution to the Cauchy problem
\begin{equation}\label{241}
{\mathcal L}H = \left[H^2 + \frac{t^2 \dot{H}^2}{n-1}\right]^\frac{n-2}{2} \left(H^2-t^2\dot{H}^2\right)\equiv 1\, , \; \; \; \; H(1)=1
\end{equation}
The general idea behind our method of solving this equation is as
follows. Differentiating   (\ref{241}) yields a first order equation
for the elasticity function,
\begin{equation}\label{242}
\frac{(1+\eta^2)\, \dot{\eta}}{\left(1-\eta^2\right)\left(\eta^2+n-1\right)}=\frac{1}{t}\, , \hskip2cm |\eta (t)| <1
\end{equation}
This can also be directly derived from (\ref{E211}). With the
equation (\ref{242}) at hand  we now proceed to the explicit
computation. Consider the function $\Gamma_+ = \Gamma_+ (s)$ defined for $-1<s<1$ by the rule
\begin{eqnarray}\label{244}
\begin{split}
\Gamma_+ (s) & = \mbox{\large{exp}}\int_0^s \frac{(1+\tau^2)\,
d\tau}{(1-\tau^2) (\tau^2+n-1)}\\ & = \mbox{$\sqrt[n]{{\tiny
\frac{1+s}{1-s}}}$} \,  \mbox{\large{exp}}\left[\mbox{${\tiny
\frac{2-n}{n\sqrt{n-1}}}$}\tan^{-1}\mbox{${\tiny
\frac{s}{\sqrt{n-1}}}$}\right]
\end{split}
\end{eqnarray}
Obviously $\Gamma_+$ is strictly increasing from zero to infinity,
thus invertible.

\begin{figure}[!h]
\begin{center}
\psfrag{1}{\tiny ${ 1}$} \psfrag{-1}{\tiny ${ -1}$}
\psfrag{s}{\small${ s}$} \psfrag{Gam}{\small ${ \Gamma_+}$}
\psfrag{Fir}{$\,$} \psfrag{Sec}{$\,$} \psfrag{Thi}{$\,$}
\includegraphics*[height=1.6in]{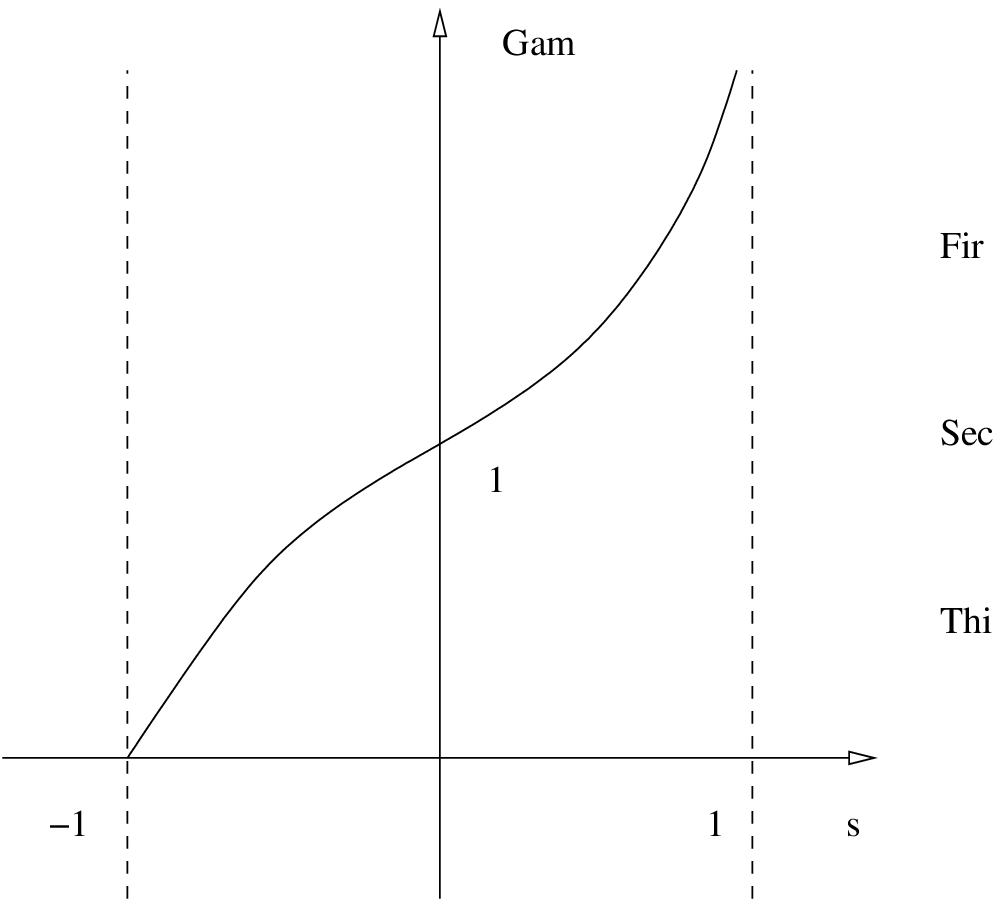}
\caption{$\Gamma_+$ is strictly increasing and ${\mathscr
C}^\infty$-smooth.}
\end{center}
\end{figure}

Note the identities
\begin{equation}
 \begin{array}{ccc} & \Gamma_+(0)=1  \\
  & \dot{\Gamma}_+ (0)=\frac{1}{n-1} \\
& \Gamma_+ (-s)\, \Gamma_+ (s)=1
\end{array}
\end{equation}
We then examine the inverse function defined for $0<t< \infty$,
\begin{equation}
u(t)=u_+(t)=\Gamma_+^{-1}(t)\, , \; \; u(1)=0\, , \; \;
\dot{u}(1)=n-1
\end{equation}
The identity $\Gamma_+ (-s)\, \Gamma_+ (s)\equiv 1$ translates  as the antisymmetry rule for
$u$
\begin{equation}\label{246}
u \left(\frac{1}{t}\right)= - u(t)\, , \; \; \; \; \; \; \; \; \;
t>0
\end{equation}
Implicit differentiation of the equation $\Gamma_+
\big(u(t)\big)\equiv t$ yields
\begin{equation}
\dot{u}\, \dot{\Gamma}_+ (u)=1\, , \hskip2cm t\,
\dot{u}(t)=\frac{\Gamma_+ (u)}{\dot{\Gamma}_+(u)}
\end{equation}
We shall now  introduce $H_+=H_+(t)$, $0<t< \infty$, and prove that
$H_+$ is the principal solution, namely
\begin{equation}\label{245}
H_+(t)= \left(1+\frac{u^2}{n-1}\right)^{\frac{1}{n}-\frac{1}{2}} \left(1-u^2\right)^{-\frac{1}{n}}\geqslant 1
\end{equation}
From this formula it  follows  that $u$ is the elasticity function
of $H_+$. Indeed,
\begin{equation}
\eta=\eta_+(t)=\frac{t\dot{H}_+}{H_+}= \frac{u\cdot t \, \dot{u}\,
(1+u^2)}{(1-u^2)(n-1+u^2)}=u\cdot t \dot{u}\,
\frac{\dot{\Gamma}_+(u)}{\Gamma_+ (u)}=u(t)
\end{equation}

\begin{figure}[!h]
\begin{center}
\psfrag{1}{\tiny ${ 1}$}
\psfrag{-1}{\tiny ${ -1}$}
\psfrag{t}{${ t}$}
\psfrag{eta}{${\tiny \eta}$}
\psfrag{eta2}{${\tiny \eta=\eta_{_+}(t)}$}
\includegraphics*[height=1.4in]{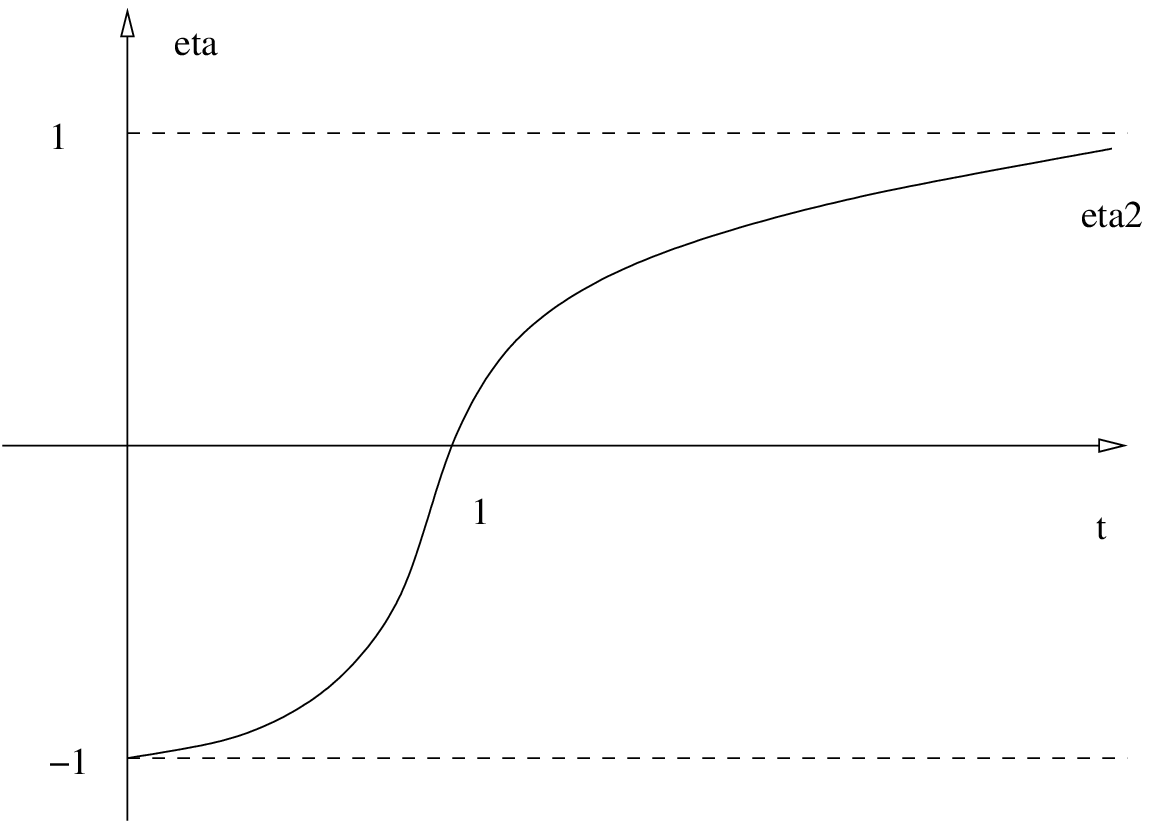}
\caption{The elasticity function of $H_+$.}
\end{center}
\end{figure}

\begin{figure}[!h]
\begin{center}
\psfrag{1}{\tiny ${ 1}$}
\psfrag{H}{${\tiny H}$}
\psfrag{t}{${ t}$}
\psfrag{H+}{${\tiny H = H_+(t)}$}
\includegraphics*[height=1.6in]{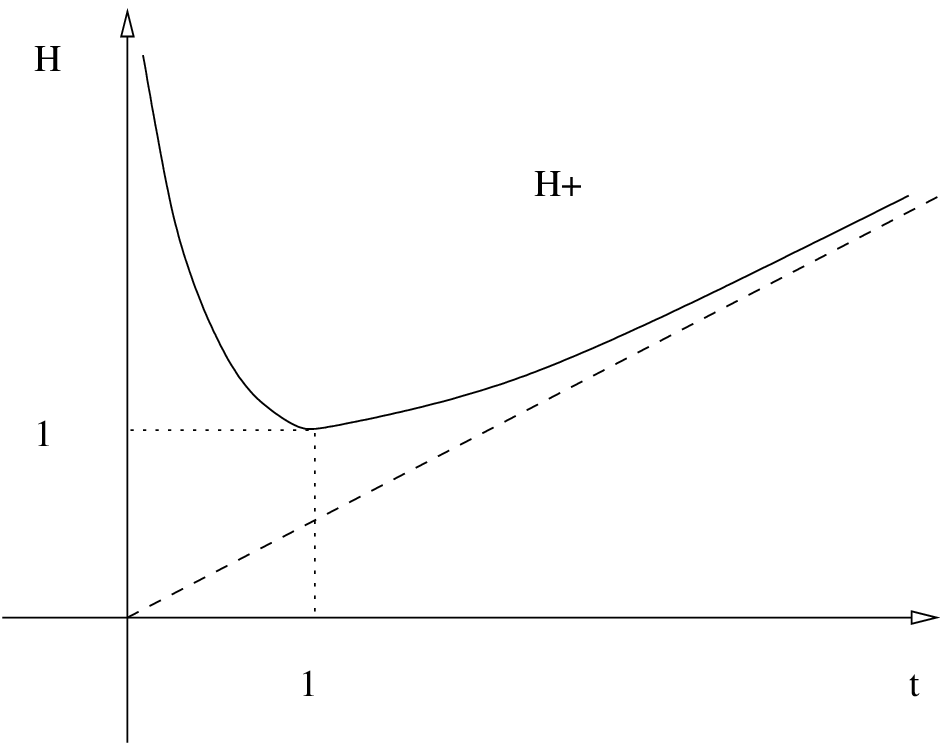}
\caption{The principal solution $H_+$.}\label{Fig10}
\end{center}
\end{figure}

Now, the equation ${\mathcal L}H_+\equiv 1$ can easily be verified

\begin{eqnarray}\label{247}
{\mathcal L}H_+&=& \left|H_+ \right|^n \left[1+\frac{\eta^2}{n-1}\right]^\frac{n-2}{2}\left(1-\eta^2\right)=\nonumber \\
&\, &  \left|H_+ \right|^n \left[1+\frac{u^2}{n-1}\right]^\frac{n-2}{2}\left(1-u^2\right)\equiv 1
\end{eqnarray}
This last step is immediate from (\ref{245}). It is worth noting
that $H_+ \in {\mathscr C}^\infty (0, \infty)$ and $\dot{H}_+(1)=0$.
Antisymmetry rule at (\ref{246}) results in a symmetry rule for
$H_+$
\begin{equation}\label{248}
H_+\left(\frac{1}{t}\right) = H_+\left( t\right)
\end{equation}
Furthermore, it follows from (\ref{243}) that $\ddot{H}_+(t)>0$. Thus $H_+$ is strictly convex, see Figure \ref{Fig10}.

Let us take a look at the  behavior of $H_+$ near zero and
infinity. It has an asymptote at infinity with the slope
\begin{equation}
\Theta = \Theta_+ = \left(1-\frac{1}{n}\right)^\frac{n-2}{2n} 4^{-\frac{1}{n}}\,   \mbox{\large{exp}}\left[\frac{n-2}{n\sqrt{n-1}}\tan^{-1}\frac{1}{\sqrt{n-1}}  \right]
\end{equation}
Indeed, under the notation above,  we have the following equation
for $t=\Gamma (s)$
\[
H^n(t)-\Theta^n\, t^n =
\frac{1}{1-s^2}\left(1+\frac{s^2}{n-1}\right)^{1-\frac{n}{2}}-\Theta^n
\, \Gamma_+^n (s)\bydef \frac{A(s)}{s-1}
\]
Note that $A=A(s)$ is  ${\mathscr C }^\infty$-smooth on $(-1,
\infty)$. It vanishes at $s=1$. Application of L'Hospital rule shows
that
\begin{equation}
\lim\limits_{t\rightarrow \infty} \left[H^n(t)-\Theta^n
t^n\right]=A^\prime (1)>0\, ,
\end{equation}
whence the inequality
\begin{equation}
0<H(t)-\Theta\, t \, \leqslant \,  \frac{C}{t^{n-1}}\, , \; \; \; \;
\; \; C= \frac{A^\prime (1)}{n\Theta^{n-1}} >0
\end{equation}
In particular by (\ref{248}) we conclude that
\begin{equation}\label{E104}
\lim\limits_{t\rightarrow \infty} \frac{H(t)}{t}=\Theta = \lim\limits_{t\rightarrow 0} t \, H(t)
\end{equation}
\section{The principal solution $H_-$ (conformal expansion)}
The principal solution $H=H_- \in {\mathscr C}^\infty (0, \infty)$
is obtained by solving the following Cauchy problem for the
characteristic equation
\begin{equation}\label{251}
{\mathcal L}H = \left[H^2 + \frac{t^2 \dot{H}^2}{n-1}\right]^\frac{n-2}{2} \left(H^2-t^2\dot{H}^2\right)\equiv -1\, , \; \; \; \; H(1)=0
\end{equation}
Obviously $\dot{H}$ cannot vanish. For an explicit computation we introduce
the following function
\begin{eqnarray}\label{252}
 \Gamma_- (s) =
\sqrt{\frac{1+s}{1-s}} \,  \mbox{\large{exp}}\left[ \frac{n-2}{n\sqrt{n-1}}\tan^{-1} \left(s\sqrt{n-1}\right)\right]
\end{eqnarray}
where $-1<s<1$. Its logarithmic derivative is computed as
\begin{equation}\label{253}
\frac{\dot{\Gamma}_-(s)}{\Gamma_-
(s)}=\frac{1+s^2}{(1-s^2)\left[1+(n-1)s^2\right]}>0
\end{equation}
This shows that $\Gamma_- : (-1,1)\rightarrow (0, \infty)$ is
strictly increasing. Consequently,  we  consider its inverse
function
\begin{equation}\label{254}
u(t)=u_-(t)=\Gamma^{-1}(t)\, , \hskip2cm 0<t<\infty
\end{equation}
The graphs of $\Gamma_-$ and $u_-$ are very similar to those of $\Gamma_+$ and $u_+$, though they are not the same. However, this time the same identity $\Gamma (-s)\, \Gamma (s)=1$ yields different transformation rule for $u$,
\begin{equation}\label{255}
u\left(\frac{1}{t}\right)=-u(t)\, , \hskip2cm 0<t<\infty
\end{equation}
It is at this point that  the analogy with $H_+$ ends.  We define the principal solution $H_-$ as
\begin{equation}\label{256}
H(t)=H_-(t)=u\, \left(u^2+\frac{1}{n-1}\right)^{\frac{1}{n}-\frac{1}{2}} \left(1-u^2\right)^{-\frac{1}{n}}
\end{equation}
Note that $H_-$ changes sign at $t=1$, exactly where $u$ vanishes.
Basic  properties of $H_-(t)$ can be derived directly from~(\ref{256}). Let us state some of them. We observe
the first difference  between   $H_+$  and $H_-$ in the antisymmetry formula
\begin{equation}\label{257}
H_-\left(\frac{1}{t}\right)=-H_-(t)
\end{equation}
which is immediate from (\ref{255}). Next for $t\not= 1$ we compute
the elasticity function for $H_-$,
\begin{eqnarray}\label{ES108}
\eta = \eta_{_-}(t)&= & \frac{t\, \dot{H}_-(t)}{H_-(t)} = \frac{(1+u^2)\, t\dot{u}}{u\, (1-u^2)\left[1+(n-1)u^2\right]}\nonumber \\
&=& \frac{t\, \dot{u}}{u}\cdot \frac{\dot{\Gamma}_-(u)}{\Gamma_-
(u)} = \frac{1}{u}
\end{eqnarray}

\begin{figure}[!h]
\begin{center}
\psfrag{1}{\tiny ${1}$}
\psfrag{-1}{\tiny ${ -1}$}
\psfrag{t}{\small ${ t}$}
\psfrag{eta}{\small ${ \eta}$}
\psfrag{eta2}{\small ${ \eta=\eta_{_-}(t)}$}
\includegraphics*[height=1.6in]{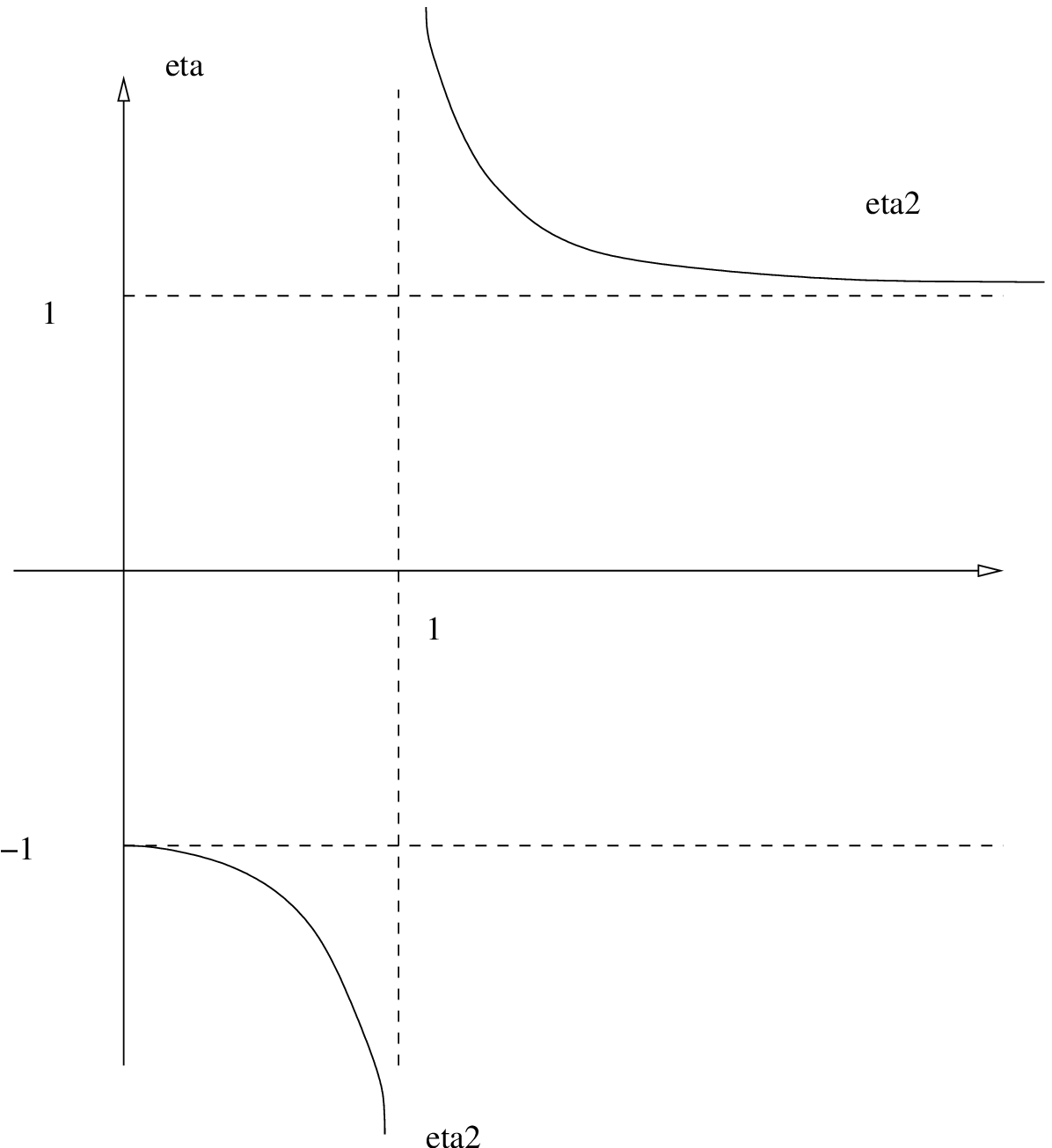}
\caption{The graph of the elasticity function for $H_-$.}\label{Figuree8}
\end{center}
\end{figure}

Now the characteristic equation for $H=H_-(t)$ is easily verified
\begin{eqnarray}
{\mathcal L}H_- &=& \left|H_-\right|^n \left[1+\frac{\eta^2}{n-1}\right]^\frac{n-2}{2}\left(1-\eta^2\right)\nonumber \\
&=& \left|\frac{H_-}{u}\right|^n (u^2-1) \left[u^2+\frac{1}{n-1}\right]^\frac{n-2}{2}\equiv  -1
\end{eqnarray}
Our analysis of convexity properties of $H_-$ is based on the equation (\ref{243}), which we write as
\begin{equation}
(1+u^2)t^2\, \ddot{H}=(u-1)\left[(n-1)u^2+(n-2)u+1\right]\frac{H}{u}
\end{equation}
for all $t\not= 1$. Since the quotient $\frac{H}{u}$ is positive we see that
\begin{equation}
\mbox{sgn}\,  \ddot{H}_{-}=- \mbox{sgn}  \left[(n-1)u^2+(n-2)u+1\right]
\end{equation}
The key observation is that for $n=2,3,4,5,6$ the polynomial $P=P_n(u)= (n-1)u^2+(n-2)u+1$ has no roots, meaning that $H_-$ is concave. Rather unexpectedly, $H_-$ is no longer concave when $n \geqslant 7$, because $P_n$ has two roots in the interval $u \in (-1, 0)$.

\begin{figure}[!h]
\begin{center}
\psfrag{-1}{\tiny ${ -1}$}
\psfrag{u}{\small ${ u}$}
\psfrag{Q}{\small $P$}
\psfrag{Qn}{\small $P=P_n(u)$}
\includegraphics*[height=1.6in]{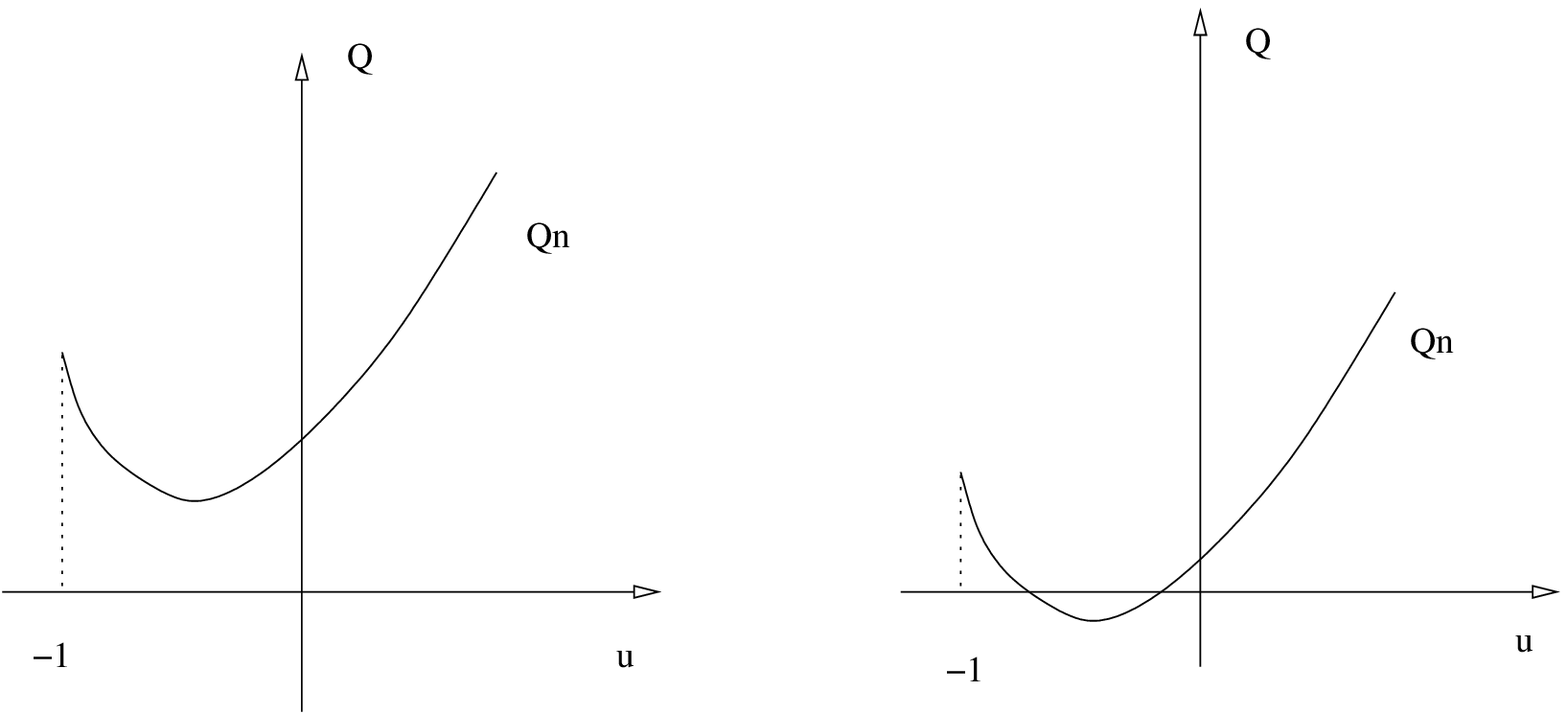}
\caption{$P_n$ is positive if $2\leqslant n \leqslant 6$ and has two roots if $n \geqslant 7$.}
\end{center}
\end{figure}

As a consequence, the graphs of $H_-$ exhibit two inflection points when $n \geqslant 7$.

\begin{figure}[!h]
\begin{center}
\psfrag{1}{\tiny ${ 1}$}
\psfrag{t}{\small ${ t}$}
\psfrag{H}{\small $H$}
\psfrag{H-}{\small $H=H_-(t)$}
\psfrag{n1}{\small $n=2,3,4,5,6$}
\psfrag{n2}{\small $n=7,8,9,...$}
\includegraphics*[height=1.8in]{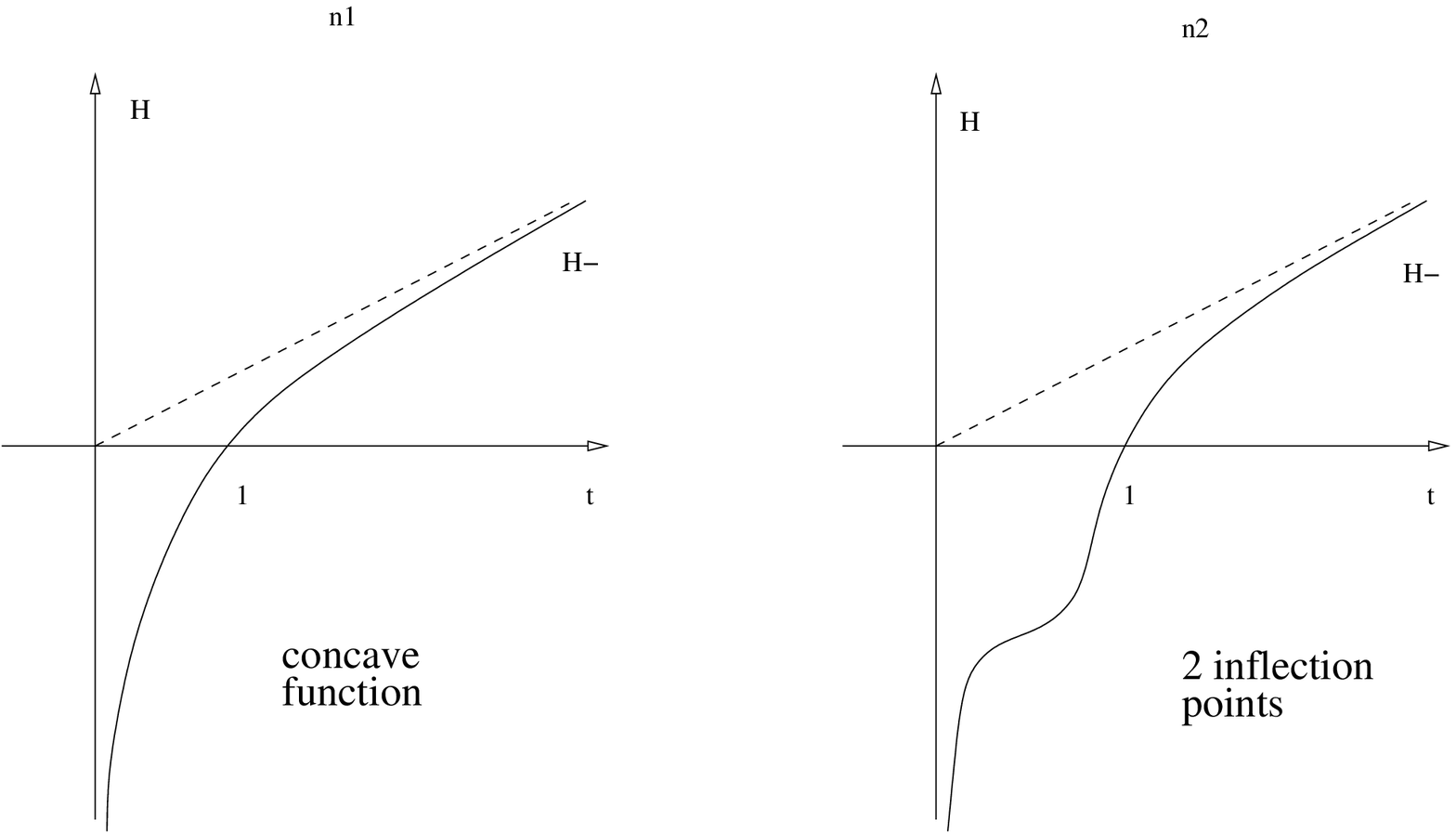}
\caption{The principal solution $H_-$.}\label{figur10}
\end{center}
\end{figure}

Analysis of the asymptotic behavior of $H_-$ is much the same as
that for $H_+$. We only  state the results. The slope of the
asymptote at infinity equals
\begin{equation}
\Theta = \Theta_- = \left(1-\frac{1}{n}\right)^\frac{n-2}{2n} 4^{-\frac{1}{n}}\,   \mbox{\large{exp}}\left[\frac{2-n}{n\sqrt{n-1}}\tan^{-1}\sqrt{n-1}  \right]
\end{equation}
and we have
\begin{equation}
0<\Theta\, t - H_-(t) \leqslant \frac{C}{t^{n-1}}
\end{equation}
In particular, for $H=H_-(t)$ it holds
\begin{equation}\label{E118}
\lim\limits_{t\rightarrow \infty} \frac{H(t)}{t}=\Theta = -\lim\limits_{t\rightarrow 0} t \, H(t)
\end{equation}
\section{The boundary value problem for radial $n$-harmonics}
Given $\alpha, \beta \in {\mathbb R}$ one may wish to look for
functions $H\in {\mathscr C}^2 (a,b)\cap {\mathscr C} [a,b]$, $0<a<b<\infty$, which  satisfy the conditions
\begin{equation}\label{equ561}
\begin{cases}{\mathcal L}H\equiv c & \\
H(a)=\alpha &\\
H(b)=\beta &
\end{cases}
\end{equation}
Here the constant $c$ is also viewed as unknown, otherwise the system ~\eqref{equ561} would be overdetermined (ill-posed). Let us begin with
the easy case $\beta =0$. The solution is given by the formula
\begin{equation}
H(t)=\lambda H_-(kt)\, , \hskip2cm k=\frac{1}{b} \hskip1cm \lambda = \frac{\alpha}{H\left(\frac{a}{b}\right)}
\end{equation}
From now on we assume that $\beta \not= 0$. There are five cases to consider.\\
{\bf Case 1.} Suppose that
\begin{equation}
-\infty < \frac{\alpha}{\beta}< \frac{a}{b}
\end{equation}
The solution will be found in the form $H(t)=\lambda H_-(kt)$.
One has to show that there exist real number $\lambda$ and a positive number $k$ such that
\begin{equation}\label{261}
\begin{cases}\lambda H_-(ka)=\alpha &\\
\lambda H_-(kb)=\beta &
\end{cases}
\end{equation}
We eliminate $\lambda$ by dividing the equations
\begin{equation}\label{ES119}
Q(k)= \frac{H_-(ka)}{H_-(k b)}\, , \hskip2cm \mbox{for }
k>\frac{1}{b}
\end{equation}
Note that $\lim\limits_{k \rightarrow \frac{1}{b}} Q(k) = -\infty$. On the other hand using (\ref{E118}) we see that
\begin{equation}\label{ES120}
\lim\limits_{k \rightarrow \infty}  Q(k) = \frac{a}{b}\left(
\lim\limits_{k \rightarrow \infty} \frac{\; \; \frac{H_-(ka)}{ka}\;
\; }  {\; \; \frac{H_-(kb)}{kb}\; \; } \right) = \frac{a\,
\Theta_-}{b\, \Theta_-} = \frac{a}{b}
\end{equation}
By Mean Value Theorem there exists $\infty > k_\circ > \frac{1}{b}$ such that
\begin{equation}
Q(k_\circ )=\frac{H_-(k_\circ a)}{H_-(k_\circ b)}=\frac{\alpha}{\beta}
\end{equation}
Finally, conditions (\ref{261}) are satisfied with
$\lambda = \frac{\beta}{H_-(k_\circ b)}$.

The next  case is obvious; a  linear function is a solution.\\
{\bf Case 2.} Suppose
\begin{equation}
\frac{\alpha}{\beta}=\frac{a}{b}
\end{equation}
Then the solution is given by
\begin{equation}
H(t)= {\small \mbox{$\frac{\alpha}{\beta}$}} \, t=\alpha \, H_\circ \left({\small \mbox{$\frac{t}{a}$}}\right)
\end{equation}
{\bf Case 3.}  Suppose now that
\begin{equation}
\frac{a}{b} < \frac{\alpha}{\beta}< \frac{b}{a}
\end{equation}
We are looking for the solution in the form $ H(t)=\lambda H_+(kt)$.
As before, one needs to find $k_\circ$ such that
\begin{equation}\label{263}
\frac{H_+(k_\circ a)}{H_+(k_\circ b)}=\frac{\alpha}{\beta}
\end{equation}
For this reason we consider the function
\begin{equation}
Q(k)= \frac{H_+(k a)}{H_+(k b)}\, , \hskip2cm k>0
\end{equation}
With the aid of the asymptotic formula (\ref{E104}) we see that
$\lim\limits_{k\rightarrow \infty} Q(k) = \frac{a}{b}$ and
$\lim\limits_{k\rightarrow 0} Q(k) = \frac{b}{a}$. By Mean
Value Theorem the equation (\ref{263}) holds for some $k_\circ >0$.
Then  the parameter $\lambda$ can be chosen to yield the boundary
conditions at (\ref{261}).

Next case is immediate.\\
{\bf Case 4.}   Suppose
\begin{equation}
\frac{\alpha}{\beta}=\frac{b}{a}
\end{equation}
The solution is given by the formula $H(t)=\frac{\alpha
a}{t}=\alpha\,  H_\infty \left(\frac{t}{a}\right)$.\\
{\bf Case 5.}  Suppose finally that
\begin{equation}
\frac{\alpha}{\beta}> \frac{b}{a}
\end{equation}
We are looking for  the solution in the form $H(t)=\lambda H_-(kt)$, so we examine  the function
\begin{equation}
Q(k)= \frac{H_-(k a)}{H_-(k b)}\, , \hskip2cm \mbox{ with } k<\frac{1}{b}
\end{equation}
By virtue of  formula (\ref{E118}) we see that
$\lim\limits_{k\rightarrow \frac{1}{b}} Q(k) = +\infty$ and  $
\lim\limits_{k\rightarrow 0} Q(k) = \frac{b}{a}$. This case is
completed by again invoking Mean Value Theorem;  there is $k_\circ <
\frac{1}{b}$ such that
\begin{equation}
\frac{H_-(k_\circ a)}{H_-(k_\circ b)}=\frac{\alpha}{\beta}
\end{equation}
Then, we take $\lambda$   to ensure that $\lambda H_- (k_\circ b)=\beta$.
\vskip0.5cm

In summing up this section, we note that the foregoing analysis,
combined with the uniqueness of the Dirichlet problem for the
$n$-harmonic equation, reveals that in fact the principal $n$-harmonics
generate all radial $n$-harmonics. As an unexpected bonus we deduce
that every radial $n$-harmonic map in a ring domain extends
$n$-harmonically to the entire punctured space $\,\mathbb R^n_\circ = \mathbb R^n \setminus \{0\}\,$ and continuously to the M\"obius space $\,\hat{\mathbb R}^n = \mathbb R^n \cup \{\infty\}\,$.

\chapter{Vector calculus on annuli}
In this chapter we briefly review differential calculus on annuli.
\section{Radial and spherical derivatives}
In studying deformations of spherical rings one must consider the
radial and spherical components of a differential map. Through
every point $x\in \rn_\circ$ there passes a sphere
\begin{equation}
{\mathbb S}^{n-1}_t =\big\{y\in \rn \, ; \; \; \; \; |y|=t\big\}\, , \hskip2cm t=|x|
\end{equation}
Its tangent hyperplane at $x$ is given by
\begin{equation}
{\mathbb T}_x = \big\{\xi\, ; \; \; \; \langle x, \xi \rangle =t^2 \big\}
\end{equation}
Fix an orthonormal basis for ${\mathbb T}_x$, say ${\bold T}=\{T_2,T_3, ..., T_n\}$.\footnote{It is not always possible to choose such a basis continuously depending on $x\in \rn_\circ$} Let $\Omega \subset \rn_\circ$ be a domain and $h: \Omega \rightarrow \rn$ a mapping having first order partial derivatives defined at $x\in \Omega$. The radial (or normal) derivative of $h$ at $x$ is a vector defined by the rule
\begin{equation}
h_{_{\bold N}} (x)= Dh(x)\cdot {\bold N}=\frac{x_1 h_{x_1}+ ... +
x_n h_{x_n}}{|x|}
\end{equation}
where ${\bold N}= \frac{x}{|x|}$ is called the  radial (or normal) vector. In an exactly similar way we define the spherical (or tangential) derivatives
\begin{equation}
h_{\bold T}(x) =\left[h_{T_2}, ..., h_{T_n}\right]\, , \; \;
h_{T_\nu}=Dh(x) \cdot T_\nu \, ,\; \; \nu = 2,3, ..., n
\end{equation}
This is an $(n-1)$-tuple of vectors in $\rn$, conveniently
considered as column vectors of an $n \times (n-1)$-matrix.
Continuing in this fashion we view the differential of $Dh$ as an $n
\times n$-matrix whose column vectors are $h_{\bold N}, h_{T_2},
..., h_{T_n}$. Of course, such a matrix representation of $Dh$
depends on choice of the orthonormal frame ${\bold N}, T_2, ...,
T_n$
\begin{equation}\label{311}
Dh(x) \approx \left[ \begin{array}{cccc}  | & | & \; \; \; & | \\
 | & | & \; \; \; & | \\
h_{{\bold N}}& h_{T_2} &  \; \; \; &  h_{T_n} \\
 | & | & \; \; \; & | \\
 | & | & \; \; \; & |
\end{array}\right]
\end{equation}
However, some differential expressions do not depend on the frame. For
example,  the Hilbert-Schmidt norm of the
differential matrix,
\begin{equation}
\norm Dh \norm^2 = \left|h_{\bold N}\right|^2 +
\left|h_{T_2}\right|^2+ ... + \left|h_{T_n}\right|^2
\end{equation}
Also, the spherical component of $\norm Dh \norm$, which we  define
by the formula
\begin{eqnarray}
\left|h_{\bold T}\right|=  \left(\frac{ \left|h_{T_2}\right|^2+ ...
+ \left|h_{T_n}\right|^2}{n-1}\right)^\frac{1}{2} \geqslant  \Big(
\left|h_{T_2}\right|\cdots \left|h_{T_n}\right|\Big)^\frac{1}{n-1}
\end{eqnarray}
is frame free. Equality holds if and only if $\left|h_{T_2}\right|=
... = \left|h_{T_n}\right|$. Thus we have
\begin{equation}\label{E142}
\norm Dh \norm^2 = \left|h_{\bold N}\right|^2 + (n-1) \left|h_{\bold
T}\right|^2
\end{equation}
Another differential quantity of interest to us is the Jacobian determinant
\begin{equation}
J(x,h)= \mbox{det}\, Dh = \langle h_{\bold N}\, ,\,  h_{T_2}\times ... \times h_{T_n} \rangle
\end{equation}
where the cross product of spherical derivatives is controlled by
using the point-wise Hadamard's inequality.
\begin{equation}
\left| h_{T_2}\times ... \times h_{T_n}  \right| \leqslant \left| h_{T_2} \right| \cdots  \left| h_{T_n} \right|
\end{equation}
Hence
\begin{equation}\label{145}
J(x,h) \leqslant  \left| h_{\bold N} \right| \cdot  \left| h_{\bold
T} \right|^{n-1}
\end{equation}
Equality occurs if and only if the vectors $ h_{\bold N},  h_{T_2},
..., h_{T_n}$ are mutually orthogonal, positively oriented and
$\left|h_{T_2}\right|=...=\left|h_{T_n}\right|$. This amounts to
saying that the Cauchy-Green tensor of $h$ must be a diagonal matrix
\begin{equation}\label{E777}
{\bold C}(x,h) \bydef D^\ast h\cdot  Dh = \left[ \begin{array}{cccc}  \left| h_{N} \right|^2 & 0 & \cdots & 0  \\
 0 & \left| h_{T} \right|^2 &   \cdots & 0\\
& &\ddots &  \\
0 & 0 & \cdots & \left| h_{T} \right|^2
\end{array}\right]
\end{equation}
The notion of elasticity that we have introduced for radial
stretchings can be extended to all weakly differentiable mappings.
Let $h : {\mathbb A} \rightarrow {\mathbb A}^\ast$ be any mapping of
annuli, and $x_\circ \in {\mathbb A}$ be a point where $h$ is
differentiable.\footnote{Homeomorphisms  of class ${\mathscr
W}^{1,n}$ are differentiable at almost every point.} We say that, at
this point, $h$ is
\begin{itemize}
\item{{\it expanding}, if $\left|\, h_N(x_\circ)\, \right| > \left|\, h_T(x_\circ)\, \right|$\, , \hskip0.7cm $x_\circ$ is a point of
expansion}
\item{{\it contracting}, if $\left|\, h_N(x_\circ)\, \right| < \left|\, h_T(x_\circ)\, \right|$\, ,  \hskip0.5cm $x_\circ$ is a point of contraction}
\end{itemize}
In what follows we shall be concerned only with the regular points,
that is $\norm Dh (x_\circ ) \norm \not= 0$. We define the  {\it elasticity module} of $h$  at the regular points by
the rule
\begin{equation}
\eta = \eta_{h} (x) = \frac{\left|\, h_N(x)\, \right|}{\left|\,
h_T(x)\, \right|}\in \hat{{\mathbb R}}
\end{equation}
This is a measurable function taking values in the extended half line $[0, \infty]$. Note the following inequalities,
\begin{eqnarray}
K_{_O} (x,h)&= & \frac{|Dh(x)|^n}{J(x,h)} \geqslant \frac{\left|\, h_N\, \right|^n}{\left|\, h_N\, \right|\, \left|\, h_T\, \right|^{n-1}}= [\eta_{h}(x)]^{n-1} \\
{\mathbb K}_{_O} (x,h) &=& \frac{\norm Dh
\norm^n}{n^\frac{n}{2}J(x,f)} \geqslant \frac{\left(\left|\, h_N\,
\right|^2 + (n-1)\left|\, h_T\, \right|^2
\right)^\frac{n}{2}}{n^\frac{n}{2}\left|\, h_N\, \right|\, \left|\,
h_T\, \right|^{n-1}}\nonumber \\ &=&
\frac{1}{\eta_{h}}\left(\frac{n-1}{n} + \frac{1}{n}
\eta^2_{h}\right)^\frac{n}{2}
\end{eqnarray}
We say that the deformation  $h : {\mathbb A} \rightarrow {\mathbb
A}^\ast$ at the given point $x_\circ$ is:
\begin{itemize}
\item{{\it contracting}, if $\eta_{h}(x_\circ) <1$}
\item{{\it expanding}, if $\eta_{h}(x_\circ) >1$}
\end{itemize}
It will be clear later that the extremal deformations  have the
same elasticity type throughout the region. Furthermore, the
elasticity module will not vanish inside the region. Thus for the
extremal deformations the equation at (\ref{E777}) takes the form of
a Beltrami system
\begin{equation}
D^\ast h \cdot Dh = J(x,h)^\frac{2}{n}{\bold K}(x)
\end{equation}
where the distortion tensor ${\bold K}$ is a diagonal matrix of
determinant 1
\begin{equation}
 {\bold K} = \left[ \begin{array}{cccc} \eta^{2-\frac{2}{n}} & 0 & \cdots & 0  \\
 0 & \eta^{-\frac{2}{n}} &   \cdots & 0\\
& &\ddots &  \\
0 & 0 & \cdots & \eta^{-\frac{2}{n}}
\end{array}\right]
\end{equation}
For the radial stretching
\begin{equation}
h(x)=H(|x|)\frac{x}{|x|}
\end{equation}
we find that
\begin{equation}\label{E147}
h_{\bold N}= \dot{H}(|x|) \cdot {\bold N}\, , \; \; \; \; \; \;
h_{T_i}= \frac{H(|x|)}{|x|}\cdot T_i
\end{equation}
Hence, the cross product of the spherical derivative is parallel to the radial vector.
\begin{equation}
 h_{T_2}\times ... \times h_{T_n} = \left(\frac{H}{|x|}\right)^{n-1} \frac{x}{|x|}
\end{equation}
The square of Hilbert-Schmidt norm of the differential is given by
\begin{equation}
\norm Dh \norm^2 = \dot{H}^2 + (n-1) \frac{H^2}{t^2}\, , \hskip2cm t=|x|
\end{equation}
and the Jacobian determinant equals:
\begin{equation}
J(x,h) = \frac{\dot{H}\, H^{n-1}}{t^{n-1}}
\end{equation}
\section{Some differential forms}
Integration of nonlinear differential expressions is best handled by
using differential forms. Let us begin with the function $t=|x|$ and
its differential
\begin{equation}
dt = \frac{x_1 dx_1+...+ x_ndx_n}{|x|}
\end{equation}
The Hodge star of $\d t$ is an $(n-1)$-form.
\begin{equation}
\star \, dt = \frac{1}{|x|} \sum_{i=1}^n (-1)^i x_i\, dx_1 \wedge ... \wedge dx_{i-1}\wedge d x_{i+1} \wedge ... \wedge dx_n
\end{equation}
This form integrates naturally on spheres centered at the origin. We
wish to normalize  this form in such a way that the integrals  will be
independent of the sphere. In this way we come to what is known as
standard area form on ${\mathbb S}^{n-1}$
\begin{equation}\label{E153}
\omega (x) = \frac{\star \, \d t}{t^{n-1}}= \sum_{i=1}^n (-1)^{i}\,
\frac{x_1 dx_1 \wedge ... \wedge dx_{i-1}\wedge dx_{i+1}\wedge ...
\wedge dx_n}{|x|^n}
\end{equation}
Viewing $\omega$ as a differential form on punctured space $\rn_\circ$, we find that $d\omega =0$. By Stokes' theorem
\begin{equation}\label{omef}
\int_{|x|=t}\omega = \omega_{n-1}\, ,  \hskip2cm \mbox{for all } t>0
\end{equation}
Next, given any mapping $h=(h^1, ..., h^n) : \Omega \rightarrow
\rn_{\, \circ}$ of  class ${\mathscr W}^{1,1}_{\loc} (\Omega ,
\rn_\circ)$, we  consider the pullback of $\omega$ via $h$
\begin{equation}
h^\sharp \omega = \sum_{i=1}^n (-1)^i \frac{h^i dh^1 \wedge ... \wedge dh^{i-1}\wedge dh^{i+1}\wedge ... \wedge dh^n}{|h|^n}
\end{equation}
Under suitable regularity hypothesis, for instance if  $h\in
{\mathscr W}^{1,n-1}_{\loc} (\Omega , \rn)$, and $|h(x)|\geqslant
\mbox{\it{const}}>0$, this form is also closed, meaning that $d
\left(h^\star \omega \right)= h^\star \left(d \omega \right) =0$. We
view the exterior derivative of $h$ as an $n$-tuple of 1-forms,
$dh=\left[dh^1, ..., dh^n\right]$. With such a view $dh \wedge dt$
becomes an $n$-tuple of 2-forms. Further notation is self
explanatory. For example,
\begin{equation}
\left|dh \wedge dt\right|^2 = \left|dh^1 \wedge dt\right|^2 + ... + \left|dh^n \wedge dt\right|^2
\end{equation}
Now the decomposition at (\ref{E142}) takes the form
\begin{equation}\label{320}
\norm Dh \norm^2 =|dh|^2 = |dh \wedge \star dt|^2 + |dh \wedge dt|^2
\end{equation}
The reader may verify this formula by the following
general algebraic identity $|{\mathfrak a}|^2 = |{\mathfrak a}
\wedge \star N|^2 + |{\mathfrak a} \wedge  N|^2 $  for all covectors
${\mathfrak a} \in \bigwedge^1 (\rn)$. Just apply it to ${\mathfrak
a}^i = \d h^i$, $i=1,2, ...,n$ and $N=dt$, then add the resulting
identities.  A moment of reflection about (\ref{320}) reveals that
\begin{equation}\label{321hal}
dh \wedge \star dt =h_N\, dx
\end{equation}
and
\begin{equation}\label{321}
dh \wedge dt =\frac{1}{2|x|} \sum_{1\leqslant i<j \leqslant n} \left(x_ih_{x_j}-x_jh_{x_i}\right)dx_i \wedge dx_j
\end{equation}
Hence
\begin{equation}\label{322}
\left| dh \wedge  dt\right|^2 =(n-1)|h_T |^2
\end{equation}
\begin{proposition}\label{pr321}
The following point-wise estimates hold
\begin{equation}\label{323}
\left|h_N\right| \geqslant \left|\, d|h|\wedge \star dt\,  \right|
\end{equation}
\begin{equation}\label{324}
\left|h_T\right|^{n-1} \geqslant |h|^{n-1} \left|dt\wedge h^\sharp \omega \right|
\end{equation}
\end{proposition}
\proof Formula (\ref{321hal}) holds for scalar functions in place of
$h$ we apply it to $|h|$
\begin{equation}
\left|h_N \right|\,  \geqslant \left|\, |h|_N\, \right|= \left|\, d|h|\wedge \star dt\, \right|
\end{equation}
establishing the inequality (\ref{323}).

The proof of the estimate (\ref{324}) is more involved. Using
(\ref{322}), we are reduced to showing that
\begin{equation}
\left|\frac{dh \wedge dt }{n-1}\right|^{n-1} \geqslant \left|\sum_{i=1}^n (-1)^i \frac{h^i}{|h|}dh^1 \wedge ... \wedge dh^{i-1}\wedge dh^{i+1}\wedge ... \wedge dh^n \wedge dt\right|
\end{equation}
By Schwarz inequality it suffices to prove that
\begin{equation}
\left|\frac{dh \wedge dt }{n-1}\right|^{n-1} \geqslant \left(\sum_{i=1}^n \left|dh^1 \wedge ... \wedge dh^{i-1}\wedge dh^{i+1}\wedge ... \wedge dh^n \wedge dt\right|^2 \right)^\frac{1}{2}
\end{equation}
Note that here in both sides  only spherical derivatives of $h$ are
significant; the terms containing $h_N\,  dt$  vanish after wedging
them with $dt$. This observation permits us to replace $dh^i$ by the
covectors ${\mathfrak a}_i=dh^i -h_N^i dt$. We view ${\mathfrak a}_i
$ as elements of the space $\bigwedge^1({\mathbb R}^{n-1})$. Once
this interpretation is accepted, the proof continues via  an
algebraic inequality.
\begin{lemma}\label{Le311}
For every $n$-tuple of covectors in ${\mathbb R}^{n-1}$, ${\mathfrak a}_1, ..., {\mathfrak a}_n \in \bigwedge^1 ({\mathbb R}^{n-1})$, we have
\begin{equation}\label{E151}
\left|{\mathfrak a}_1 \right|^2 +...+ \left|{\mathfrak a}_n \right|^2 \geqslant (n-1) \left(\sum_{i=1}^n \left|{\mathfrak a}_1\wedge ... \wedge {\mathfrak a}_{i-1} \wedge {\mathfrak a}_{i+1} \wedge ... \wedge {\mathfrak a}_n   \right|^2 \right)^\frac{1}{n-1}
\end{equation}
\end{lemma}
\proof
We look at the matrix $M$ of size $n \times (n-1)$ whose rows are made of the covectors ${\mathfrak a}_1, ..., {\mathfrak a}_n$. Let ${\mathfrak b}_2, ..., {\mathfrak b}_n\in {\mathbb R}^{n}$ denote the column vectors of~$M$
\begin{equation}
\left[ \begin{array}{ccc} &- \; \; - \; \; {\mathfrak a}_1 \; \; - \; \; -  & \\
 & - \; \; - \; \; {\mathfrak a}_2 \; \; - \; \; - & \\
 &\ddots &  \\
 & - \; \; - \; \; {\mathfrak a}_n \; \; - \; \; - &
\end{array}\right] = M = \left[ \begin{array}{ccc}  &|& |  \hskip0.8cm | \\
&|& | \hskip0.8cm | \\
&{\mathfrak b}_2 & {\mathfrak b}_3   \cdots    {\mathfrak b}_n \\
 &|& |  \hskip0.8cm | \\
&|& |  \hskip0.8cm |
\end{array}\right]
\end{equation}
In the left hand side of (\ref{E151}) we are actually dealing with
the square of the Hilbert-Schmidt norm of $M$
\begin{equation}
\left|{\mathfrak a}_1 \right|^2 + ... + \left|{\mathfrak a}_n \right|^2 = \mbox{Tr}\, \left( M^\ast M\right) = \left|{\mathfrak b}_2 \right|^2 + ... + \left|{\mathfrak b}_n \right|^2
\end{equation}
On the other hand, the sum in the right hand side of (\ref{E151}) is
the Hilbert-Schmidt square of the norm of the cross product of the
column vectors
\[
\sum_{i=1}^n  \left|{\mathfrak a}_1\wedge ... \wedge {\mathfrak a}_{i-1} \wedge {\mathfrak a}_{i+1} \wedge ... \wedge {\mathfrak a}_n   \right|^2  = \left| {\mathfrak b}_2 \times ... \times {\mathfrak b}_n\right|^2 \leqslant \left| {\mathfrak b}_2 \right|^2 \cdots  \left| {\mathfrak b}_n \right|^2
\]
by Hadamard's inequality. The inequality between arithmetic and geometric means yields
\begin{equation}
 \left| {\mathfrak b}_2 \right|^2 \cdots  \left| {\mathfrak b}_n \right|^2 \leqslant \left( \frac{\left| {\mathfrak b}_2 \right|^2 + ... + \left| {\mathfrak b}_n \right|^2}{n-1}\right)^{n-1}
\end{equation}
From here Lemma \ref{Le311} is immediate. This also completes the proof of Proposition \ref{pr321}.
\chapter{Free Lagrangians}
We are about to introduce the basic concept of this paper. Given two domains $\Omega$ and $\Omega^\ast$ in $\rn$, we shall
consider orientation preserving homeomorphisms $h: \Omega
\rightarrow \Omega^\ast$ in a suitable Sobolev class ${\mathscr
W}^{1,p}(\Omega, \Omega^\ast)$ so that a given energy integral
\begin{equation}\label{elasE}
I[h]=\int_\Omega {\mathcal E}(x, h, Dh)\, dx
\end{equation}
is well defined. The term {\it free Lagrangian} pertains to a
differential $n$-form ${\mathcal E}(x,h, Dh)\, dx$ whose integral
depends only on the homotopy class of $h : \Omega
\overset{\textnormal{\tiny{onto}}}{\longrightarrow} \Omega^\ast$.
The Jacobian determinant is pretty obvious example; we have,
\begin{equation}
 \int_\Omega J(x,h)\, dx = \left|\Omega^\ast \right|
\end{equation}
This identity holds for all orientation preserving homeomorphisms of
Sobolev class ${\mathscr W}^{1,n}(\Omega, \Omega^\ast)$. Many more
differential expressions enjoy a property such as this. In the next
three lemmas we collect examples of free Lagrangians for
homeomorphisms $h :{\mathbb A}
\overset{\textnormal{\tiny{onto}}}{\longrightarrow} {\mathbb
A}^\ast$ of annuli ${\mathbb A}=\{ x\, ; \; \; \; r <|x| < R\}$ and
${\mathbb A}^\ast =\{ x\, ; \; \; \; r_\ast < |x| < R_\ast \}$.
\begin{lemma}\label{Le331}
Let $\Phi : [r_\ast , R_\ast] \rightarrow {\mathbb R}$ be any
integrable function. Then the $n$-form
\begin{equation}
\Phi(|h|)\, dh^1 \wedge ... \wedge dh^n
\end{equation}
is  a free Lagrangian. Precisely, we have
\begin{equation}\label{Jaco}
 \int_{\mathbb A} \Phi(|h|) J(x,h)\, dx = \omega_{n-1}\int_{r_\ast}^{R_\ast} \tau^{n-1} \, \Phi(\tau)\, d\tau
\end{equation}
for every orientation preserving homeomorphism  $h\in {\mathscr
W}^{1,n}({\mathbb A}, {\mathbb A}^\ast)$.
\end{lemma}
This is none other than a  general formula of integration by
substitution.
\begin{lemma}\label{Le332}
The following differential $n$-form
\begin{equation}
\sum_{i=1}^n \frac{x_i\, dx_1 \wedge ... \wedge dx_{i-1}\wedge d|h|\wedge dx_{i+1}\wedge ... \wedge d x_n}{|h|\, |x|^n} =  \frac{\left(d|h|\right)\wedge \star dt}{|h|\, t^{n-1}}
\end{equation}
is  a free Lagrangian in the class of all homeomorphisms $h\in {\mathscr
W}^{1,1}({\mathbb A}, {\mathbb A}^\ast)$ that preserve the order of the
boundary components of the annuli ${\mathbb A}$ and ${\mathbb
A}^\ast$. Precisely, we have
\begin{equation}\label{Eq178}
\int_{\mathbb A} \frac{d|h|\; \wedge \star dt}{|h|\,
t^{n-1}}=\textnormal{Mod}\, {\mathbb A}^\ast
\end{equation}
\end{lemma}
\proof  First observe that the function $|h|: {\mathbb A}
\rightarrow (r_\ast, R_\ast)$ extends continuously to the closure of
${\mathbb A}$ so that $|h(x)|=r_\ast$ for $|x|=r$ and
$|h(x)|=R_\ast$ for $|x|=R$. The point is that  our arguments will
actually work for arbitrary continuous mappings $h :
\overline{{\mathbb A}}\rightarrow \rn_{\, \circ}$ of Sobolev class $
{\mathscr W}^{1,1}({\mathbb A}, \rn)$ which satisfy the above
boundary conditions, not necessarily homeomorphisms. The advantage
in such a generality is that we can apply smooth approximation of
$h$ without worrying about injectivity.  This extends Lemma
\ref{Le332} to all permissible mappings. If $h$ is smooth it is
legitimate to apply Stokes' theorem. Recalling the form $\omega$ and
(\ref{omef}), we write the integral at (\ref{Eq178}) as:
\begin{eqnarray}
 \int_{\mathbb A} \left(d \log |h|\right) \wedge \omega = \int_{\mathbb A} d\left(\omega\,  \log |h|\right)
&=& \left(\log R_\ast \right)\int_{|x|=R}\omega -  \left(\log r_\ast
\right)\int_{|x|=r}\omega\nonumber \\ &=& \omega_{n-1} \log
\frac{R_\ast}{r_\ast} = \mbox{Mod}\, {\mathbb A}^\ast
\end{eqnarray}
The same proof works for more general integrals of the form
\begin{equation}\label{Eq178p}
\int_{\mathbb A} \Phi^\prime \left(\, |h|\, \right) \frac{d|h|\wedge
\star dt}{t^{n-1}}=\omega_{n-1} \left[\Phi(R_\ast) -
\Phi(r_\ast)\right]
\end{equation}
with $\Phi \in {\mathscr C}^1 [r_\ast, R_\ast]$. This also reads as
\begin{equation}\label{Eq188p}
\int_{\mathbb A} \frac{\Phi^\prime (\,
|h|\,)\left|h\right|_N}{|x|^{n-1}} = \omega_{n-1}\left[
\Phi(R_\ast)- \Phi(r_\ast)\right]
\end{equation}
\vskip0.2cm

Another example of a free Lagrangian rests   on the topological degree. Let
$\Psi : S_t^{n-1} \rightarrow {\mathbb S}^{n-1}$ be any smooth
mapping of a sphere ${S}_t^{n-1} =\{x\, ; \; |x|=t\}$ into the unit
sphere $\mathbb S^{n-1} \subset \rn$. Then there is an integer, denoted by $\mbox{deg}\,
\Psi$, and called the {\it degree} of $\Psi$, such that
\begin{equation}
\int_{S_t^{n-1}} \Psi^\sharp (\omega)=\omega_{n-1}\, \mbox{deg}\,
\Psi
\end{equation}
where $\omega$  stands for the standard surface measure of  the unit sphere
$ {\mathbb S}^{n-1} \subset \rn$, as given in~(\ref{E153}). Note
a general fact that $\mbox{deg}\, \Psi$ remains unchanged under
small perturbations of $\Psi$. If $\Psi$ is an
orientation-preserving diffeomorphism, then $\mbox{deg}\, \Psi =1$.
Now a free Lagrangian somewhat  dual to that in Lemma \ref{Le332}, is
obtained as follows.
\begin{lemma}\label{Le333}
The following  differential $n$-form
\begin{equation}
 \sum_{i=1}^n \frac{h^i\, dh^1 \wedge ... \wedge dh^{i-1}\wedge d|x|\wedge dh^{i+1}\wedge ... \wedge d h^n}{|x|\, |h|^n}
\end{equation}
is a free Lagrangian in the class of all orientation preserving homeomorphism
mappings $h \in {\mathscr W}^{1,n-1}({\mathbb A}, {\mathbb
A}^\ast)$. Precisely, we have
\begin{equation}\label{331}
\int_{\mathbb A} \frac{dt}{t}\wedge h^\sharp \omega =  \textnormal{Mod}\, {\mathbb A}
\end{equation}
\end{lemma}
\proof As a first step, we approximate $h$ by smooth mappings $h_\nu$
converging  to $h$ $c$-uniformly on ${\mathbb A}$ and in the norm of
$ {\mathscr W}^{1,n-1}({\mathbb A},\rn)$. For $\nu $ sufficiently
large we consider the mappings
\begin{equation}\label{E181}
\Psi =\frac{h_\nu}{|h_\nu|}: {\mathbb A}^\prime \rightarrow {\mathbb
S}^{n-1}\, , \; \; \; \; \nu=1,2,...
\end{equation}
defined on a slightly smaller annulus ${\mathbb A}^\prime =\{\, x \,
; \; \; \; r^\prime \leqslant |x| \leqslant R^\prime\}$, where
$r<r^\prime < R^\prime < R$. The degree of $\Psi$ restricted to any
concentric sphere in ${\mathbb A}^\prime$ is equal to 1. Moreover,
we have the point-wise identity $\Psi^\sharp \omega = h_\nu^\sharp
(\omega )$. Now we  integrate as follows
\begin{equation}
\int_{\mathbb A^\prime} \frac{dt}{t}\wedge \Psi^\sharp \omega =
\int_{\mathbb A^\prime} d \left(\log |t|\,  \Psi^\sharp \omega
\right)
\end{equation}
because $d\left(\Psi^\sharp \omega \right)=\Psi^\sharp \left(d
\omega \right)=0 $. By Stokes' theorem, the right hand side is equal
to:
\begin{eqnarray}
\int_{|x|=R^\prime} \log |x|\,  \Psi^\sharp \omega -
\int_{|x|=r^\prime} \log |x|\,  \Psi^\sharp \omega =
\omega_{n-1}\log R^\prime\,  -\, \omega_{n-1}\log r^\prime =
\mbox{Mod}\, {\mathbb A}^\prime \nonumber
\end{eqnarray}
Passing to the limit as $\nu \rightarrow \infty$ we conclude with
the formula
\begin{equation}
\int_{\mathbb A^\prime} \frac{dt}{t}\wedge h^\sharp \omega =
\textnormal{Mod}\, {\mathbb A^\prime}
\end{equation}
Finally letting $r^\prime \rightarrow r$ and $R^\prime \rightarrow
R$ yield (\ref{331}).
\begin{remark}
The reader may wish to generalize  this lemma to mappings $h:
{\mathbb A} \rightarrow \rn_\circ$ which are not necessarily
homeomorphisms. The degree of $\frac{h}{|h|}$ will emerge as a
factor in front of the right hand side of (\ref{331}).
\end{remark}
The following corollary from Lemmas \ref{Le331}, \ref{Le332},
\ref{Le333} provide us with three inequalities, which we call
free Lagrangian estimates.
\begin{corollary}
Let $h$ be a homeomorphism between spherical rings ${\mathbb A}$ and
${\mathbb A}^\ast$ in the Sobolev class  $ {\mathscr
W}^{1,n}({\mathbb A}, {\mathbb A}^\ast)$. Then
\begin{equation}\label{E184}
 \int_{\mathbb A} \Phi \big( |h| \big ) \left|h_N\right|\,  \left|h_T\right|^{n-1} \geqslant \omega_{n-1} \int_{r_\ast}^{R_\ast} \tau^{n-1} \Phi(\tau)\, d\tau
\end{equation}
whenever  $\Phi$ is  integrable in $[r_\ast, R_\ast]$. We have the
equality in (\ref{E184}) if and only if $\left|h_N\right|\,
\left|h_T\right|^{n-1}=J(x,h)$. Furthermore,
\begin{equation}\label{E185}
\int_{\mathbb A} \frac{\left|h_N\right|}{|h|\, |x|^{n-1}} \geqslant
\textnormal{Mod}\, {\mathbb A}^\ast
\end{equation}
\begin{equation}\label{E186}
\int_{\mathbb A} \frac{ \left|h_T\right|^{n-1}}{|h|^{n-1}\, |x|} \geqslant   \textnormal{Mod}\, {\mathbb A}
\end{equation}
\end{corollary}
Note that we have equalities  if $h$ is a radial mapping. Other cases
of equalities are also possible.
\chapter{Some  estimates of free Lagrangians}\label{cha7}
Before we proceed to the general proofs, it is instructive to take
on  some estimates which are at the heart of our application of free Lagrangians.

\section{The $\mathscr F_h$-energy integral with operator norm}
The extremal problem  is remarkably simpler if we use the operator
norm of the differential. Let $h$ be a permissible map of ${\mathbb
A}$ onto ${\mathbb A}^\ast$. The following point-wise inequality is
straightforward,
\begin{equation}
|Dh|\geqslant \max \{|h_N|\, , \, |h_T|\}
\end{equation}
Indeed, we have
\begin{equation}
|Dh|\geqslant \max \{|h_N| ,  |h_{T_2}|, ..., |h_{T_n}|\} \geqslant
\max \{|h_N|\, , \, |h_T|\}
\end{equation}
Certainly, equality occurs when $h$ is a radial stretching. All that
is needed is to apply the lower bounds at  (\ref{E185}) and
(\ref{E186}).
\begin{equation}\label{E428}
{\mathscr F}_h = \int_{\mathbb A} \frac{|Dh|^n}{|h|^n} \geqslant
\max \left\{ \int_{\mathbb A} \frac{|h_N|^n}{|h|^n}\, ,
\int_{\mathbb A} \frac{|h_T|^n}{|h|^n}   \right\}
\end{equation}
Here we estimate each term by H\"older's inequality,
\begin{eqnarray}
 \int_{\mathbb A} \frac{|h_N|^n}{|h|^n} &\geqslant &  \left(  \int_{\mathbb A} \frac{|h_N|}{|h|\, |x|^{n-1}} \right)^n \left(\int_{\mathbb A} \frac{dx}{|x|^n}\right)^{1-n} \nonumber \\
&\geqslant & \left(\mbox{Mod}\, {\mathbb A}^\ast \right)^n
\left(\mbox{Mod}\, {\mathbb A} \right)^{n-1} =\alpha^n \,
\mbox{Mod}\, {\mathbb A}
\end{eqnarray}
where
\begin{equation}
\alpha \bydef \frac{\mbox{Mod}\, {\mathbb A}^\ast}{\mbox{Mod}\,
{\mathbb A}}
\end{equation}
Similarly, we obtain
\begin{eqnarray}
 \int_{\mathbb A} \frac{|h_T|^n}{|h|^n} &\geqslant &  \left(  \int_{\mathbb A} \frac{|h_T|^{n-1}}{|h|^{n-1}|x|} \right)^\frac{n}{n-1} \left(\int_{\mathbb A} \frac{dx}{|x|^n}\right)^\frac{1}{1-n} \nonumber \\
&\geqslant & \left(\mbox{Mod}\, {\mathbb A}\right)^\frac{n}{1-n}
\left(\mbox{Mod}\, {\mathbb A} \right)^\frac{1}{n-1} =\mbox{Mod}\,
{\mathbb A}
\end{eqnarray}
Substituting these two estimates into (\ref{E428}) we conclude with
the following Theorem
\begin{theorem}\label{th4}
Let ${\mathbb A}$ and ${\mathbb A}^\ast$ be spherical rings in
$\rn$, $n\geqslant 2$. Then for every $h\in {\mathscr P}({\mathbb
A}, {\mathbb A}^\ast)$ we have
\begin{equation}\label{IE11}
\int_{{\mathbb A}} \frac{| Dh |^n}{|h|^n} \geqslant   \max\left\{1,
\alpha^n\right\}  \textnormal{Mod}\, {\mathbb A}
\end{equation}
As for the sharpness of this estimate we note that equality holds
for the power stretching $h (x) =|x|^{\alpha -1} x$.
\end{theorem}
There are, however, other cases of equality in~\eqref{IE11} if
$\mbox{Mod}\, {\mathbb A}^\ast \neq \mbox{Mod}\, {\mathbb A}$.
\section{Radial symmetry}\label{Sec52}
Suppose that the radial stretching
\begin{equation}
h(x)=H(|x|)\, \frac{x}{|x|}\, , \; \; \; H\in \mathscr
A\mathcal{C}\, [r,R]
\end{equation}
maps homeomorphically  an annulus ${\mathbb A}=\{x \, ; \; \; r <
|x| < R \}$ onto  ${\mathbb A}^\ast=\{y \, ; \; \; r_\ast <|y| <
R_\ast \}$, where $\mbox{Mod}\, {\mathbb A}^\ast = \alpha\,
\mbox{Mod}\, {\mathbb A}$, with some $\alpha
>0$. We may  assume that $r_\ast = r^\alpha$ and $R_\ast =
R^\alpha$. Thus, in particular, $h$ preserves the order of the
boundary components. The power stretching $h^\alpha = |x|^{\alpha
-1} x$ serves as an example of such a homeomorphism. We shall prove
that $h^\alpha$ is in fact the minimizer of $\mathscr F_h$
among all radial stretchings.
\begin{proposition}\label{Pr41}
For each radial stretching we have
\begin{equation}\label{E411}
{\mathscr F}_h= \int_{\mathbb A} \frac{\norm Dh \norm^n}{|h|^n}
\geqslant \int_{\mathbb A} \frac{\norm Dh^\alpha \norm^n}{|h^\alpha
|^n} = \left(\alpha^2 + n-1\right)^\frac{n}{2} \textnormal{Mod}\, {\mathbb
A}
\end{equation}
Equality holds modulo rotation  of ${\mathbb A}$.
\end{proposition}
\proof We express the energy of $h$ as
\begin{equation}
{\mathscr F}_h= \int_{\mathbb A} \frac{\left(\left|h_N\right|^2 +
(n-1) \left|h_T \right|^2 \right)^\frac{n}{2}}{|h|^n}
\end{equation}
where the radial and spherical derivatives of $h$ are found in terms
of the strain function $H$ by using formulas  (\ref{E147});
\begin{equation}
\left|h_N \right| = \left|\dot{H}(t) \right| \; \; \; \; \mbox{ and }\; \; \; \; \left|h_T \right|= \frac{\left|H(t)\right|}{t}
\end{equation}
Integration in  polar coordinates leads to a line integral
\begin{eqnarray}
{\mathscr F}_h &=& \omega_{n-1} \int_r^R \left(\frac{t^2 \dot{H}^2}{H^2}+n-1\right)^\frac{n}{2}\, \frac{dt}{t} \nonumber \\
&=&  \omega_{n-1} \log \frac{R}{r} \Aint \left(X^2 +n-1\right)^\frac{n}{2}\, d \mu
\end{eqnarray}
where  the integral average is taken with respect to the measure
$d\mu = \frac{dt}{t}$. The key is that $X(t)= \frac{t \dot H}{H}$ is
a one dimensional free Lagrangian. Indeed, we see that the integral mean
\begin{equation}
\Aint_r^R X \, d\mu = \frac{1}{\log \frac{R}{r}}\int_r^R
\frac{\dot{H}(t)\, dt}{H(t)} = \frac{\log \frac{R_\ast
}{r_\ast}}{\log \frac{R}{r}}= \alpha
\end{equation}
does not depend on $H$. Next using Jensen's inequality for the convex function $X \rightarrow (X^2 + n-1)^\frac{n}{2}$ yields
\begin{eqnarray}
{\mathscr F}_h&\geqslant &  \left[\left(\Aint_r^R X\, d\mu \right)^2 + n-1\right]^\frac{n}{2}\mbox{Mod}\, {\mathbb A} \nonumber \\
&=&    \left[\alpha^2 + n-1\right]^\frac{n}{2} \mbox{Mod}\, {\mathbb
A}
\end{eqnarray}
as desired. Examining these arguments backwards we obtain
the uniqueness statement. \vskip0.3cm

As a corollary, letting $\alpha$ go to $0$, we observe that among all radial deformations $h: {\mathbb A} \rightarrow {\mathbb S}^{n-1}$, the $n$-harminic energy
\begin{equation}
{\mathscr E}_h = \int_{\mathbb A} \norm Dh \norm^n \geqslant
(n-1)^\frac{n}{2} \mbox{Mod}\, {\mathbb A}
\end{equation}
assumes its minimal value for the mapping $h(x)= \frac{x}{|x|}$.
This will be later shown to be true in a larger class of mappings.

In dimensions $n=2,3$ the assertion of Proposition \ref{Pr41}
remains true for all permissible mappings in the Sobolev class  $
{\mathscr W}^{1,n}({\mathbb A}, {\mathbb A}^\ast)$.\footnote{In the
next chapter we shall extend this result to all dimensions, but it
will require an upper bound of the modulus of ${\mathbb A}^\ast$.}

\section[Proof of Theorem 1.14]{Proof of Theorem \ref{th3}}
The proof in case $n=2$ is accomplished by an elegant use of
H\"older's inequality,
\begin{eqnarray}
{\mathscr F}_h&= &\int_{\mathbb A} \frac{\norm Dh \norm^2}{|h|^2} = \int_{\mathbb A} \left(\left|\frac{h_N}{h}\right|^2 +\left|\frac{h_T}{h}\right|^2 \right) \nonumber \\
&\geqslant & \left[\left(\int_{\mathbb A} \frac{\left|h_N\, \right|}{|x|\, |h|}\right)^2 + \left(\int_{\mathbb A} \frac{\left|h_T\, \right|}{|x|\, |h|}\right)^2\right] \left(\int_{\mathbb A} \frac{dx}{|x|^2}\right)^{-1}
\end{eqnarray}
The  lower bounds at (\ref{E185}) and (\ref{E186}) yield
\[
{\mathscr F}_h \geqslant \left[\left(\mbox{Mod}\, {\mathbb A}^\ast
\right)^2 + \left(\mbox{Mod}\, {\mathbb A}
\right)^2\right]\left(\mbox{Mod}\, {\mathbb A}  \right)^{-1} =
\left(\alpha ^2 +1\right)\, \mbox{Mod}\, {\mathbb A}  = {\mathscr
F}_{h^\alpha }
\]
as claimed.

The proof in dimension $n=3$ involves application of  H\"older's inequality
twice
\begin{eqnarray}
{\mathscr F}_h &=& \int_{\mathbb A} \frac{\norm Dh \norm^3}{|h|^3} = \int_{\mathbb A} \left(\left|\frac{h_N}{h}\right|^2 +2 \left|\frac{h_T}{h}\right|^2 \right)^\frac{3}{2} \nonumber \\
&\geqslant & \left[\int_{\mathbb A} \left( \left|\frac{h_N}{h}\right|^2 +2 \left|\frac{h_T}{h}\right|^2 \right)\frac{dx}{|x|} \right]^\frac{3}{2} \left(\int_{\mathbb A} \frac{dx}{|x|^3}\right)^{-\frac{1}{2}} \nonumber \\
&\geqslant & \left[\int_{\mathbb A} \left( \left|\frac{h_N}{h}\right| \frac{dx}{|x|^2}\right)^2  \left(\int_{\mathbb A} \frac{dx}{|x|^3}\right)^{-1}  +2   \int_{\mathbb A}  \left|\frac{h_T}{h}\right|^2 \frac{dx}{|x|} \right]^\frac{3}{2} \left(\int_{\mathbb A} \frac{dx}{|x|^3}\right)^{-\frac{1}{2}}\nonumber
\end{eqnarray}
Now using  lower bounds at  (\ref{E185}) and (\ref{E186}) we find that
\begin{eqnarray}
{\mathscr F}_h &\geqslant & \left[\left(\mbox{Mod}\, {\mathbb
A}^\ast \right)^2 \left(\mbox{Mod}\, {\mathbb A} \right)^{-1}+ 2\,
\mbox{Mod}\, {\mathbb A}\right]^\frac{3}{2}\left(\mbox{Mod}\,
{\mathbb A}  \right)^{-\frac{1}{2}} \nonumber \\ &=&\left(\alpha^2
+2\right)^\frac{3}{2}\,  \mbox{Mod}\, {\mathbb A}  =
\mathscr{F}_{h^\alpha }
\end{eqnarray}
For the uniqueness we refer to Section \ref{SecUniSpe} where such
problems are dealt  in all dimensions. The lower bounds of free Lagrangians,  for
the energy integrand $\frac{\norm Dh \norm^n}{|h|^n}$ in higher
dimensions,  are more sophisticated. We found optimal bounds if
$\mbox{Mod}\, {\mathbb A}^\ast \leqslant \alpha_n \mbox{Mod}\,
{\mathbb A} $, where $1< \alpha_n < \infty$ is a solution to the
algebraic equation
\begin{equation}\label{E42}
\left(\alpha^2_n+n-1\right)^{n-2} \left(\alpha_n^2-1\right)^2
=\alpha_n^{2n} \; \Rightarrow \; \alpha_n < \sqrt{\frac{n-1}{n-3}}.
\end{equation}

\chapter[Proof of Theorem 1.15]{Proof of Theorem \ref{th3half}}\label{Sec6}
It suffices to consider homeomorphisms $h: {\mathbb A} \rightarrow
{\mathbb A}^\ast$ of class ${\mathscr W}^{1,n}( {\mathbb A},
{\mathbb A}^\ast )$. Recall that the target annulus is not too fat.
Precisely we  assume that
\begin{equation}
\mbox{Mod}\, {\mathbb A}^\ast < \alpha_n\,  \mbox{Mod}\,
{\mathbb A}
\end{equation}
where the critical factor $\alpha_n >1$ is determined by the
equation  (\ref{E42}), which we write as
\begin{equation}\label{E222st}
\frac{(\alpha_n^2+n-1)^\frac{n-2}{2} (\alpha_n^2-1)}{\alpha_n^n}=1
\end{equation}
\begin{lemma}\label{le7}
Let $X,Y \geqslant 0$ and $1 \leqslant \alpha < \alpha_n$. Then
\begin{equation}\label{ES208}
a=a(\alpha)\bydef \frac{(\alpha^2+n-1)^\frac{n-2}{2}
(\alpha^2-1)}{\alpha^n}<1
\end{equation}
and, we have
\begin{equation}\label{E212}
\left[X^2 +(n-1)Y^2\right]^\frac{n}{2} \geqslant a\,  X^n + b\,  X
Y^{n-1}
\end{equation}
where
\begin{equation}\label{defb}
b= \frac{n\left(\alpha^2 + n-1\right)^\frac{n-2}{2}}{\alpha}
\end{equation}
Equality holds if and only if $X = \alpha Y$.
\end{lemma}
\proof Because of homogeneity we may assume that $Y=1$; the case
$Y=0$ is obvious. We are reduced to proving an  inequality with one
variable $X\geqslant 0$,
\begin{equation}\label{228}
\left[X^2 +(n-1)\right]^\frac{n}{2} \geqslant a X^n +bX
\end{equation}
Consider the function
\begin{equation}
\varphi (X) = \left(X^2 + n-1\right)^\frac{n}{2} - aX^n -bX
\end{equation}
and its two derivatives
\begin{eqnarray}
\varphi^\prime (X) &=& n \left(X^2+n-1\right)^\frac{n-2}{2}X-naX^{n-1}-b \\
\varphi^{\prime \prime} (X) &=& n(n-1) \left[\left(X^2+n-1\right)^\frac{n-4}{2}(X^2+1)-aX^{n-2}\right]>0
\end{eqnarray}
This letter estimate is guaranteed by our hypothesis that $a<1$.
Thus $\varphi$ is strictly convex and, therefore, has at most one
critical point. The coefficient $b$ has been defined above exactly
in  a way to ensure that $\varphi^\prime (\alpha )=0$. The
coefficient $a$ has been chosen so that  $\varphi (\alpha)=0$. This
completes the proof of Lemma \ref{le7}. \vskip0.5cm

\begin{lemma}\label{le5}
Let $X,Y \geqslant 0$ and $0 \leqslant  \alpha \leqslant 1$. Then
\begin{equation}
\left[X^2 +(n-1)Y^2\right]^\frac{n}{2} \geqslant a\,  Y^n +b\, X
Y^{n-1}
\end{equation}
where
\begin{equation}
a= (n-1) \left(\alpha^2 +n-1\right)^\frac{n-2}{2} \left(1-
\alpha^2\right)
\end{equation}
and
\begin{equation}\label{eqb}
b=  n \alpha \left(\alpha^2 +n-1 \right)^\frac{n-2}{2}
\end{equation}
Equality holds if and only if $X = \alpha Y$.
\end{lemma}
\proof Because of homogeneity we may assume that $Y=1$, the case
$Y=0$ is obvious. We are reduced to the inequality with one
variable $X\geqslant 0$,
\begin{equation}\label{239}
\left[X^2 +(n-1)\right]^\frac{n}{2} \geqslant a  +bX
\end{equation}
Consider the function
\begin{equation}
\varphi (X) = \left(X^2 + n-1\right)^\frac{n}{2} -bX
\end{equation}
and its two derivatives
\begin{eqnarray}
\varphi^\prime (X) &=& n \left(X^2+n-1\right)^\frac{n-2}{2}X-b \\
\varphi^{\prime \prime} (X) &=& n(n-1) \left(X^2+n-1\right)^\frac{n-4}{2}(X^2+1)>0
\end{eqnarray}
As before, $\varphi$ is strictly convex with one critical point. The
coefficient $b$ has been defined exactly to ensure that
$\varphi^\prime (\alpha)=0$, whereas  $a=\varphi (\alpha)$. This
completes the proof of Lemma \ref{le5}.


We now proceed to the proof of Theorem \ref{th3half}. Consider two cases.
\section{The  case of expanding pair} We  apply Lemma
\ref{le7} with $X= \frac{\left|h_N\right|}{|h|}$, $Y=
\frac{\left|h_T\right|}{|h|}$ and $\alpha= \frac{\textnormal{Mod}\,
{\mathbb A}^\ast }{\textnormal{Mod}\, {\mathbb A}} \geqslant 1$, to
obtain the point-wise inequality
\begin{eqnarray}\label{487}
\frac{\norm Dh \norm^n}{|h|^n} &=& \left( \frac{\left|h_N\right|^2}{|h|^2} +(n-1) \frac{\left|h_T\right|^2}{|h|^2} \right)^\frac{n}{2} \nonumber \\
&\geqslant & a  \frac{\left|h_N\right|^n}{|h|^n} + b \frac{\left|h_N\right| \, \left|h_T\right|^{n-1}}{|h|^n}
\end{eqnarray}
Then we integrate (\ref{487}) over the annulus ${\mathbb A}$. For
the last term  we may apply the  lower bound at
(\ref{E184}). The first term in the right hand side of (\ref{487})
needs some adjustments before using the lower bound at (\ref{E185}).
These adjustments are easily accomplished by H\"older's inequality.
\begin{eqnarray}
\int_{\mathbb A} \frac{\norm Dh \norm^n}{|h|^n} &\geqslant & a \left(\int_{\mathbb A} \frac{\left|h_N\right|}{|h|\, |x|^{n-1}}\right)^n  \left(\int_{\mathbb A} \frac{dx}{|x|^n}\right)^{1-n} + b\,  \mbox{Mod}\, {\mathbb A}^\ast \nonumber \\
&\geqslant & a   \left(\mbox{Mod}\, {\mathbb A}^\ast \right)^n  \left(\mbox{Mod}\, {\mathbb A}\right)^{1-n} +  b \, \mbox{Mod}\, {\mathbb A}^\ast \nonumber \\
&=&  \left(a \alpha^n + b\alpha  \right)  \mbox{Mod} \, {\mathbb A} \nonumber \\
&=& \left(\alpha^2 + n-1 \right)^\frac{n}{2}\, \mbox{Mod}\, {\mathbb
A}
\end{eqnarray}
because we have equality  at (\ref{228}) for $X=\alpha$. We
see that the right hand side is the energy of the power stretching
\begin{equation}
h^\alpha (x)= r^\ast r^{-\alpha} |x|^{\alpha -1}x \, , \hskip2cm
\alpha = \frac{\mbox{Mod}\, {\mathbb A}^\ast}{\mbox{Mod}\, {\mathbb
A}}
\end{equation}
\section{The  case of contracting pair}
This time we use Lemma \ref{le5} with the same data $\,X, Y, \alpha\,$ as were used in the case of expanding pair; that is,
with $X= \frac{\left|h_N\right|}{|h|}$, $Y=
\frac{\left|h_T\right|}{|h|}$ and $\alpha= \frac{\textnormal{Mod}\,
{\mathbb A}^\ast }{\textnormal{Mod}\, {\mathbb A}} \leqslant 1$. This gives the following point wise inequality
\begin{eqnarray}\label{488}
\frac{\norm Dh \norm^n}{|h|^n} &=& \left( \frac{\left|h_N\right|^2}{|h|^2} +(n-1) \frac{\left|h_T\right|^2}{|h|^2} \right)^\frac{n}{2} \nonumber \\
&\geqslant & a \, \frac{\left|h_T\right|^n}{|h|^n} + b\,
\frac{\left|h_N\right| \, \left|h_T\right|^{n-1}}{|h|^n}
\end{eqnarray}
We integrate it over the annulus ${\mathbb A}$. For the last term we
apply the lower bound at (\ref{E184}). As in previous case the  first term in
the right hand side of (\ref{488}) needs some adjustments before
using the lower bound at (\ref{E186}). These adjustments are easily
accomplished by H\"older's inequality.
\begin{eqnarray}
\int_{\mathbb A}\frac{\norm Dh \norm^n}{|h|^n} &\geqslant & a \left(\int_{\mathbb A} \frac{\left|h_T\right|^{n-1}}{|x|\, |h|^{n-1}}\right)^\frac{n}{n-1}  \left(\int_{\mathbb A} \frac{dx}{|x|^n}\right)^\frac{-1}{n-1} + b\, \mbox{Mod}\, {\mathbb A}^\ast \nonumber \\
&\geqslant &  a\left(\mbox{Mod}\, {\mathbb A} \right)^\frac{n}{n-1}   \left(\mbox{Mod}\, {\mathbb A}\right)^\frac{-1}{n-1} +  b\,  \mbox{Mod}\, {\mathbb A}^\ast  \nonumber \\
&=&   \left(a  + b\,\alpha  \right)\, \mbox{Mod}\, {\mathbb A} \nonumber \\
&=&\left(\alpha^2 + n-1 \right)^\frac{n}{2}\mbox{Mod}\, {\mathbb A}
\end{eqnarray}
because we have equality in  (\ref{239}) for $X = \alpha$. We easily
see that the right hand side is the energy of the power stretching
\begin{equation}
h^\alpha (x)= r^\ast r^{-\alpha} |x|^{\alpha -1}x \, , \hskip2cm
\alpha = \frac{\mbox{Mod}\, {\mathbb A}^\ast}{\mbox{Mod}\, {\mathbb
A}}
\end{equation}
\section{Uniqueness}\label{SecUniSpe}
The reader might want to compare this proof with Chapter~\ref{SecUni}, in which the uniqueness question is  dealt in greater generality. Let $h : {\mathbb A} \rightarrow {\mathbb
A}^\ast$ be any extremal mapping. First observe that in both cases
we have used a general inequality
\begin{equation}
J(x,h) \leqslant \left|h_N\right| \, \left|h_T\right|^{n-1}
\end{equation}
For $h$ to be extremal, we  must have equality. In view of
(\ref{145}) the vectors $h_N, h_{T_2}, ..., h_{T_n}$ are mutually
orthogonal and $\left| h_{T_2} \right|= ... = \left| h_{T_n} \right|
= \left| h_{T} \right|$. Using the  matrix representation of $Dh$ at
(\ref{311}) we find the Cauchy-Green tensor of $h$ to be a diagonal
matrix
\begin{equation}
D^\ast h \; Dh =\left[ \begin{array}{cccc} \left| h_{N} \right|^2 & 0 & \cdots & 0  \\
 0 & \left| h_{T} \right|^2 &   \cdots & 0\\
& &\ddots &  \\
0 & 0 & \cdots & \left| h_{T} \right|^2
\end{array}\right]
\end{equation}
Another necessary condition for $h$ to be extremal is that $\left|
h_{N} \right| = \alpha \left| h_{T} \right|$, because of the
equality cases in Lemmas \ref{le7} and \ref{le5}. In this way,
we arrive at the Beltrami type system for the extremal mapping
\begin{equation}
D^\ast h \; Dh = J(x,h)^\frac{2}{n} {\bold K}
\end{equation}
where ${\bold K}$ is a constant diagonal matrix
\begin{equation}
{\bold K}=\left[ \begin{array}{cccc} \alpha^{1-\frac{1}{n}} & 0 & \cdots & 0  \\
 0 & \alpha^{-\frac{1}{n}} &   \cdots & 0\\
& &\ddots &  \\
0 & 0 & \cdots &  \alpha^{-\frac{1}{n}}
\end{array}\right]
\end{equation}
The power mapping $h^\alpha = h^\alpha (x)$ is one of the
homeomorphic solutions to this system. It is well known and easy to verify that other
solutions are obtained by composing this particular one with a
conformal transformation \cite{BI, Reb, Rib}. Thus $h$
takes the form
\begin{equation}
h=g\circ h^\alpha
\end{equation}
where $g : {\mathbb A}^\ast \rightarrow {\mathbb A}^\ast$ is a
conformal automorphism of the target  annulus ${\mathbb A}^\ast = A(r_\ast, \R_\ast)$ onto
itself. Thus, up to rotation $g(y)=y$ or $g(y)=r_\ast R_\ast \frac{y}{\abs{y}^2}$. In either case $h$ is a power stretching up to the rotation.  The proof of
uniqueness is complete.

\part{The $n$-Harmonic Energy}

A study of the extremal problems for the conformal energy
\begin{equation}
{\mathscr E}_h= \int_{\mathbb A} \norm Dh (x) \norm^n\, \d x
\end{equation}
is far more involved. Even when the solutions turn out to be radial
stretchings, they no longer represent elementary functions such as
$h^\alpha = |x|^{\alpha -1}x$. And this is not the only difficulty;
there are  new phenomena in case of non-injective solutions.  In view of these concerns, it is  remarkable that the
method of free Lagrangians is still effective  for non-injective solutions. First we examine planar mappings.

 Concerning  mappings minimizing the $\mathscr L^p$-norm of the gradient, we referee the interested reader to \cite{HL} and references there.

\chapter[Harmonic mappings between planar annuli]{Harmonic mappings between planar annuli,\\ proof of Theorem~\ref{thn=2}}\label{secn2}

Let $h : {\mathbb A} \rightarrow {\mathbb A}^\ast$ be a
homeomorphism between annuli in the Sobolev class ${\mathscr
W}^{1,2} ({\mathbb A}, {\mathbb A}^\ast)$. We view $h$ as a complex
valued function. Let us recall the formulas
\begin{equation}
\norm Dh \norm^2 = \left|h_N \right|^2 + \left|h_T \right|^2
\end{equation}
and
\begin{equation}
\mbox{det}\, Dh = \mbox{Im}\, \left(h_T\, \overline{h_N}\right) \leqslant  \left|h_N \right|\, \left|h_T \right|
\end{equation}
{\bf {\it Case 1.} The expanding pair, $\frac{R}{r}\leqslant
\frac{R^\ast}{r^\ast}$.}

Thus the target annulus is
conformally fatter than the domain. We find a unique number $\omega
\leqslant 0$ such that
\begin{equation}\label{51star}
\frac{R}{r} = \frac{R_\ast+\sqrt{R_\ast^2-\omega}}{r_\ast+
\sqrt{r_\ast^2-\omega}}
\end{equation}
Without loss of generality we may assume, by rescaling the annulus $\mathbb A$ if necessary,   that
\begin{equation}\label{51dstars}
R = R_\ast + \sqrt{R_\ast^2 -\omega} \; \; \mbox { and } \; \; r =
r_\ast + \sqrt{r_\ast^2 - \omega}
\end{equation}
We begin with the inequality
\begin{equation}\label{515}
\left(\frac{|h|\,  |h_N|}{\sqrt{|h|^2-\omega}} -   \left|h_T \right|\right)^2 \geqslant 0
\end{equation}
Equivalently,
\begin{eqnarray}
\norm Dh \norm^2 &\geqslant & \frac{- \omega}{|h|^2 - \omega}  \left|h_N \right|^2 + \frac{2\, |h|}{ \sqrt{|h|^2 -\omega}}  \left|h_N \right|\,  \left|h_T \right| \nonumber \\
&\geqslant &  \frac{- \omega}{|h|^2 - \omega}  \left(\, |h|_N\,  \right)^2 + \frac{2\, |h|}{ \sqrt{|h|^2 -\omega}}\, \left|h_N \right|\,  \left|h_T \right| \nonumber \\
&=& - \omega \left\{ \left[\log \left(|h|+ \sqrt{|h|^2 - \omega}\right)\right]_N\right\}^2   + \frac{2\, |h|}{ \sqrt{|h|^2 -\omega}}\,  \left|h_N \right|\,  \left|h_T \right| \nonumber
\end{eqnarray}
Here we have used an elementary fact that $\left|\, h_N\,  \right|
\geqslant \left|\, |h|_N\,  \right|$, equality occurs if and only if
$\frac{h_N}{h}$ is a real valued function. Let us integrate this
estimate. We apply H\"older's inequality to the first term and the
estimate (\ref{E184}) to the second term,
\begin{eqnarray}
\int_{\mathbb A} \norm Dh \norm^2 &\geqslant&  - \omega \left[\int_{\mathbb A} \left[\log \left(|h|+ \sqrt{|h|^2 - \omega}\right)\right]_N \, \frac{dx}{|x|}\right]^2 \cdot \left(\int_{\mathbb A} \frac{dx}{|x|^2}\right)^{-1}   + \nonumber \\ &\; & +\, 4\pi \int_{r_\ast}^{R_\ast} \frac{\tau^2}{ \sqrt{\tau^2 -\omega}}\, d \tau \nonumber \\
&=& - \omega \left|2\pi \log \frac{R_\ast+\sqrt{R_\ast^2-\omega}}{r_\ast+ \sqrt{r_\ast^2-\omega}}\right|^2 \left(2\pi \log \frac{R}{r}\right)^{-1} + \nonumber \\
&\; &  +\, 2\pi \bigg[ \tau \sqrt{\tau^2 - \omega} + \omega \log (\tau+ \sqrt{\tau^2 -\omega})\bigg]_{\tau=r_\ast}^{R_\ast} \nonumber \\
&=& 2 \pi R_\ast \sqrt{R_\ast^2-\omega}\,  -\,   2 \pi r_\ast
\sqrt{r_\ast^2-\omega}
\end{eqnarray}
Elementary inspection reveals that equality holds for the Nitsche
mapping $$h(z)= \frac{1}{2} \left( z + \frac{\omega}{
\bar{z}}\right)$$ {\bf {\it Case 2.} The contracting pair, $ \frac{R_\ast}{r_\ast}
< \frac{R}{r} \leqslant  \frac{R_\ast}{r_\ast} +
\sqrt{\frac{R_\ast^2}{r_\ast^2}-1} $.}

 In this case the target
annulus is conformally thinner than the domain, but not too thin. We
express this condition by using the equation (\ref{51star}), where
this time $0< \omega \leqslant r_\ast^2$. Again we may assume that
the relations at (\ref{51dstars}) hold. The same inequality
(\ref{515}) can be expressed in a somewhat different form
\begin{equation}
\norm Dh \norm^2 \geqslant  \omega  \left|\frac{h_T }{h}\right|^2 +
2 \left|h_N \right|\,  \left|h_T \right| \sqrt{1-
\frac{\omega}{|h|^2}}
\end{equation}
As before, we   apply H\"older's inequality, which together with the
estimate (\ref{E184}) yields
\[
\int_{\mathbb A} \norm Dh \norm^2 \geqslant \omega
\left(\int_{\mathbb A} \frac{\left|h_T \right|}{|x|\, |h|}\right)^2
\cdot \left(\int_{\mathbb A} \frac{\d x}{|x|^2}\right)^{-1}   + 4
\pi \int_{r_\ast}^{R_\ast} \sqrt{\tau^2- \omega}\, d \tau
\]
Next, the estimate (\ref{E186}) gives
\begin{eqnarray}\label{5117}
\int_{\mathbb A} \norm Dh \norm^2 &\geqslant & \omega \left(2 \pi  \log \frac{R}{r}\right)^2 \cdot \left(2 \pi  \log \frac{R}{r}\right)^{-1} +  4 \pi \int_{r_\ast}^{R_\ast} \sqrt{\tau- \omega}\, d \tau \nonumber \\
&=&  2 \pi \omega \log \frac{R}{r} + 2 \pi \bigg[\tau \sqrt{\tau^2 -
\omega }- \omega \log (\tau + \sqrt{\tau^2 - \omega})\bigg]_{\tau =
r_\ast}^{R_\ast}\; \; \;
\end{eqnarray}
In view of (\ref{51dstars}) we find that
\begin{equation}
\int_{\mathbb A} \norm Dh \norm^2 \geqslant 2 \pi R_\ast
\sqrt{R_\ast^2-\omega}\, -\, 2 \pi r_\ast \sqrt{r_\ast^2-\omega}
\end{equation}
Again equality holds for the Nitsche mapping $$h(z)= \frac{1}{2}
\left( z + \frac{\omega}{ \bar{z}}\right)$$ {\bf The borderline
case.} Taking $\omega =r_\ast^2$  we obtain what is called the
critical Nitsche map with $R=R_\ast + \sqrt{R_\ast^2-r_\ast^2}$ and
$r=r_\ast$,
\begin{equation}
h^\ast (z) = \frac{1}{2} \left(z + \frac{r_\ast^2}{\bar{z}}\right)
\end{equation}

\begin{figure}[!h]
\begin{center}
\psfrag{r=}{\tiny ${r=r_\ast}$} \psfrag{R=}{\tiny ${R= R_\ast +
\sqrt{R_\ast^2 - r_\ast^2} }$} \psfrag{Rp}{\tiny ${ R_\ast }$}
\psfrag{hsha}{\small ${ h^\ast }$} \psfrag{rp}{\tiny ${ r_\ast }$}
\includegraphics*[height=1.2in]{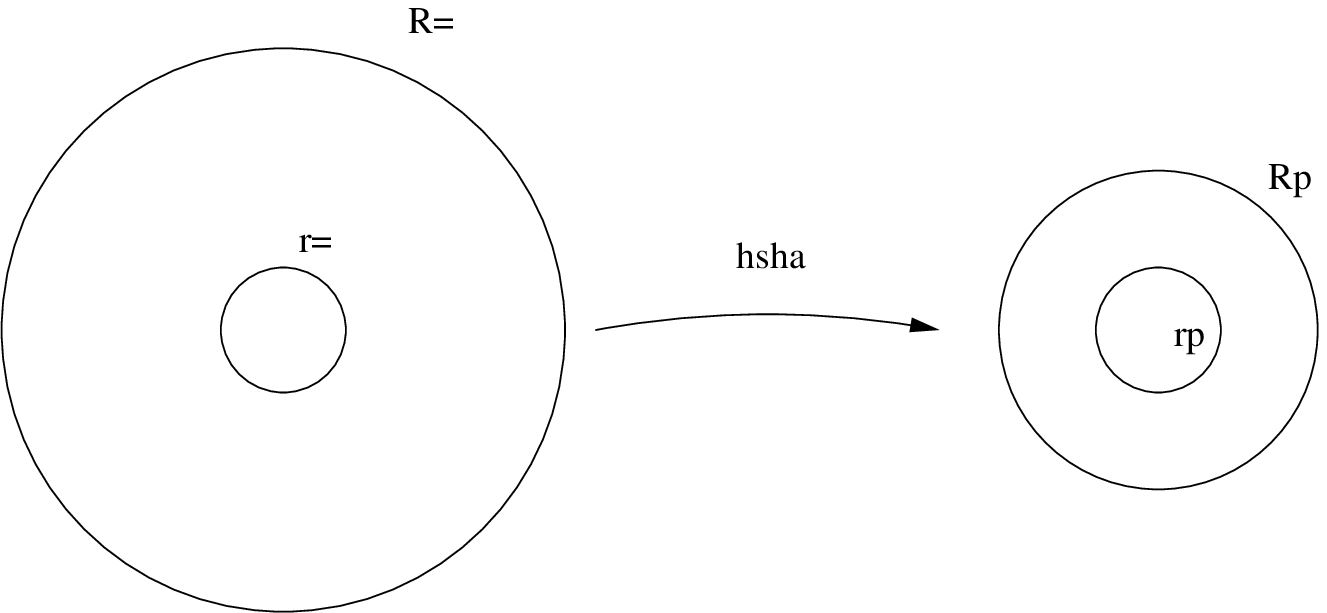}
\caption{The Jacobian of $h^\ast$ vanishes on the inner circle.}
\end{center}
\end{figure}

The conformal energy of the critical Nitsche map equals:
\begin{equation}\label{equ288}
\int_{\mathbb A} \norm Dh^\ast \norm^2 = 2 \pi  R_\ast
\sqrt{R_\ast^2-r_\ast^2}
\end{equation}
We are now in a position to consider the case:\\
{\bf {\it Case 3.} Below the lower Nitsche bound,  $ \frac{R_\ast}{r_\ast} +
\sqrt{\frac{R_\ast^2}{r_\ast^2}-1} < \frac{R}{r}$.}

 The target
annulus ${\mathbb A}^\ast$ is too thin. We shall see that an inner part of $\mathbb A$
has to be hammered flat to the inner circle of ${\mathbb A}^\ast$.
We can certainly rescale the annulus $\mathbb A$ to have $R=R_\ast + \sqrt{R_\ast^2-r_\ast^2}$. This together with the hypothesis of this case, yields $r<
r_\ast$.

\begin{figure}[!h]
\begin{center}
\psfrag{r}{\tiny ${r}$} \psfrag{R=}{\tiny ${R= R_\ast +
\sqrt{R_\ast^2 - r_\ast^2} }$} \psfrag{Rp}{\tiny ${ R_\ast }$}
\psfrag{h}{\small ${ h }$} \psfrag{A}{\small ${\mathbb A}$}
\psfrag{Ap}{\small ${\mathbb A}^\ast$} \psfrag{rp}{\tiny ${ r_\ast
}$}
\includegraphics*[height=1.4in]{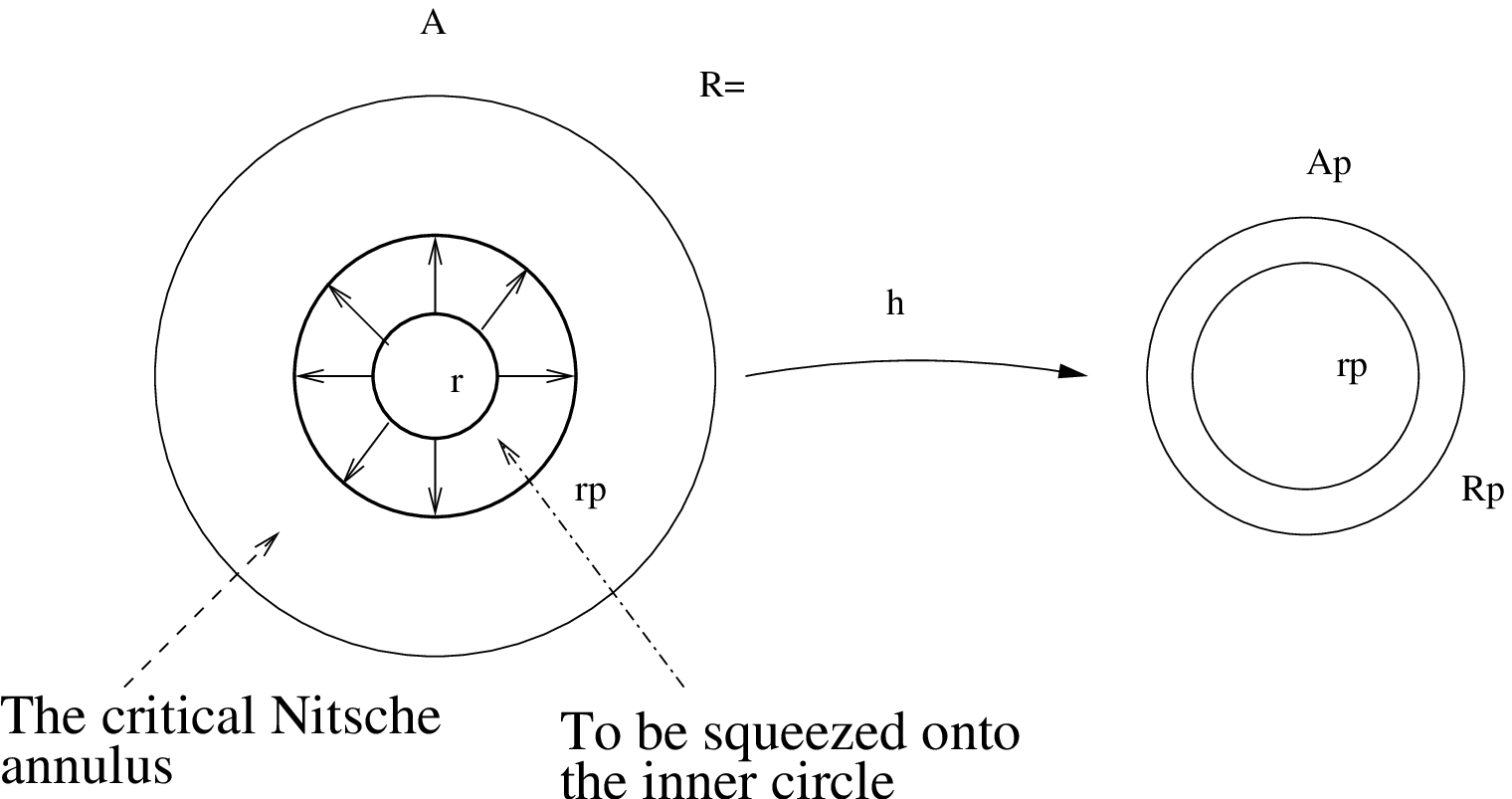}
\caption{Below the lower Nitsche bound.}
\end{center}
\end{figure}

For every permissible map $h : {\mathbb A} \rightarrow {\mathbb
A}^\ast$ we still have the estimate (\ref{5117}), in which we let $\omega =
r_\ast^2$. Hence
\begin{eqnarray}
\int_{\mathbb A} \norm Dh \norm^2 &\geqslant &   2 \pi r_\ast \log \frac{R}{r} + 2 \pi \bigg[\tau \sqrt{\tau^2 - r_\ast }- r_\ast \log (\tau + \sqrt{\tau^2 - r_\ast})\bigg]_{\tau = r_\ast}^{R_\ast} \nonumber \\
&=& 2 \pi R_\ast \sqrt{R_\ast^2-r_\ast^2}\,  +\,   2 \pi r_\ast^2 \log \frac{r_\ast}{r}
\end{eqnarray}
The first term represents conformal energy of $h^\ast \colon \mathbb A (r_\ast, R) \to \mathbb A(r_\ast, R_\ast)$, by~\eqref{equ288}. The second term turns out to be exactly the energy of the hammering mapping  $g(z)= r_\ast \frac{z}{\abs{z}}$, which takes the remaining part $\mathbb A(r,r_\ast) \subset \mathbb A$ onto the inner circle of $\mathbb A(r_\ast, R_\ast)$.
\begin{equation}
{\mathscr E}_g= \int_{r<|z|<r_\ast} \norm Dg \norm^2 = 2 \pi
r_\ast^2 \log \frac{r_\ast}{r}.
\end{equation}
Thus the energy of the mapping $h^\circ \colon \mathbb A \onto \mathbb A^\ast$
\[h^\circ = \begin{cases} r_\ast \frac{z}{\abs{z}}, & \quad r< \abs{z} \le r_\ast \\
\frac{1}{2} \left(z + \frac{r_\ast^2}{\bar z}\right), & \quad r_\ast < \abs{z}< R  \end{cases}\]
is smaller than the energy of any homeomorphism $h \colon \mathbb A (r,R) \onto \mathbb A (r_\ast, R_\ast)$. It is clear that $h^\circ$ is a $\W^{1,2}$-limit of such homeomorphisms, completing the proof of Theorem~\ref{thn=2}.
 \vskip0.3cm
 Now we
proceed to higher dimensions. Different estimates will be required
for the contracting pairs of annuli than for the expanding pairs. Thus
we devote a separate section for each case.
\chapter{Contracting pair, $ \mbox{Mod}\, {\mathbb A}^{\! \ast} \leqslant \mbox{Mod}\, {\mathbb
A}$}\label{Sec22}
First we consider the case when  $\mathbb A^{\! \ast}$ is not too
thin.
\section[Proof of Theorem 1.12]{Proof of Theorem \ref{ThNe3}}
Here we assume the Nitsche bound, ${\mathcal N}_\dag
(\mbox{Mod}\, {\mathbb A})\leqslant \mbox{Mod}\, {\mathbb A}^{\!
\ast}$, see~\eqref{equ34p} for ${\mathcal N}_\dag$. This bound means precisely that there is a radial
$n$-harmonic homeomorphism
\begin{equation}\label{E312}
h^\circ : {\mathbb A} \rightarrow {\mathbb A}^{\! \ast}\, ,
\hskip1cm h^\circ (x) = H\big(|x|\big)\, \frac{x}{|x|}
\end{equation}
Recall the characterictic equation at (\ref{E231}) for $H=H(t)$
\begin{equation}\label{E954376}
\left[H^2 + \frac{t^2 \dot{H}^2}{n-1}  \right]^\frac{n-2}{2} \left(H^2 - t^2 \dot{H}^2\right) \equiv c
\end{equation}
where $c$ is a positive constant determined uniquely by the rings
${\mathbb A}$ and ${\mathbb A}^{\! \ast}$. Equivalently,
\begin{equation}\label{E314}
\left(1 + \frac{\eta_{_H}^2}{n-1}\right)^\frac{n-2}{2} \left(1-
\eta_{_H}^2 \right) = \frac{c}{H^n} \leqslant 1
\end{equation}
In particular, $c\leqslant [H(r)]^n = r_\ast^n$. Equation
(\ref{E314}) suggests that we should consider the nonnegative solution
$\eta=\eta(t)$  to the equation
\begin{equation}
\left(1+\frac{\eta^2}{n-1}\right)^\frac{n-2}{2} (1-\eta^2)=
\frac{c}{t^n} \leqslant 1\, , \; \; r_\ast < t < R_\ast
\end{equation}
There is exactly one such solution. The values of $\eta$ lie in the interval
$[0, 1]$.



Now, let $h : {\mathbb A}
\overset{\textnormal{\tiny{onto}}}{\longrightarrow} {\mathbb A}^{\!
\ast}$ be any  homeomorphism of annuli preserving both orientation and
order of the boundary components. We assume that $h \in {\mathscr
W}^{1,n}({\mathbb A}, {\mathbb A}^{\! \ast} )$. For each $x\in
{\mathbb A}$,  we apply Lemma \ref{le5} with $X=
\left|h_{_N}(x)\right|$, $Y= \left|h_{_T}(x)\right|$ and $\alpha=
\eta \big(|h(x)|\big)\leqslant 1$, to obtain the point-wise
inequality
\begin{eqnarray}\label{45tah}
\norm Dh \norm^n &=& \left[\left|h_N\right|^2 +(n-1) \left|h_T\right|^2 \right]^\frac{n}{2} \nonumber \\
&\geqslant & (n-1)^\frac{n}{2}\,  c
\,\frac{\left|h_T\right|^n}{|h|^n} +b\big(|h|\big)
\left|h_N\right|\, \left|h_T\right|^{n-1}
\end{eqnarray}
The coefficient $b \big(|h|\big)$ comes from (\ref{eqb}) where we take
$\alpha= \eta \big(|h|\big)$. An important fact about
$b=b\big(|h|\big)$ is that we have  equality at (\ref{45tah})  if
$\left|h_N \right|= \eta (|h|)\, \left|h_T \right|$. This is exactly
happening for the radial $n$-harmonic map at (\ref{E312}), by the
definition of the constant $c$. Let us integrate (\ref{45tah}) over
the ring ${\mathbb A}$. For the last term  we may apply the
 lower bound at (\ref{E184}). For the first term in
the right hand side of (\ref{45tah}) we use H\"older's inequality, and then~\eqref{E186}.
\begin{eqnarray}\label{E305}
\int_{\mathbb A} \norm Dh \norm^n &\geqslant & (n-1)^\frac{n}{2} c \left(\int_{\mathbb A} \frac{\left|h_T\right|^{n-1}}{|x|\, |h|^{n-1}}\right)^\frac{n}{n-1}  \left(\int_{\mathbb A} \frac{dx}{|x|^n}\right)^\frac{-1}{n-1} + \nonumber \\ &\; & +\,  \omega_{n-1} \int_{r_\ast}^{R_\ast} \tau^{n-1} b(\tau ) \, d \tau   \\
&\geqslant &    (n-1)^\frac{n}{2} \, c \, \mbox{Mod}\, {\mathbb A}\,
+ \, \omega_{n-1} \int_{r_\ast}^{R_\ast} \tau^{n-1} b(\tau ) \, d
\tau \nonumber
\end{eqnarray}
Finally, observe that we have equalities in all estimates for the
radial stretchings. Thus
\begin{equation}\label{E273}
\int_{\mathbb A} \norm Dh \norm^n \geqslant \int_{\mathbb A} \norm
Dh^\circ \norm^n
\end{equation}
as claimed. A proof of the uniqueness statement is postponed until Chapter~\ref{SecUni}. We only record further use  that the equality at (\ref{E273})
yield
\begin{equation}\label{E273p}
\left|h_N\right|\, \left|h_T\right|^{n-1}=J(x,h)
\end{equation}
\vskip0.3cm

Next we turn to the case when ${\mathbb A}^{\! \ast}$ is
significantly  thinner than ${\mathbb A}$.

\section[Proof of Theorem 1.13]{Proof of Theorem \ref{ThNe4}}
Here we assume that the lower Nitsche bound fails; that is, $$\mbox{Mod}\, {\mathbb
A}^{\! \ast} < {\mathcal N}_\dag(\mbox{Mod}\, {\mathbb A}) $$ We split $\mathbb A$  into two concentric annuli
\begin{equation}
{\mathbb A} = {\mathbb A}(r,R)={\mathbb A}(r,1) \cup {\mathbb A}[1
,R)
\end{equation}
where
\begin{equation}
\mbox{Mod}\, {\mathbb A}^{\! \ast} = {\mathcal N}_\dag(\mbox{Mod}\,
{\mathbb A}[1, R))
\end{equation}
 Let $\aleph : {\mathbb A}[1, R) \rightarrow {\mathbb A}^{\!
\ast}$ denote by critical Nitsche map. This is the radial
$n$-harmonic function
\begin{equation}
\aleph (x)= H\big(|x|\big) \frac{x}{|x|}
\end{equation}
determined by the characteristic equation
\begin{equation}
\left[H^2 + \frac{t^2 \dot{H}^2}{n-1}  \right]^\frac{n-2}{2}
\left(H^2 - t^2 \dot{H}^2\right) \equiv c
\end{equation}
Evaluating it at $t=1$ yields $c=[H(1)]^n=r_{\! \ast}^n$. Now
consider any permissible map  $h: \mathbb A \rightarrow {\mathbb
A}^{\! \ast}$ in ${\mathscr P}({\mathbb A}, {\mathbb A}^\ast)$. We
may use inequality (\ref{E305}) with $c=r_{\! \ast}^n$ to obtain
\begin{equation}\label{E278}
\int_{\mathbb A} \norm Dh \norm^n \geqslant
(n-1)^\frac{n}{2}r_\ast^n \; \mbox{Mod}\, {\mathbb A} + \omega_{n-1}
\int_{r_\ast}^{R_\ast} \tau^{n-1} b(\tau)\, d\tau
\end{equation}
On the other hand we have the following identity for the critical
Nitsche map
\begin{equation}
\int_{\mathbb A[1, R]} \norm D\aleph \norm^n =
(n-1)^\frac{n}{2}r_{\! \ast}^n\, \mbox{Mod}\; {\mathbb A}[1, R) +
\omega_{n-1} \int_{r_\ast}^{R_\ast} \tau^{n-1} b(\tau)\, d\tau
\end{equation}
Hence (\ref{E278}) takes the form
\begin{equation}
\int_{\mathbb A} \norm Dh \norm^n \geqslant
(n-1)^\frac{n}{2}r_\ast^n\; \mbox{Mod}\, {\mathbb A}(r, 1) +
\int_{\mathbb A[1, R)} \norm D\aleph \norm^n
\end{equation}
As a final step we notice that the first term is precisely equal to
the energy of the hammering map
\begin{equation}
(n-1)^\frac{n}{2}r_\ast^n\, \mbox{Mod}\; {\mathbb A}(r, 1)=
\int_{\mathbb A(r, 1)} \norm Dg \norm^n\, , \; \; g(x)=r_\ast
\frac{x}{|x|}
\end{equation}
We now glue $\aleph$ and $g$ along the sphere $|x|=1$ to obtain a
 map $h^\circ : {\mathbb A} \rightarrow {\mathbb A}^\ast$
\begin{equation}
h^\circ=\begin{cases}g(x) & \; \; \mbox{ on } \mathbb A(r, 1) \\
\aleph (x) & \; \; \mbox{ on } \mathbb A[1, R)
\end{cases}
\end{equation}
This map minimizes the conformal energy. Indeed,
\begin{equation}
\int_{\mathbb A} \norm Dh \norm^n \geqslant \int_{\mathbb A} \norm
Dh^\circ \norm^n
\end{equation}
for every $h\in {\mathscr P}({\mathbb A}, {\mathbb A}^\ast)$. The
question of uniqueness is discussed in Chapter \ref{SecUni}.
\chapter{Expanding pair, $\mbox{Mod}\, {\mathbb A}^{\! \ast} >
\mbox{Mod}\, {\mathbb A}$}\label{Sec33}

Such annuli determine uniquely a radial $n$-harmonic map $h^\circ :
{\mathbb A} \rightarrow {\mathbb A}^\ast$ in the class ${\mathcal
H}_-$,
\begin{equation}
h^\circ (x)= H\big(|x|\big)\, \frac{x}{|x|} \hskip2cm \mbox{ where }
\; {\mathcal L}H \equiv c <0
\end{equation}
If ${\mathbb A}^\ast$ is too fat, it will latter become clear that in
dimensions $n \geqslant 4$ this radial $n$-harmonic mapping is not
the minimum energy solution.
\section{Within the  bounds, $\mbox{Mod}\, {\mathbb A} <\mbox{Mod}\, {\mathbb A}^{\! \ast}
\leqslant {\mathcal N}^\dag (\mbox{Mod}\, {\mathbb A})$ }
Here the upper Nitsche function ${\mathcal N}^\dag = {\mathcal N}^\dag (t)$, $t>0$, is determined uniquely by requiring that the inequality $\mbox{Mod}\, {\mathbb A}^{\! \ast}
\leqslant {\mathcal N}^\dag (\mbox{Mod}\, {\mathbb A})$ be equivalent to

\begin{equation}\label{equ310}
\frac{\left(\eta^2_{_H} +n-1\right)^\frac{n-2}{2} \left(\eta^2_{_H}
-1\right)}{\eta^n_{_H}} \leqslant 1
\end{equation}
where  $\eta_{_H} = \frac{t \dot{H}}{H}$ is the elasticity function
of $h^\circ : {\mathbb A} \rightarrow {\mathbb A}^\ast$. A somewhat explicit formula for ${\mathcal N}^\dag$ is given in~\eqref{equ320p} after an analysis of~\eqref{equ310}. It may be worth noting in advance that $\mathcal N^\dag \equiv \infty$ in dimensions $n=2,3$. Thus the condition $\mbox{Mod}\, {\mathbb A}^{\! \ast}
\leqslant {\mathcal N}^\dag (\mbox{Mod}\, {\mathbb A})$  is actually void for $n=2,3$. By the definition of $\alpha_n$ at (\ref{E222st}), condition~\eqref{equ310} is equivalent to
\begin{equation}\label{eristar}
\eta_{_H} (t) \leqslant \alpha_n \, , \hskip2cm \mbox{for all } r
\leqslant t \leqslant R
\end{equation}
Since we are in the expanding case the function $\eta_{_H}$ is
decreasing. We therefore need only assume that
\begin{equation}\label{E522}
1< \eta_{_H} (r) \leqslant \alpha_n
\end{equation}
Recall that $\alpha_2=\alpha_3=\infty$, so  this
condition poses no restriction on ${\mathbb A}^\ast$ in dimensions
$n=2$ and $3$.

\begin{figure}[!h]
\begin{center}
\psfrag{n1}{\tiny ${ \frac{n-1}{n-3} }$} \psfrag{an}{\tiny
${\alpha_n }$} \psfrag{s}{\small ${ s }$} \psfrag{func}{\tiny ${
\frac{(s^2+n-1)^2(s^2-1)}{s^n} }$} \psfrag{1}{\tiny ${ 1 }$}
\includegraphics*[height=1.4in]{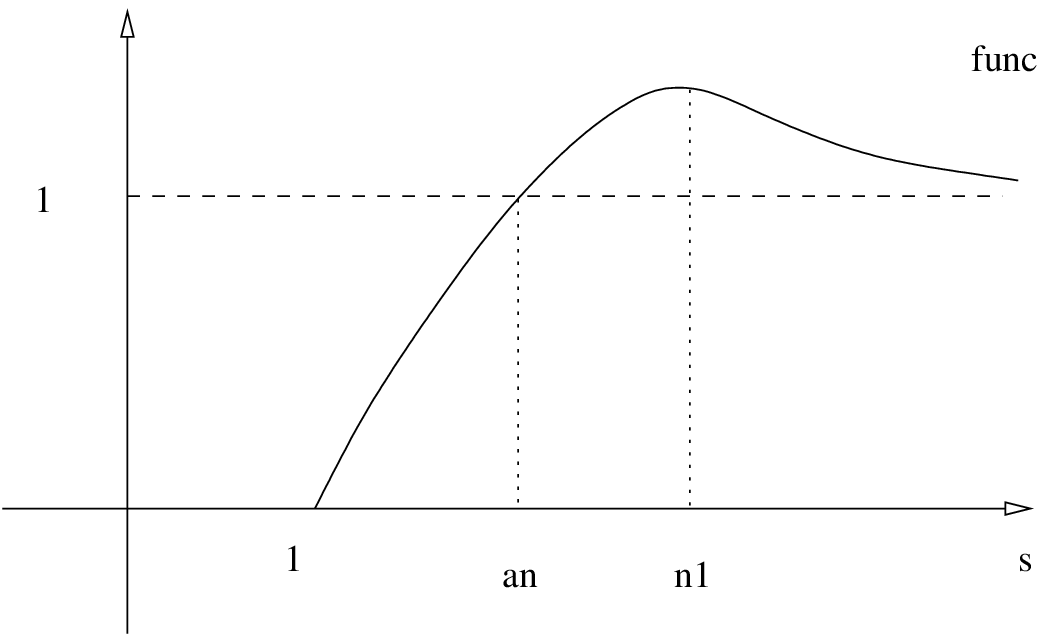}
\caption{The bounds for the  elasticity of $h^\circ$.}\label{Fig17}
\end{center}
\end{figure}

Although it is not immediately clear the condition  (\ref{E522})
 is in fact a condition on the moduli of ${\mathbb
A}$ and ${\mathbb A}^\ast$ alone, precisely
\begin{equation}\label{E523}
H_- (\gamma_n) \cdot \mbox{\small{$\frac{R_\ast}{r_\ast}$}}
\leqslant H_- \left(\gamma_n  \mbox{\small {$\frac{R}{r}$}}
\right)\; \; \footnote{In dimensions $n=2,3$ we have $\alpha_n =
\infty$, hence $\gamma_n =1$. In this case $H_- (\gamma_n)=0$ and,
therefore, we impose no upper bound for $\mbox{Mod}\, {\mathbb
A}^\ast$.}
\end{equation}
where the number $\gamma_n >1$ (for $n\geqslant 4$) is determined by
the equation
\begin{equation}\label{ES252}
\gamma_n = \Gamma_- \left(\mbox{{\small $\frac{1}{\alpha_n
}$}}\right)
\end{equation}
{\it Proof of (\ref{E523}).} Since $h^\circ$ lies in the class
${\mathcal H}_-$ its  strain function takes the form
\begin{equation}
H(t)= \lambda\, H_-(kt)\, , \hskip0.5cm k>\frac{1}{r}\, ,
\hskip0.5cm \mbox{for } r\leqslant t \leqslant R
\end{equation}
The boundary constrains $H(r)=r_\ast$ and $H(R)=R_\ast$ yield a
system of equations for the parameters $\lambda$ and $k$
\begin{equation}
\lambda\, H_-(kr)=r_\ast \, \hskip0.5cm \mbox{ and } \hskip0.5cm
\lambda\, H_-(kR)=R_\ast
\end{equation}
Once we eliminate $\lambda$ the parameter $k$ is determined
(uniquely) from the equation
\begin{equation}\label{ES251st}
\frac{R_\ast}{r_\ast} = \frac{H_-(kR)}{H_-(kr)}
\end{equation}
To solve this equation   we look at the function
\begin{equation}
Q(k)=\frac{H_-(kR)}{H_-(kr)}\, , \hskip1cm k>\frac{1}{r}
\end{equation}
see (\ref{ES119}). Elementary computation shows that
\begin{equation}
\frac{k\, \dot{Q} (k)}{Q(k)} = \eta_{_-}(Rk)- \eta_{_-}(rk) < 0
\end{equation}
In view of (\ref{ES120}) we see that $Q(k)$ decreases from $+\infty$
to $\frac{R}{r}$ as $k$ runs from $\frac{1}{r}$ to $\infty$. Recall
that we are in the expanding case, $\frac{R_\ast}{r_\ast} >
\frac{R}{r}$. Then by Mean Value Theorem there is unique
$k>\frac{1}{r}$ satisfying (\ref{ES251st}). We now rewrite condition
(\ref{ES252}) as follows:
\begin{eqnarray}
\eta_{_H}(r) &\leqslant & \alpha_n \; \; \;  \Leftrightarrow \; \;
\; \eta_{_-}(k\, r)
\leqslant \alpha_n \; \; \; \Leftrightarrow \nonumber \\
k\, r&\geqslant & \eta^{-1}_{_-} (\alpha_n)=
\Gamma_-\left(\mbox{\small{$\frac{1}{\alpha_n}$}}\right)= \gamma_n
\; \; \; \Leftrightarrow \; \; \;  k \geqslant \frac{\gamma_n }{r}
\end{eqnarray}
Here we have used the identity $\Gamma_-(s)=\gamma_-^{-1}
\left(\mbox{\small{$\frac{1}{s}$}}\right)$, see (\ref{254}) and
(\ref{ES108}). The latter inequality is equivalent to
$Q\left(\mbox{\small{$\frac{\gamma_n}{r}$}}\right) \geqslant Q(k)=
\frac{R_\ast}{r_\ast}$, which is the same as (\ref{E523}). We are now in a position to define
\begin{equation}\label{equ320p}
\mathcal N^\dag (t)= \omega_{n-1} \log H_- \left(\gamma_n \mbox{exp}\, \frac{t}{\omega_{n-1}}\right)- \omega_{n-1} \log H_- (\gamma_n).
\end{equation}

\begin{proposition}\label{pr51}
Let  $h : {\mathbb A} \rightarrow {\mathbb A}^\ast$ be a permissible
map in ${\mathscr P}({\mathbb A}, {\mathbb A}^\ast)$ where
\begin{equation}\label{E51star}
\frac{R}{r} <  \frac{R_\ast}{r_\ast} < \frac{ H_- \left(\gamma_n
\frac{R}{r} \right) }{H_-(\gamma_n)}
\end{equation}
Then
\begin{equation}\label{E51dstar}
\int_{\mathbb A} \norm Dh \norm^n \geqslant \int_{\mathbb A} \norm
Dh^\circ \norm^n
\end{equation}
Equality holds if and only if $h(x)=h^\circ (x)$ up to a conformal
automorphism of ${\mathbb A}$.
\end{proposition}
\proof The characteristic equation (\ref{E231}) for the
map $h^\circ : {\mathbb A} \rightarrow {\mathbb A}^\ast$ defines a
positive constant $q=-(n-1)^\frac{n-2}{2}\, {\mathcal L}H$; that is,
\begin{equation}\label{E258}
q\equiv \left[t^2 \dot{H}^2 + (n-1)H^2 \right]^\frac{n-2}{2}
\left(t^2 \dot{H}^2 -H^2\right) >0
\end{equation}
Equivalently,
\begin{equation}\label{ES256}
\left(n-1 + \eta_{_H}^2\right)^\frac{n-2}{2} \left(\eta_{_H}^2
-1\right) = \frac{q}{H^n}
\end{equation}
Now, consider an arbitrary permissible map  $h: {\mathbb A}
\rightarrow {\mathbb A}^\ast$. We may assume that $h$ is a
homeomorphism preserving orientation and the order of the boundary
components of the annuli. We introduce a function $\eta = \eta (t)$
implicitly defined for $r_\ast \leqslant t \leqslant R_\ast$ by the
equation,
\begin{equation}\label{E517}
\left(n-1 + \eta^2\right)^\frac{n-2}{2} \left(\eta^2 -1\right) =
\frac{q}{t^n}
\end{equation}
Note that the equation (\ref{ES256}) for $\eta_{_H}$ at $t=r$
coincides with (\ref{E517}) for $\eta$ at $t=r_\ast$. Hence
$\eta(r_\ast)=\eta_{_H}(r)$. Since $\eta = \eta (t)$ is decreasing
we find that  for every $r_\ast \leqslant t \leqslant R_\ast$,
\begin{equation}
1< \eta (t) \leqslant \eta (r_\ast ) = \eta_{_H} (r) \leqslant
\alpha_n
\end{equation}
Next we recall the function $a=a(\alpha)$, $1\leqslant \alpha <
\infty$, from the formula (\ref{ES208}). By the definition of
$\alpha_n$ at (\ref{E42}) it follows that
\begin{equation}
a(\eta)= \frac{(n-1+ \eta^2)^\frac{n-2}{2} (\eta^2-1)}{\eta^{\, n}}
\leqslant 1 \hskip1cm \mbox{ for all } r_\ast < t < R_\ast
\end{equation}

\begin{figure}[!h]
\begin{center}
\psfrag{eta}{\tiny ${ \eta(|h|) }$} \psfrag{an}{\tiny ${\alpha_n }$}
\psfrag{rn}{\tiny ${ a(|h|) }$} \psfrag{1}{\tiny ${ 1 }$}
\psfrag{aalph}{\tiny $a=a(\alpha)$}
\includegraphics*[height=1.6in]{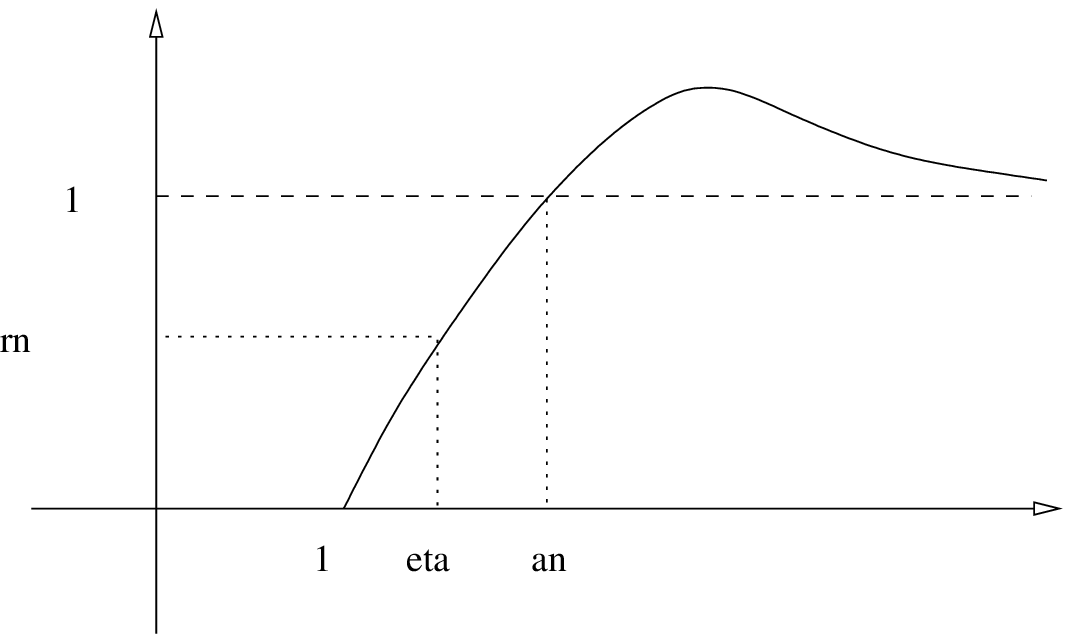}
\caption{The coefficient $a(|h|)$.}
\end{center}
\end{figure}

Now Lemma \ref{le7} applies to $X= \left|h_N \right|$, $Y= \left|h_T
\right|$ and $\alpha = \eta (|h|)$,
\begin{eqnarray}\label{E519}
\norm Dh \norm^n &=& \left[\left|h_N\right|^2 +(n-1) \left|h_T\right|^2 \right]^\frac{n}{2} \nonumber \\
&\geqslant & a\big(|\eta|\big) \left|h_N\right|^n  + b\big(|\eta|\big) \left|h_N\right|\, \left|h_T\right|^{n-1} \nonumber \\
&=& q \left[\frac{\left|h_N\right|}{|h|\,  \eta
\big(|h|\big)}\right]^n + \Phi \big(|h|\big)  \left|h_N\right|\,
\left|h_T\right|^{n-1}
\end{eqnarray}
where $\Phi(|h|)=b\big(\eta (|h|)\big)$ and $b$ is given by
(\ref{defb}). According to Lemma \ref{le7}, equality holds at a
given point $x$ if and only if $\left| h_N (x)\right|= \eta
\big(|h(x)|\big)\, \left|h_T(x)\right|$. In particular, it holds
almost everywhere for $h=h^\circ (x)$, because $\left|h_N^\circ
\right|=\eta_H \left|h_T^\circ \right|$. We now integrate over the
annulus ${\mathbb A}$. The last term at (\ref{E519}) is easily
handled by the  lower bound at (\ref{E184}),
\begin{eqnarray}\label{E300}
\int_{\mathbb A}  \Phi \big(|h|\big)  \left|h_N\right|\,
\left|h_T\right|^{n-1} &\geqslant &  \omega_{n-1}
\int_{r_\ast}^{R_\ast}
\tau^{n-1} \Phi(\tau)\, d\tau \nonumber \\
&=& \int_{{\mathbb A}} \Phi \big(|h^\circ |\big)\, \left|h_N^\circ
\right|  \left|h_T^\circ \right|^{n-1}
\end{eqnarray}
The first term in the right hand side of (\ref{E519}) will be treated by H\"older's inequality in order to apply the  lower bound at
(\ref{E185}). This should be done in a way so one obtains  equality for $h^\circ$. Taking into an account the
 identities
\begin{equation}
\frac{ \left|h^\circ_N\right|}{ \left|h^\circ \right|\, \eta
\big(|h^\circ |\big)} = \frac{\dot{H}}{H\, \eta_{_H}} =
\frac{1}{|x|}
\end{equation}
We proceed as follows
\begin{equation}\label{E307}
\int_{\mathbb A} \left[\frac{ \left|h_N\right|}{ \left|h\right|\,
\eta \big(|h|\big)}\right]^n \geqslant \left(  \int_{\mathbb A}
\frac{ \left|h\right|_N\, dx }{ \left|h\right|\, \eta
\big(|h|\big)\, |x|^{n-1}} \right)^n \left(\int_{\mathbb A}
\frac{dx}{|x|^n}\right)^{1-n}
\end{equation}
The first term in the right hand side is none other than a free Lagrangian  at (\ref{Eq188p}).  This allows us to quickly evaluate
the first term, without getting buried under massive computation.
\begin{equation}
  \int_{\mathbb A} \frac{ \left|h\right|_N\, dx }{ \left|h\right|\, \eta \big(|h|\big)\, |x|^{n-1}} = \int_{\mathbb A}
  \frac{\left|h^\circ\right|_N\, dx}{\left|h^\circ\right|\, \eta \big(|h^\circ |\big)\, |x|^{n-1}}= \int_{\mathbb A} \frac{\dot{H}\big(|x|\big)\, dx }{H \, \eta_{_H} \, |x|^{n-1}} =  \int_{\mathbb A}
  \frac{dx}{|x|^n} \nonumber
\end{equation}
Hence
\begin{equation}\label{E308}
\int_{\mathbb A} \left[\frac{ \left|h_N\right|}{ \left|h\right|\,
\eta \big(|h|\big)}\right]^n \geqslant  \int_{\mathbb A}
\frac{dx}{|x|^n} =\mbox{Mod}\, {\mathbb A}
\end{equation}
In conclusion,
\begin{eqnarray}
 \int_{\mathbb A} \norm Dh \norm^n \geqslant  q \, \mbox{Mod}\, {\mathbb A} + \omega_{n-1}  \int_{r_\ast}^{R_\ast}
 \tau^{n-1} \Phi(\tau)\, d\tau
\end{eqnarray}
with equality occuring for $h^\circ$, as claimed.

The uniqueness results will be proven in more unified manner  in Chapter~\ref{SecUni}.

\subsection{Proof of Theorem~\ref{ThNe2}}
Indeed, the condition $\mbox{Mod}\, \mathbb A^\ast \leqslant \mathcal N^\dag (\mbox{Mod}\, \mathbb A)$ is equivalent to~\eqref{E51star}, in view of the explicit formula for $\mathcal N^\dag$ given in~\eqref{equ320p}.

\subsection{Proof of Theorem \ref{ThNe1}}

The proof is immediate from Proposition~\ref{pr51}  once we recall that $\mathcal N^\dag \equiv \infty$, for $n=2,3$ posing no upper bound for $\mbox{Mod}\, \mathbb A^\ast$.

\chapter{The Uniqueness}\label{SecUni}
There is a comprehensive approach to all our
questions concerning  uniqueness of the conformal energy
\begin{equation}\label{Equ1}
{\mathscr E}_h=\int_{\mathbb A}\norm Dh(x) \norm^n \d x \quad \mbox{ -minimal mappings $h: {\mathbb A}
\rightarrow {\mathbb A}^\ast$. }
\end{equation}
We have already found that within the Nitsche bounds
\begin{equation}\label{Equ2}
{\mathcal N}_\dag (\mbox{Mod}\, {\mathbb A}) \leqslant \mbox{Mod}\,
{\mathbb A}^\ast \leqslant {\mathcal N}^\dag(\mbox{Mod}\, {\mathbb
A})
\end{equation}
the radial $n$-harmonics
\begin{equation}\label{Equ3}
h^\circ (x)= H\big(|x|\big)\, \frac{x}{|x|}
\end{equation}
are among the extremal solutions. We shall show that
\begin{theorem}\label{ThU1}
Below the upper Nitsche bound; that is,
\begin{equation}
\textnormal{Mod}\, {\mathbb A}^\ast   \leqslant {\mathcal N}^\dag
(\textnormal{Mod}\, {\mathbb A})
\end{equation}
every permissible minimizer $h: {\mathbb A} \rightarrow {\mathbb
A}^\ast$ coincides with the radial extremal map modulo conformal
automorphisms of ${\mathbb A}$.
\end{theorem}
Let $h\in {\mathscr P}({\mathbb A}, \, {\mathbb A}^\ast)$ be any
permissible extremal mapping. In all the preceding cases we came to
the following equation as one of the  necessary conditions for $h$ to minimize the energy
\begin{equation}\label{E346}
\left|h_N\right|=\eta\big(|h|\big)\,  \left|h_T\right|
\end{equation}
see (\ref{487}), (\ref{488}), (\ref{45tah})  and (\ref{E519}). Here $\eta =
\eta(\tau)$ is a nonnegative function defined for $r_\ast \leqslant
\tau \leqslant R_\ast$ by the rule
\begin{equation}\label{Equ5}
\tau^n \left(\eta^2+n-1 \right)^\frac{n-2}{2}\left(\eta^2 -1\right)
\equiv c
\end{equation}
Here the constant $c$ comes from the characteristic equation for the
radial extremal map $h^\circ(x)= H\big(|x|\big)\, \frac{x}{|x|} $.
Precisely, we have
\begin{equation}\label{Equ6}
c\equiv
[H(t)]^n\left(\eta^2_{_H}+n-1\right)^\frac{n-2}{2}\left(\eta_{_H}^2-1\right)\,
, \; \; \; \; \mbox{ where } \eta_{_H}(t)=\frac{t\,
\dot{H(t)}}{H(t)}
\end{equation}
Note that  $\eta$ is strictly positive if $c>0$. However, if $c<0$
then $\eta(\tau)$ is strictly increasing and, therefore, can vanish
only at the endpoint $\tau=r_\ast$. This latter situation arises
when  $\eta_{_H}\equiv
0$ in the hammering part of $h^\circ$. Another necessary condition
for a mapping to be extremal takes the form
\begin{equation}\label{mmm}
\left|h_N\right| \, \left|h_T \right|^{n-1} = J(x,h)
\end{equation}
because in each case we used the estimate  (\ref{E184}) , which is sharp only when (\ref{mmm}) holds.   Consequently, the
Cauchy-Green tensor of $h$ becomes a diagonal matrix,
\begin{equation}
 D^\ast h\cdot  Dh = \left[ \begin{array}{cccc}  \left| h_{N} \right|^2 & 0 & \cdots & 0  \\
 0 & \left| h_{T} \right|^2 &   \cdots & 0\\
& &\ddots &  \\
0 & 0 & \cdots & \left| h_{T} \right|^2
\end{array}\right] \bydef {\bold C} (x)
\end{equation}
see (\ref{E777}). Whenever  H\"older's inequalities were used   we always
arrived at one of the following conditions for $h$ to become extremal
\begin{eqnarray}
\left|\, h_T\, \right| = \frac{C_1 \, |h|}{|x|}\, , \; \mbox{ or }
\left|\, h_N\, \right| = \frac{C_2 \, |h|}{|x|}\eta \big(|h|\big)
\end{eqnarray}
where $C_1$ and $C_2$ are constants, see (\ref{E305}) and
(\ref{E307}). No matter which case we look at, the conclusion is that
\begin{eqnarray}\label{abc}
\left|\, h_T\, \right| = \frac{\lambda \, |h|}{|x|} \; \; \mbox{ and
} \; \; \left|\, h_N\, \right| = \frac{\lambda \, |h|}{|x|}\eta
\big(|h|\big)
\end{eqnarray}
for some constant $\lambda >0$, because of (\ref{E346}). This
constant $\lambda$ is the same for all extremal solutions; in fact, $\,\lambda = 1\,$. This can be observed as follows. Since equality holds in either   (\ref{E305}) or  (\ref{E308}), at least one of the integrals
$$
 \int_\mathbb A \frac{|h_T|^{n-1}}{|x|\, |h|^{n-1}} \; \quad \textnormal {and}\;\;\;\;  \int_\mathbb A \Big[\frac{|h_N|}{|x|\, \eta(|h|)}\,\Big]^n
$$
is equal to $\,\textnormal{Mod}\, \mathbb A\,$.  In view of (\ref{abc}) this yields $\,\lambda = 1$.   Hence
\begin{equation}\label{E64}
 D^\ast h\cdot  Dh = \frac{\lambda^2\, |h|^2}{|x|^2}\left[ \begin{array}{cccc} \eta^2 \big(|h|\big) & 0 & \cdots & 0  \\
 0 & 1 &   \cdots & 0\\
& &\ddots &  \\
0 & 0 & \cdots & 1
\end{array}\right]={\bold G}(x,|h|)
\end{equation}
We emphasize that the metric tensor ${\bold G}= {\bold G}(x, \tau)$, viewed as a
function of two variables $x\in {\mathbb A}$ and $\tau \in (r_\ast ,
R_\ast)$, is the same for all extremals. A key step in establishing
 uniqueness  is:

\section{The Point-Cauchy Problem}

We shall consider the Cauchy-Green equation
\begin{equation}\label{Eq1}
D^\ast h(x)\, Dh(x)= {\bold G}(x,|h|)
\end{equation}
for mappings $h:\Omega \rightarrow {\mathbb A}^\ast$ defined in a
domain $\Omega \subset \rn$ and valued in an annulus ${\mathbb
A}^\ast={\mathbb A}\, (r_\ast,R_\ast) \subset \rn$. Here the
function
\begin{equation}
{\bold G}: \Omega \times (r_\ast , R_\ast) \rightarrow {\mathbb
R}^{n \times n}_+ = \begin{cases}\mbox{the space of symmetric}
 \\ \mbox{positive definite matrices}
\end{cases}
\end{equation}
is assumed to be ${\mathscr C}^\infty$-smooth. As for the regularity
of the solutions, we initially assume that $h$ is only continuous. It
then follows from the formula
\begin{equation}
\norm Dh \norm^2 = \mbox{Tr}\, {\bold G}(x,|h|)\in {\mathscr
L}^\infty_{\loc}(\Omega)
\end{equation}
that $h$ is actually locally Lipschitz continuous. One can
look at $\Omega$ as a Riemannian manifold equipped with the positive
definite element of arclength
$$\d s^2=G_{ij}(x)\, \d x^i \otimes \d x^j \; \; \; \; \mbox{ where }\; \left[G_{ij}(x)\right]={\bold G}(x,|h(x)|)$$
In this way $h:\Omega \rightarrow {\mathbb A}^\ast$ becomes a local
isometry with respect to this metric tensor on $\Omega$. At this
point we may appeal to the well known regularity result due to
Calabi and Hartman \cite{CH}. It tells us that if $\, D^\ast h \cdot Dh
\in {\mathscr C}^{k, \alpha}(\Omega, {\mathbb R}_+^{n \times n})$
for some $0< \alpha <1$ and $k=0,1,...$, then $h\in {\mathscr
C}^{k+1, \alpha}(\Omega, {\mathbb R}^{n})$. This result can be applied
repeatedly to infer that in fact $h\in {\mathscr C}^\infty(\Omega, {\mathbb
A}^\ast)$. One more observation is that the system (\ref{Eq1}) is
invariant under linear isometries of the target annulus ${\mathbb A}^\ast$.
Precisely, if $h$ solves this system so does the mapping $T h$ for
every linear isometry $T:\rn \rightarrow \rn$. A priori the system
(\ref{Eq1}) may not admit any local solution. But if it does, the
Riemann curvature tensor of ${\bold G}(x, |h(x)|)$ must vanish. We
shall take advantage of the classical  computation of curvature. We express
the second derivatives of $h$ in terms of its first order
derivatives
\begin{equation}
\frac{\partial h}{\partial x_j \, \partial x_k}= \sum_{\nu =1}^n
\Gamma_{jk}^\nu \frac{\partial h}{\partial x_\nu}
\end{equation}
Here $ \Gamma_{jk}^\nu$ are the  Christoffel symbols; explicitly,
\begin{equation}
 \Gamma_{jk}^\nu = \frac{1}{2} \sum_{i=1}^n G^{i \nu} \left(\frac{\partial {G}_{ij}}{\partial x_k} + \frac{\partial { G}_{ik}}{\partial x_j } - \frac{\partial { G}_{jk}}{\partial x_i}\right)
\end{equation}
where $G^{i\nu}$ are the entries of the inverse matrix to ${\bold
G}$, see for instance \cite[p. 37]{IMb}. It should be noted that the
partials $\frac{\partial {\bold G}}{\partial x_k}$ are computed in accordance with
the chain rule; so in fact we have
\begin{equation}
\frac{\partial {\bold G}}{\partial x_k}= {\bold G}_{x_k}+{\bold
G}_\tau \frac{\partial |h|}{\partial x_k}
\end{equation}
where ${\bold G}_{x_k}$ and ${\bold G}_\tau$ stand for the partial
derivatives of the function $(x, \tau) \rightarrow {\bold G}(x,
\tau)$. Thus, in particular, the Christoffel symbols depend linearly
on $Dh$.  Let us state the general form of the second order
equations obtained in this way
\begin{equation}\label{Eq2}
D^2h(x)=\Phi(x,\, |h|, \, Dh)
\end{equation}
where $\Phi : \Omega \times (r_\ast, R_\ast)\times {\mathbb R}^{n
\times n} \rightarrow {\mathbb R}^{n \times n \times n}$ is a given
${\mathscr C}^\infty$-smooth function. This is in fact a quadratic
polynomial with respect to $Dh$. It should be noted that (\ref{Eq2}) is frame indifferent; that is, it holds for $\,Th\,$ if it holds for $\,h\,$. With these equations at hand we
can now prove the following
\begin{lemma}\label{LeU}$\textnormal{\textsc{(Uniqueness in the point Cauchy problem)}}$
Suppose we are given two solutions $h^\circ$ and $h$ to the
Cauchy-Green equation (\ref{Eq1}) such that
$\left|h(a)\right|=\left|h^\circ (a)\right|$ for some point $a\in
\Omega$. Then there is an isometry $T: {\mathbb R}^n \rightarrow
{\mathbb R}^n$ such that
\begin{equation}
h(x)=T\, h^\circ (x) \; \; \; \; \; \; \mbox{ for all }\; x\in
\Omega
\end{equation}
\end{lemma}
\proof It follows from the equation (\ref{Eq1}) that
\begin{equation}
D^\ast h(a)\, Dh(a)=  D^\ast h^\circ (a)\, Dh^\circ (a)
\end{equation}
Therefore, there is an isometry $T: \rn \rightarrow \rn$ such that
$$Dh(a)=T\circ Dh^\circ (a)=D\big(Th^\circ \big)(a)$$
With a suitable choice of the isometry we may assume without loss of
generality that the first order derivatives of $h$ and $h^\circ$
also coincide at $a$; namely,
\begin{equation}\label{Eq3}
\left|h(a)\right|= \left|h^\circ (a)\right|\; \; \; \; \mbox{ and
}\; \; \; \; Dh(a)=Dh^\circ (a)
\end{equation}
We are going to show that these two equations hold in the entire
domain $\Omega$. Obviously, the set where (\ref{Eq3}) hold
is relatively closed in $\Omega$. Thus, we need only show that this
set is also open. To this end we consider a small ball
$${\mathbb B}=\{\, x\in \rn \, ; \; \; \; |x-a|\leqslant \epsilon \,\}$$
center at $a$ and contained in $\Omega$. We estimate the supremum
norms of $Dh-Dh^\circ$ over ${\mathbb B}$ with $\epsilon$
approaching zero.
\begin{eqnarray}
\norm Dh-Dh^\circ  \norm_{{\mathscr L}^\infty({\mathbb B})} &\preccurlyeq & \epsilon \, \norm D^2h-D^2h^\circ \norm_{{\mathscr L}^\infty({\mathbb B})}   \nonumber \\
&= & \epsilon\, \norm \Phi(x,\, |h|, \, Dh)-\Phi(x, \, |h^\circ|, \, Dh^\circ)\norm_{{\mathscr L}^\infty({\mathbb B})}   \nonumber \\
&\preccurlyeq & \epsilon \, \norm \, |h|-|h^\circ| \,
\norm_{{\mathscr L}^\infty({\mathbb B})} + \epsilon \, \norm \,
Dh-Dh^\circ \, \norm_{{\mathscr L}^\infty({\mathbb B})} \nonumber
\end{eqnarray}
Here and subsequently, the symbol $\preccurlyeq\,$ stands for the inequality with a constant independent of $\epsilon$. This constant varies from line to line.
Noting that
\begin{equation}
\epsilon \, \norm \, |h|-|h^\circ| \, \norm_{{\mathscr
L}^\infty({\mathbb B})} \preccurlyeq  \epsilon \, \norm \,
Dh-Dh^\circ \, \norm_{{\mathscr L}^\infty({\mathbb B})}\,,
\end{equation}
we conclude that $Dh=Dh^\circ$ and $|h|=|h^\circ|$ in ${\mathbb B}$,
as desired.

The equation $Dh(x) \equiv Dh^\circ (x)$ implies that $h^\circ (x)=
h(x)- 2v$,  where $v$ is a constant vector. This combined with the
condition $|h^\circ (x)|\equiv |h(x)|$ yields that $h(x)-v$ is
orthogonal to $v$. If $v$ was not zero the image of $\Omega$ under
$h$ would be an $(n-1)$-hyperplane and, consequently, $J(x,h)\equiv
0$. This is impossible for any solution to the equation (\ref{Eq1}),
because
\begin{equation}
J(x,h)= \sqrt{\mbox{det}\, {\bold G}(x,\, |h|)} \not= 0
\end{equation}
\section[Proof of Theorem 13.1]{Proof of Theorem \ref{ThU1}}
We may assume that the extremal map $h: {\mathbb A} \rightarrow
\overline{{\mathbb A}^\ast}$ preserves both the orientation and the order
of the boundary components of the annuli. For, if not, we compose
$h$ with a suitable conformal automorphism of ${\mathbb A}$. By
virtue of Theorem \ref{ThMono} any extremal map $h$ is
monotone and
\begin{equation}
{\mathbb A}\, (r_\ast, \, R_\ast) = {\mathbb A}^\ast \subset
h({\mathbb A}) \subset \overline{{\mathbb A}^\ast}={\mathbb A}\,
[r_\ast, \, R_\ast]
\end{equation}
Thus, in particular, the set $h^{-1}({\mathbb A}^\ast)$ is a
connected subset of ${\mathbb A}$. We also consider the annulus
${\mathbb A}(\rho, R) \bydef (h^\circ)^{-1}({\mathbb A}^\ast)
\subset {\mathbb A}$. In case $\mbox{Mod}\, \mathbb A^\ast \geqslant \mathcal N_\dag (\mbox{Mod}\, \mathbb A)$  the annulus
${\mathbb A}(\rho , R)$ is  the entire annulus ${\mathbb A}={\mathbb
A}(r,R)$. However, when $\mbox{Mod}\, \mathbb A^\ast < \mathcal N_\dag (\mbox{Mod}\, \mathbb A)$  the inner radius of the
preimage $(h^\circ)^{-1}({\mathbb A}^\ast)={\mathbb A}(\rho , R)$ is
determined in such a way that
$$|h^\circ (x)|=r_\ast \; \; \; \; \mbox{ if } r<|x|\leqslant \rho$$
and
$$r_\ast <\left|h^\circ (x)\right|<R_\ast \; \; \; \; \mbox{ if }  \rho < |x| <R$$
We shall soon see that $h^{-1}({\mathbb A}^\ast)={\mathbb A}(\rho,
R)$ for every extremal map $h$. To this effect let us consider the
union
$$\Omega = h^{-1}({\mathbb A}^\ast) \cup {\mathbb A}(\rho , R) \subset {\mathbb A}$$
which is  connected because $h^{-1}({\mathbb A}^\ast) \cap {\mathbb
A}(\rho , R) \not= \emptyset$. This latter set actually  contains all
points near the outer boundary component of ${\mathbb A}$. Precisely
we have $\lim\limits_{|x|\rightarrow R}|h(x)|=R_\ast$. We now make
use of the free Lagrangian identities (\ref{Jaco})
\begin{equation}
\int_{\mathbb A} \frac{J(x,h)\, \d x}{|h(x)|^n}= \mbox{Mod}\,
{\mathbb A}^\ast = \int_{\mathbb A} \frac{J(x,h^\circ )\, \d
x}{|h^\circ(x)|^n}
\end{equation}
It follows from (\ref{E64}) that
\begin{equation}
\int_{\Omega} \frac{\eta \big(|h(x)|\big)\, \d x}{|x|^n} =
\int_{\Omega} \frac{\eta \big(|h^\circ (x)|\big)\, \d x}{|x|^n}
\end{equation}
The point to make here is that both Jacobians $J(x,h^\circ)$ and
$J(x,h)$ vanish almost everywhere outside $\Omega$. Now, by mean
value property for integrals there exists a point $a\in \Omega$ such
that
\begin{equation}\label{KaT}
\eta \big(|h(a)|\big)= \eta \big(|h^\circ (a)|\big)
\end{equation}
Let ${\mathbb U}\subset \Omega$ denote the set of all points $a\in
\Omega$ for which (\ref{KaT}) holds; ${\mathbb U}$ is certainly
relatively closed.  This set ${\mathbb U}$ is also open. To see this
we first observe that $\eta \big(|h(a)|\big)>0$ for all $a\in
{\mathbb U}$. Indeed, if $a\in {\mathbb A}(\rho, R)$ then $|h^\circ
(a)|>r_\ast$, whereas for $a\in h^{-1}({\mathbb A}^\ast)$ we have
$|h(a)|>r_\ast$. It remains to recall that $\eta(\tau)>0$, whenever
$\tau
>r_\ast$. In other words, the Cauchy-Green tensor of $h$ is positive definite near every point of ${\mathbb U}$. Now it is legitimate to appeal to Lemma \ref{LeU}. Accordingly, there exists a linear isometry
$T : \rn \rightarrow \rn$ so that $h(x)=T\, h^\circ (x)$ near the
point $a$; thus, $|h(x)|=|h^\circ (x)|$ near this point. This shows
that ${\mathbb U}$ is the entire domain $\Omega$. We then have
$$|h(x)|=|h^\circ (x)|>r_\ast\; \; \; \; \mbox{ for all }x\in \Omega$$
Again  by Lemma \ref{LeU} we infer that, upon  suitable adjustment
via isometry of the target annulus it holds
$$h(x)=h^\circ (x) \; \; \; \; \; \; \mbox{ in }\; \Omega$$
As a matter of fact we have
$$h^{-1}({\mathbb A}^\ast)={\mathbb A}(\rho , R)$$
Indeed
$$x\in h^{-1}({\mathbb A}^\ast) \; \Leftrightarrow \; h(x)\in {\mathbb A}^\ast \; \Leftrightarrow \;    h^\circ (x)\in {\mathbb A}^\ast \; \Leftrightarrow \; x\in {\mathbb A}(\rho , R)$$
Thus the uniqueness theorem is proven if ${\mathbb A}(\rho , R)
={\mathbb A}$. This is the case within the Nitsche bounds. In case
below the lower Nitsche bound we look at the remaining region
${\mathbb A}(r, \rho]={\mathbb A}\setminus \Omega$ in which
$|h(x)|\equiv r_\ast$. On the outer boundary of this ring we have
$h(x)=h^\circ (x)= r_\ast \frac{x}{|x|}$. On the other hand, it
follows from the equation (\ref{E346}) that
\begin{equation}
\left|h_N\right|= \frac{\lambda \, |h|}{|x|}\, \eta \big(|h|\big) =
\frac{\lambda\, r_\ast}{|x|} \eta (r_\ast)=0
\end{equation}
Thus $h$ is constant along each ray $t \, \frac{x}{|x|}$, for $r<t
\leqslant \rho$. This means that $h(x)=r_\ast \,
\frac{x}{|x|}=h^\circ (x)$ in ${\mathbb A}\setminus \Omega$.

\chapter{Above the upper Nitsche bound, $n \geqslant
4$}\label{SecCro} This chapter is devoted to a thorough discussion of
the minimization problem when, in dimensions $n \geqslant 4$, the
target annulus is conformally too fat. We shall see that the
extremals are not  radially symmetric. But first we need to
examine related extremal problems for mappings on spheres.
\section{Extremal deformations of the sphere } The study of
extremal deformations of ring domains by using spherical coordinates
leads us to a variational problem for mappings of the unit sphere.
This problem, to be explored later, is the following. Among all
homeomorphisms $\Phi : {\mathbb S}^{n-1} \rightarrow {\mathbb
S}^{n-1}$ of Sobolev class ${\mathscr W}^{1,n}({\mathbb S}^{n-1} ,
{\mathbb S}^{n-1})$ find the one which minimizes the energy integral
\begin{equation}\label{E101}
{\mathcal T}[\Phi]= \Aint_{{\mathbb S}^{n-1}}\left[\alpha^2
+(n-1) \normav{D\Phi}^2\right]^\frac{n}{2}
\end{equation}
where $\alpha$ is any given number. Here $\normav{D\Phi}$ stands for the
normalized Hilbert-Schmidt norm of the tangent  map
\begin{equation}
D\Phi : T_x\,  {\mathbb S}^{n-1} \rightarrow T_y \, {\mathbb
S}^{n-1}\, , \hskip2cm y=\Phi(x)
\end{equation}
That is,
\begin{equation}
\normav{D\Phi}^2=\frac{1}{n-1} \mbox{Tr}\left[\, D^\ast   \Phi\,  D\Phi\,
\right]
\end{equation}
which equals 1 for the identity map. An obvious question to ask is
whether the identity map
\begin{equation}
id\; :\;  {\mathbb S}^{n-1} \rightarrow {\mathbb S}^{n-1}
\end{equation}
is the minimizer. Naturally it is tempting  to apply Jensen's
inequality.
\begin{equation}
{\mathcal T}[\Phi] \geqslant \left[\alpha^2 +(n-1)\Aint_{{\mathbb
S}^{n-1}} \normav{D\Phi }^2 \right]^\frac{n}{2}
\end{equation}
In dimensions $n=2$ and $n=3$, one may appeal to H\"older's and
Hadamard's inequalities to find that
\begin{equation}
\Aint_{{\mathbb S}^{n-1}} \normav{D\Phi}^2 \geqslant \left(\Aint_{{\mathbb
S}^{n-1}} \normav{D\Phi}^{n-1}\right)^\frac{2}{n-1} \geqslant
\left(\Aint_{{\mathbb S}^{n-1}}J(x, \Phi)\,
dx\right)^\frac{2}{n-1}=1
\end{equation}
As usual $J(x, \Phi)$ denotes the Jacobian determinant of $D\Phi$,
the pullback via $\Phi$ of the standard  $(n-1)$-form $\omega$ on ${\mathbb
S}^{n-1}$.
\begin{equation}
\Aint_{{\mathbb S}^{n-1}}J(x, \Phi)\, dx = \mbox{deg}\, \Phi =1
\end{equation}
Unfortunately, in dimensions greater than 3, the infimum of the
integrals $$\int_{{\mathbb S}^{n-1}}\normav{D\Phi}^2$$ is not
attained for the identity map. This infimum equals zero. Here is a
computation which, in addition to illustrating this fact, provides a
method for constructing more sophisticated examples. Consider a
permissible map $\Phi^\epsilon :  {\mathbb S}^{n-1} \rightarrow
{\mathbb S}^{n-1}$, stretching a spherical cap ${\mathbb
S}^\epsilon$ of radius $\epsilon$ around the north pole onto the
entire sphere ${\mathbb S}^{n-1}$.\footnote{$\Phi^\epsilon$ is a
weak ${\mathscr W}^{1,n}$-limit of homeomorphisms of ${\mathbb
S}^{n-1}$ onto itself.} The rest of the sphere is shrank
 into the south pole.

\begin{figure}[!h]
\begin{center}
\psfrag{ep}{\small ${ \epsilon}$} \psfrag{Se}{\small ${\mathbb
S}^\epsilon $} \psfrag{pi}{${ \Pi }$} \psfrag{N}{\small ${N}$}
\psfrag{S}{\small ${S}$}
\includegraphics*[height=1.8in]{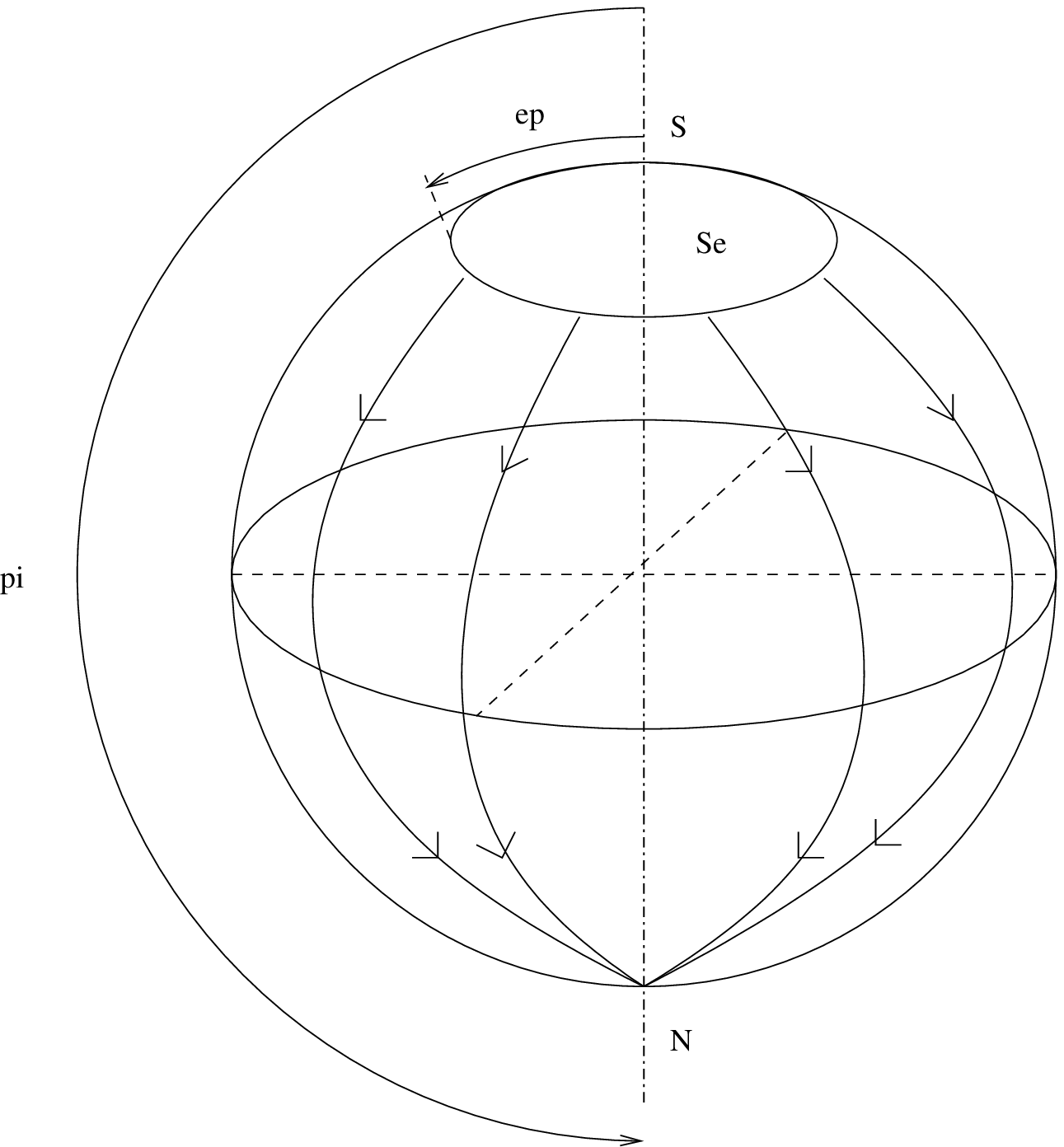}
\caption{A spherical cap to be stretched around the sphere.}
\end{center}
\end{figure}

Elementary geometric considerations show that
\begin{equation}
\normav{D \Phi^\epsilon } \leqslant \begin{cases} \frac{\pi}{\epsilon} & \; \; \; \mbox{ in } {\mathbb S}^\epsilon  \\
0 & \; \; \;  \mbox{ otherwise}
\end{cases}
\end{equation}
whereas $| {\mathbb S}^\epsilon| \approx \epsilon^{n-1}$. Hence, for
every $1\leqslant p < n-1$ we find that
\begin{equation}
\Aint_{{\mathbb S}^{n-1}} \normav{D \Phi^\epsilon (x)}^p\, dx =
O(\epsilon^{n-1-p})\rightarrow 0 \hskip0.5cm \mbox{as } \epsilon
\rightarrow 0
\end{equation}
the infimum being equal to zero. The above computation suggests that
we have to express the integrand of $\mathcal T [\Phi]$,
$$\left[\alpha^2 +(n-1)\normav{ D\Phi}^2\right]^\frac{n}{2}$$
as a convex function in $\normav{D\Phi}^{n-1}$. For this purpose we
introduce
\begin{equation}
F(s)= \left[\alpha^2 +(n-1)s^\frac{2}{n-1}\right]^\frac{n}{2}\, ,
\hskip1.5cm 0 \leqslant s < \infty
\end{equation}
Of particular interest to us will be the lower bound of $F$ by a
convex function $F^\ast$ such that
\begin{equation}
\begin{cases} & F^\ast=F^\ast(s) \leqslant F(s)  \\
& F^\ast (1)=F(1)=(\alpha^2+n-1)^\frac{n}{2}
\end{cases}
\end{equation}
\section{Random variable setting}
It is both illuminating and rewarding to consider  even more
general setting of  the variational integrals such as (\ref{E101}).
We just replace $\normav{D\Phi}^{n-1}$ by a general measurable function.

\begin{figure}[!h]
\begin{center}
\psfrag{an}{} \psfrag{f1}{}
\psfrag{F}{\small ${ F }$} \psfrag{s=}{\tiny ${ s=1 }$}
\psfrag{s}{\small ${ s}$}
\includegraphics*[height=1.6in]{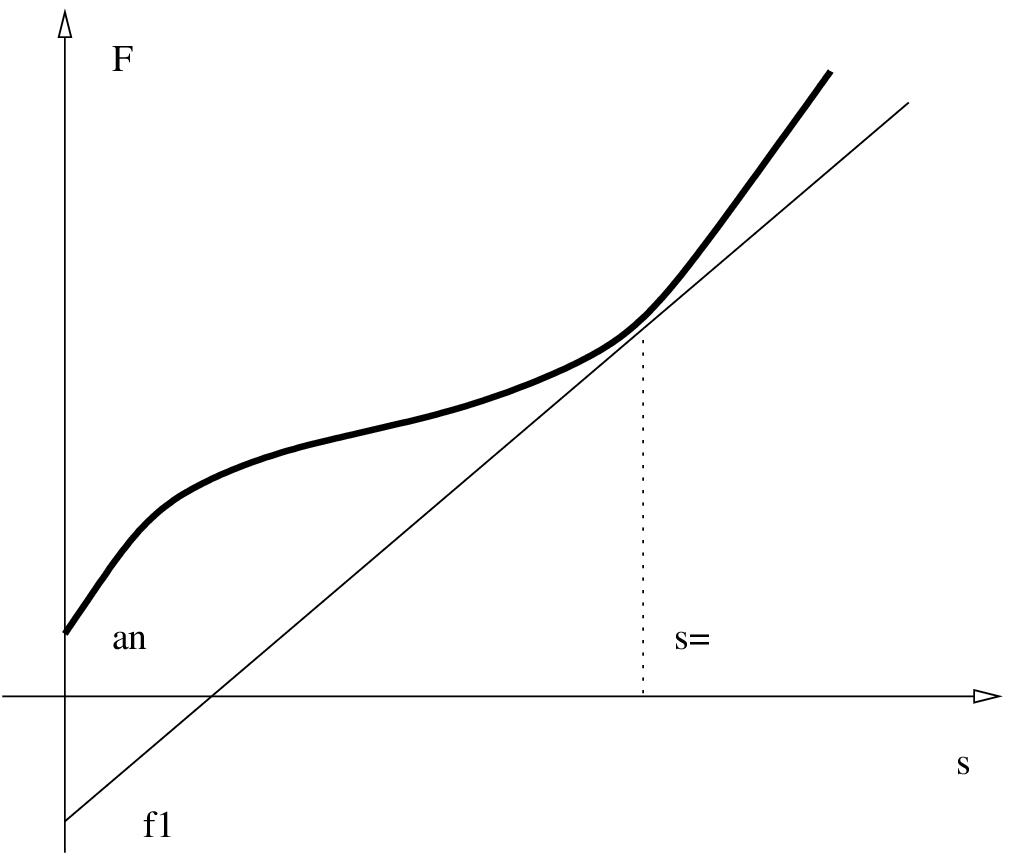}
\caption{}\label{Fig201}
\end{center}
\end{figure}

Let $({\bold S}, \mu)$ be a probability measure space. We shall
consider random variables $X: {\bold S}\rightarrow {\mathbb R}_+$
whose integral mean is at least one;
\begin{equation}
\int_{{\bold S}}X\, d\mu \geqslant 1
\end{equation}
The energy of $X$ is defined by the formula
\begin{equation}\label{E317}
{\mathcal E}[X]= \Aint_{{\bold S}}\left[\alpha^2 +(n-1)
X^\frac{2}{n-1}\right]^\frac{n}{2}\, d \mu
\end{equation}
where $0< \alpha < \infty$ and $n=2,3,...$ We look for the
parameters $\alpha$ for which the constant function $X\equiv 1$ is a
minimizer. That is,
\begin{equation}
\inf {\mathcal E}[X]={\mathcal E}[1]=\left[\alpha^2 +
n-1\right]^\frac{n}{2}
\end{equation}
This is certainly true for all $\alpha$ if $n=2,3$, because $F$ is
convex. However, in higher dimensions $F$ changes concavity. It has
an inflection point at
\begin{equation}
s_\circ = \left(\sqrt{\frac{n-3}{n-1}}\; \alpha \right)^{n-1}
\end{equation}
Precisely, $F$ is concave for $0 \leqslant s \leqslant s_\circ$ and
convex for $s \geqslant s_\circ$. We shall first examine the case
$0< \alpha \leqslant \alpha_n$, where the upper bound $\alpha_n$ is
determined by the equation
\begin{equation}\label{E106}
\left(\alpha_n^2+n-1\right)^\frac{n-2}{2}\left(\alpha_n^2-1\right)=\alpha_n^n
\end{equation}
We have
\begin{equation}
1< \alpha_n < \sqrt{\frac{n-1}{n-3}}
\end{equation}
{\bf Case 1.} $0< \alpha \leqslant \alpha_n$. This means that
\begin{equation}
\left(\alpha^2+n-1\right)^\frac{n-2}{2}\left(\alpha^2-1\right)\leqslant\alpha^n
\end{equation}
It is important to observe that in this case the tangent line of
$F=F(s)$ at the point $s=1$ lies entirely below its graph, see
Figure \ref{Fig201}.

We then consider a convex lower bound $F^\ast (X) \leqslant F(X)$
which coincides with that of $F$ for $s \geqslant 1$ and extends
along the tangent line for $0\leqslant s \leqslant 1$. Jensen's
inequality yields
\begin{eqnarray}
{\mathcal E}[X] &=& \int_{{\bold S}} F(X)\, d\mu \geqslant  \int_{{\bold S}}F^\ast (X)\, d\mu \geqslant F^\ast \left( \int_{{\bold S}} X\, d\mu \right) \nonumber \\
&\geqslant & F^\ast (1)=F(1)={\mathcal E}[1]
\end{eqnarray}
Furhermore, equality occurs if and only if $X \equiv 1$.\\
{\bf Case 2.} $\alpha> \alpha_n$. Thus  the tangent line to $F$ at
$s=1$ intersect the graph of $F=F(s)$ at some point near the origin,
provided $\alpha$ is closed to $\alpha_n$. This geometric
observation suggests that we must look for the extremals  which
assume exactly two values. The best choice turns out to be when
$X$ assumes exactly two values, $0$ and
$\left(\frac{\alpha}{\alpha_n}\right)^{n-1}$. To see this we split
the sample space ${\bold S}$ into two parts
\begin{equation}
{\bold S}= {\bold S}_1 \cup {\bold S}_2
\end{equation}
where
\begin{equation}
\mu \left( {\bold S}_1 \right) =
\left(\frac{\alpha_n}{\alpha}\right)^{n-1} \; \; \mbox{ and }\; \;
\mu \left( {\bold S}_2 \right) =1-
\left(\frac{\alpha_n}{\alpha}\right)^{n-1}
\end{equation}
and define
\begin{equation}
X_\circ = \begin{cases} \left(\frac{\alpha}{\alpha_n}\right)^{n-1} & \; \; \mbox{ on } {\bold S}_1 \\
0 &  \; \; \mbox{ on } {\bold S}_2
\end{cases}
\end{equation}
The energy of $X_\circ$ is, therefore, easily computed as
\begin{equation}
{\mathcal E}[X_\circ ]= \left[\alpha^2 + (n-1)
\frac{\alpha^2}{\alpha^2_n}\right]^\frac{n}{2}\, \mu \left({\bold
S}_1\right) + \alpha^n \, \mu \left({\bold S}_2\right)
\end{equation}
Taking into account the definition of $\alpha_n$ at (\ref{E106}), we
arrive at the formula
\begin{equation}
{\mathcal E}[X_\circ]= \alpha^n + b\,  \alpha
\end{equation}
where
\begin{equation}
b= \frac{n \left(\alpha_n^2 +n-1\right)^\frac{n-2}{2}}{\alpha_n}
\end{equation}
There remains the question as to whether $X_\circ$ possesses the
minimum energy among all random variables of mean at least one. To
this end we appeal to Lemma \ref{le7} in its borderline case when
$a=a(\alpha_n)=1$. By formula (\ref{E212}) we obtain
\begin{equation}
\left[\alpha^2 + (n-1) X^\frac{2}{n-1}\right]^\frac{n}{2} \geqslant
\alpha^n +b\, \alpha X
\end{equation}
for every random variable $X : {\bold S} \rightarrow {\mathbb R}_+$.
Equality holds if and only if $X$ assumes exactly two values $0$ and
$\left(\frac{\alpha}{\alpha_n}\right)^{n-1}$. Upon integration over
${\bold S}$ we conclude with the desired estimate
\begin{equation}
{\mathcal E}[X_\circ]\geqslant \alpha^n + b\,  \alpha
\end{equation}
For the equality to hold the distribution function of the random
variable $X$ must coincide with that of $X_\circ$. In particular,
the constant function $X \equiv 1$ is not the extremal one.

\section{Pulling back a homothety via stereographic projection}
Before we see what can happen if $\alpha > \sqrt{\frac{n-1}{n-3}} $,
let us consider an example of a non-isometry $\Phi: {\mathbb S}^{n-1}
\rightarrow {\mathbb S}^{n-1}$ with
$$\Aint_{{\mathbb S}^{n-1}} \normav{D \Phi}^{n-1}=1$$ Here is a construction
of such a map.

\begin{figure}[!h]
\begin{center}
\psfrag{N}{\small${N}$}
\psfrag{S}{\small${S}$}
\psfrag{x}{\small$x$}
\psfrag{pix}{\tiny${\Pi(x)}$}
\psfrag{y}{\small$y$}
\psfrag{la}{\tiny$\lambda \Pi (x)$}
\psfrag{0}{\tiny$0$}
\includegraphics*[height=1.4in]{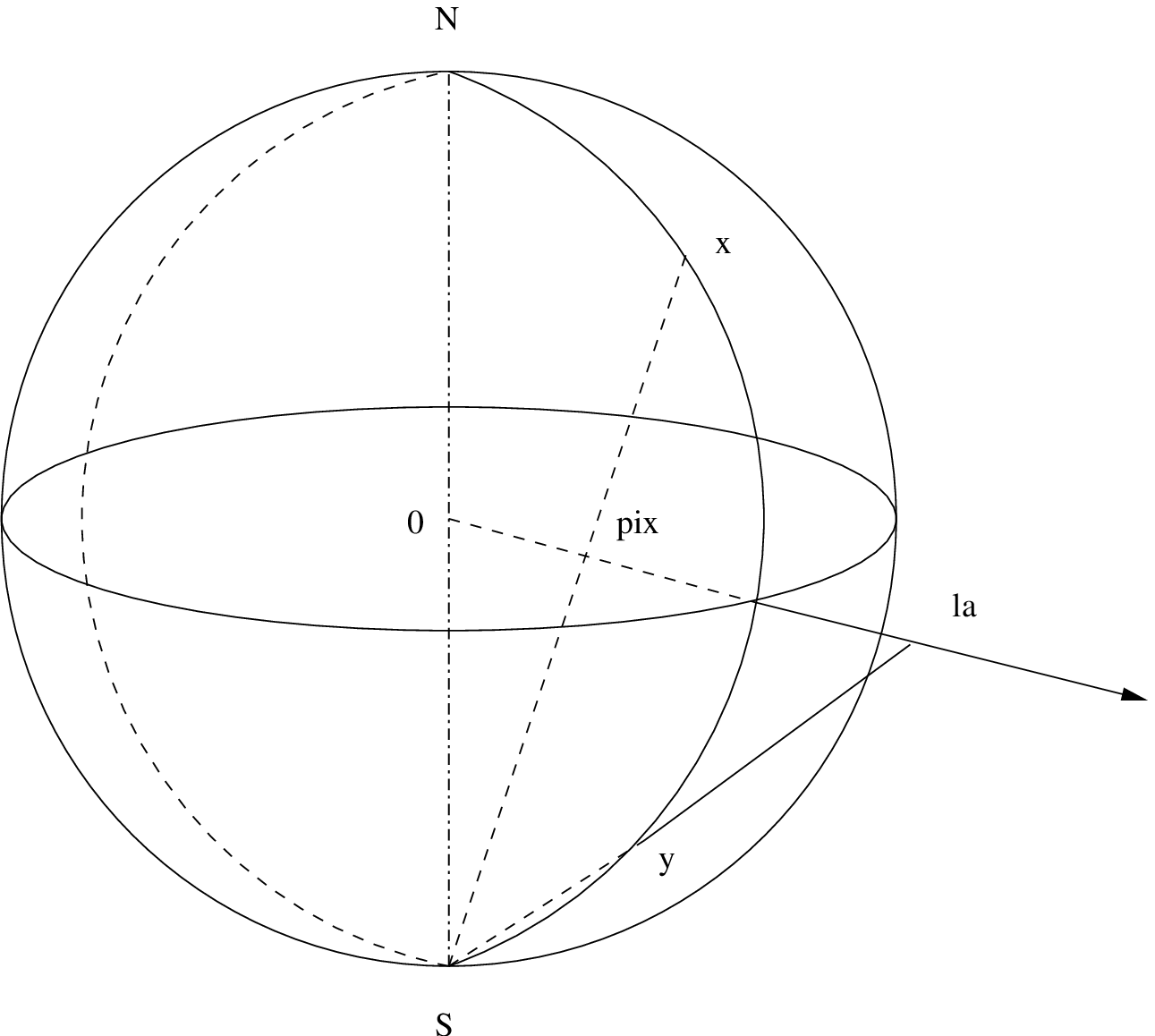}
\caption{Spherical homothety $\Phi : {\mathbb S}^{n-1}\rightarrow  {\mathbb S}^{n-1}$ via stereographic projection.}\label{Fig22}
\end{center}
\end{figure}

Let $\Pi : {\mathbb S}^{n-1} \rightarrow \widehat{\R}^{n-1}$ denote
the stereographic projection of the unit sphere ${\mathbb S}^{n-1}
\subset \rn$ through the south pole onto $\hat{{\mathbb R}}^{n-1}$.
Given any positive number $\lambda$ we consider the homothety
$f=f^\lambda : \hat{{\mathbb R}}^{n-1} \rightarrow \hat{{\mathbb
R}}^{n-1}$ defined by $f^\lambda (x)= \lambda x$ for $x\in
\hat{{\mathbb R}}^{n-1}$. Conjugate to $f^\lambda$ is a conformal
mapping $\Phi = \Phi^\lambda : {\mathbb S}^{n-1} \rightarrow
{\mathbb S}^{n-1}$
\begin{equation}\label{E327}
\Phi^\lambda = \Pi^{-1} \circ f^\lambda \circ \Pi
\end{equation}
We call $\Phi = \Phi^\lambda$ the {\it spherical homothety}. Let $x=(\cos
\theta, {\mathfrak s} \sin \theta)\in {\mathbb S}^{n-1}$ be a point
of longitude ${\mathfrak s} \in {\mathbb S}^{n-2} \subset {\mathbb
R}^{n-1}$ and meridian $0 \leqslant \theta \leqslant \pi$, see Section~\ref{seclalo}. The south
pole corresponds to $\theta = \pi$. Thus $\Phi :  {\mathbb S}^{n-1}
\rightarrow {\mathbb S}^{n-1}$ is longitude preserving. Let $\varphi
= \varphi (\theta )$ denote the meridian of $\Phi (x)= (\cos
\varphi, {\mathfrak s} \sin \varphi)$, see Figure \ref{Fig22}.

\begin{figure}[!h]
\begin{center}
\psfrag{N}{\small${N}$}
\psfrag{S}{\small${S}$}
\psfrag{th}{\small${\theta}$}
\psfrag{var}{\small${\varphi}$}
\psfrag{thh}{\tiny$\frac{\theta}{2}$}
\psfrag{t1}{\tiny$\tan \frac{\theta}{2}$}
\psfrag{t2}{\tiny$\lambda \tan \frac{\theta}{2}=\tan \frac{\varphi}{2}$}
\psfrag{0}{\tiny$0$}
\includegraphics*[height=1.6in]{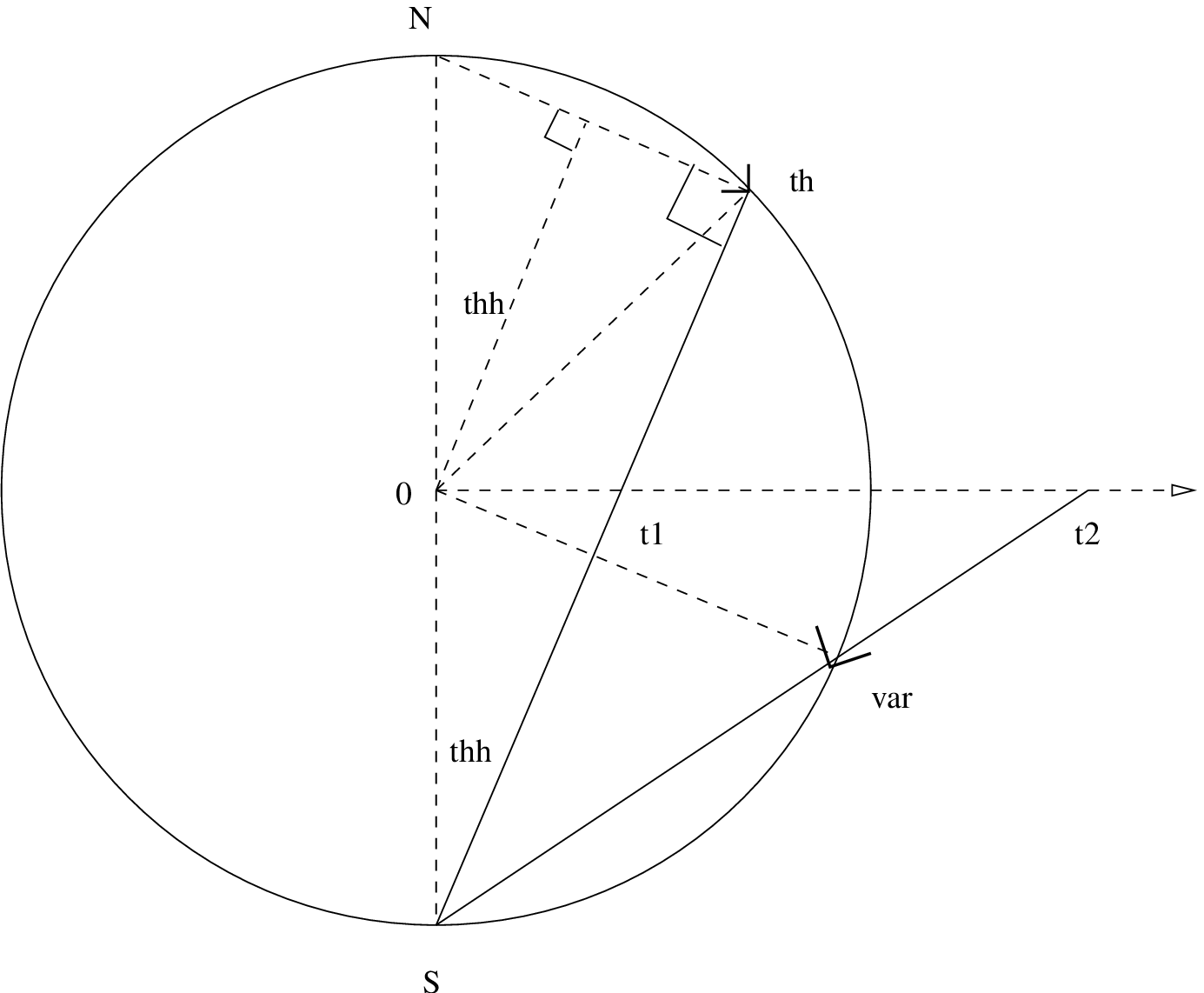}
\caption{Computation of the meridian of $x$ and $\Phi (x)$.}\label{F45}
\end{center}
\end{figure}

\begin{figure}[!h]
\begin{center}
\psfrag{th}{\small${\theta}$}
\psfrag{var}{\small${\varphi}$}
\psfrag{pi}{\tiny${\pi}$}
\psfrag{la=}{\tiny$\lambda=1$}
\psfrag{la<}{\tiny$\lambda<1$}
\psfrag{la>}{\tiny$\lambda>1$}
\includegraphics*[height=1.8in]{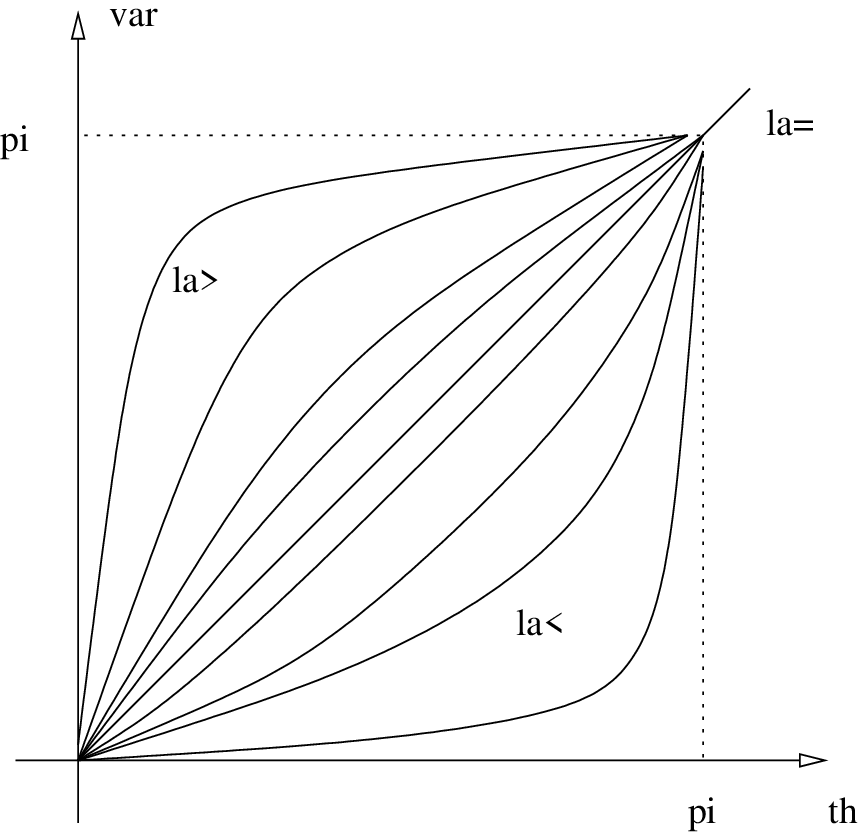}
\caption{The meridian functions.}\label{Fig24}
\end{center}
\end{figure}

Geometric considerations give the following explicit formula
\begin{equation}\label{E91}
\tan \frac{\varphi}{2} = \lambda \tan \frac{\theta}{2}
\end{equation}
see Figure \ref{F45}. Hence
\begin{equation}\label{E92}
\varphi (\theta) = 2 \tan^{-1}\left(\lambda \tan \frac{\theta}{2}\right)
\end{equation}
and
\begin{equation}\label{E93}
\dot{\varphi} (\theta) =\frac{1}{\frac{1}{\lambda} \cos^2 \frac{\theta}{2}+ \lambda \sin^2 \frac{\theta}{2}}>0
\end{equation}
$$\varphi(0)=0\, , \; \varphi(\pi )= \pi \, , \; \dot{\varphi}(0)=\lambda \, , \; \dot{\varphi} (\pi)= \frac{1}{\lambda}$$ see Figure \ref{Fig24}.

Further  differentiation shows that
\begin{equation}\label{E94}
\ddot{\varphi} (\theta) =\frac{1}{2}\left(\frac{1}{\lambda }- \lambda \right)  \left[\, \dot{\varphi} (\theta)\, \right]^2\, \sin \theta
\end{equation}
Thus $\varphi$ is convex if $0< \lambda <1$ and concave if $\lambda >1$. The extreme values of $\dot{\varphi}$ are assumed at the end-points,
\begin{equation}\label{E95}
\min \left\{\lambda, \frac{1}{\lambda}\right\} \leqslant \dot{\varphi }(\theta) \leqslant \max \left\{ \lambda, \frac{1}{\lambda}\right\}
\end{equation}
Implicit differentiation of (\ref{E91}) yields
\begin{equation}
\frac{\dot{\varphi}}{2\, \cos^2 \frac{\varphi}{2}}= \frac{\lambda }{2\, \cos^2 \frac{\theta}{2}} = \frac{\tan \frac{\varphi}{2}\big/ \tan \frac{\theta}{2}}{2\, \cos^2 \frac{\theta}{2}}
\end{equation}
Hence
\begin{equation}\label{E96}
\dot{\varphi}(\theta) = \frac{\sin \varphi}{\sin \theta}
\end{equation}
for every parameter $\lambda >0$.

Next, we return to the mapping $\Phi :  {\mathbb S}^{n-1}
\rightarrow {\mathbb S}^{n-1} $ and its tangent bundle map $D\Phi:
{\bold T}{\mathbb S}^{n-1} \rightarrow {\bold T}{\mathbb S}^{n-1}$,
\begin{equation}
\normav{D\Phi }^2 \bydef \frac{1}{n-1} \left[\dot{\varphi}^2 +(n-2) \frac{\sin^2 \varphi}{\sin^2 \theta}\right]= \frac{\sin^2 \varphi}{\sin^2 \theta}
\end{equation}
Thus
\begin{equation}\label{E97}
\normav{D\Phi}= \frac{\sin \varphi}{\sin \theta}= \dot{\varphi}
\end{equation}
In particular,
\begin{equation}\label{E98}
\min \left\{\lambda, \frac{1}{\lambda}\right\} \leqslant \normav{D\Phi^\lambda}   \leqslant \max \left\{ \lambda, \frac{1}{\lambda}\right\}
\end{equation}
Finally, let us entertain the reader with the following computation
\begin{equation}
 \normav{D\Phi^\lambda}^{n-1}\;=\; \textnormal {det}\,D\Phi^\lambda = \left(\frac{\sin \varphi }{\sin \theta}\right)^{n-1}= \frac{\sin^{n-2}\varphi }{\sin^{n-2}\theta }\; \dot{\varphi}
\end{equation}
Hence we find that
\begin{eqnarray}
\Aint_{{\mathbb S}^{n-1}} \normav{D\Phi^\lambda }^{n-1} &=& \left( \int_0^\pi \dot{\varphi }\, \frac{\sin^{n-2}\varphi} {\sin^{n-2}\theta} \cdot \sin^{n-2}\theta \, d\theta\right) \left(\int_0^\pi \sin^{n-2} \theta \, d\theta \right)^{-1} \nonumber \\
&=&  \frac{\int_0^\pi \sin^{n-2} \varphi \, d\varphi}{\int_0^\pi \sin^{n-2} \theta \, d\theta} =1\;
\end{eqnarray}
as expected.
\section{Back to the variational integral  ${\mathcal T}[\Phi]$}
Recall that
\begin{equation}
{\mathcal T}[\Phi]= \Aint_{{\mathbb S}^{n-1}} |\,\alpha^2 + (n-1) \normav{D
\Phi} ^2\,|^\frac{n}{2}
\end{equation}
for homeomorphisms $\Phi: {\mathbb S}^{n-1} \rightarrow {\mathbb
S}^{n-1}$. We will test it with $\Phi = \Phi^\lambda$, so that
\begin{equation}
\normav{D \Phi}^2 = \frac{1}{n-1} \left[\dot{\varphi}^2 +(n-2) \frac{\sin^2 \varphi}{\sin^2 \theta}\right]\;=\; \frac{\sin^2 \varphi}{\sin^2 \theta}
\end{equation}
Now we recall the function
\begin{equation}
F(s)= \left[\alpha^2 + (n-1)s^\frac{2}{n-1}\right]^\frac{n}{2}\, , \; \; s>0
\end{equation}
The first and second derivatives of $F$ are given by the formulas
\begin{equation}
{\dot F}(s)=n \left[\alpha^2 + (n-1)s^\frac{2}{n-1}\right]^\frac{n-1}{2}s^\frac{3-n}{n-1}>0
\end{equation}
and
\begin{equation}
{\ddot F}(s)=n \left[\alpha^2 + (n-1)s^\frac{2}{n-1}\right]^\frac{n-4}{2}s^\frac{4-n}{n-1}\left(s^\frac{2}{n-2}-\frac{n-3}{n-1}\alpha^2\right)
\end{equation}
Hence $F$ is convex if $n=2,3$. For $n \geqslant 4$ the function $F$ is concave in the interval
\begin{equation}
0\leqslant s \leqslant \left(\frac{n-3}{n-1}\right)^\frac{n-1}{2} \alpha^{n-1}
\end{equation}
From now on, we consider the case
\begin{equation}
\alpha^2 > \frac{n-1}{n-3}\, , \hskip2cm n \geqslant 4
\end{equation}
Given any parameter
\begin{equation}
1< \lambda \leqslant \sqrt{\frac{n-3}{n-1}}\; \alpha
\end{equation}
We examine  the spherical homothety $\Phi = \Phi^\lambda :  {\mathbb
S}^{n-1} \rightarrow {\mathbb S}^{n-1}$. First notice that
\begin{equation}
 \normav{D\Phi}^{n-1} \leqslant \left(\frac{n-3}{n-1}\,  \alpha^2\right)^\frac{n-1}{2}
\end{equation}
point-wise everywhere. Since $ \normav{D\Phi}$ is not
identically equal to one, by concavity argument we conclude with the
following strict inequality
\begin{eqnarray}\label{E349}
\Aint_{{\mathbb S}^{n-1}}\left[\alpha^2 + (n-1)  \normav{D\Phi}^2 \right]^\frac{n}{2} &=& \Aint_{{\mathbb S}^{n-1}}\left[\alpha^2 + (n-1)  \left( \normav{D\Phi}^{n-1}\right)^\frac{2}{n-1} \right]^\frac{n}{2} \nonumber \\
&<& \left[\alpha^2 + (n-1)  \left(  \Aint_{{\mathbb S}^{n-1}} \normav{D\Phi}^{n-1}\right)^\frac{2}{n-1} \right]^\frac{n}{2}  \\
&=& \left[ \alpha^2 + n-1 \right]^\frac{n}{2} \nonumber
\end{eqnarray}
Thus the identity $id:  {\mathbb S}^{n-1} \rightarrow {\mathbb S}^{n-1}  $ is not a minimizer.\\

\section[The failure of  radial symmetry, Proof of Theorem~1.11]{The failure of  radial symmetry, Proof of Theorem~\ref{n4no}}
Throughout this section we make the following standing assumption on
the moduli of ${\mathbb A}$ and ${\mathbb A}^\ast$
\begin{equation}\label{ExE1}
\mbox{Mod}\, {\mathbb A} \leqslant \mbox{Mod}\, {\mathbb A}^{\!\ast}
\end{equation}
Theorems \ref{ThNe1} and \ref{th3} tell us that in dimensions
$n=2,3$ the radial mappings are unique (up to a conformal
automorphism of ${\mathbb A}$) minimizers of both energy
\begin{eqnarray}
\mathscr{E}_h&=& \int_{\mathbb A}\norm Dh \norm^n\, , \hskip1cm h\in
{\mathscr
P}({\mathbb A}\, , \, {\mathbb A}^\ast)\\
\mathscr{F}_h&=& \int_{\mathbb A}\frac{\norm Dh\norm^n}{|h|^n}\, ,
\hskip1cm h\in {\mathscr P}({\mathbb A}\, , \, {\mathbb A}^\ast)
\end{eqnarray}
This is also true for $n\geqslant 4$ provided the target annulus
${\mathbb A}^\ast$ is not too fat, as in Theorem \ref{ThNe2} and
Theorem \ref{th3half}. Here we construct examples to show that these
results do not extend to the full range of moduli at (\ref{ExE1}).
\begin{theorem}
Suppose
\begin{equation}
\textnormal{Mod}\, {\mathbb A}^\ast > \sqrt{\frac{n-1}{n-3}}\,
\textnormal{Mod}\, {\mathbb A}\, , \; \; \;  \; \; n\geqslant 4
\end{equation}
Then
\begin{equation}\label{ExE2}
\inf_{h\in {\mathscr R} ({\mathbb A}\, , \, {\mathbb A}^\ast)
}\int_{\mathbb A} \frac{\norm Dh\norm^n}{|h|^n} > \inf_{h\in
{\mathscr P} ({\mathbb A}\, , \, {\mathbb A}^\ast) }\int_{\mathbb A}
\frac{\norm Dh\norm^n}{|h|^n}
\end{equation}
where ${\mathscr R} ({\mathbb A}\, , \, {\mathbb A}^\ast)$ stands
for the class of radial stretchings $h : {\mathbb A}
\overset{\textnormal{\tiny{onto}}}{\longrightarrow} {\mathbb
A}^\ast$.
\end{theorem}
\proof In Section \ref{Sec52} we have shown that the infimum in the
left hand side is attained for
\begin{equation}
h=h^\alpha(x)=r_{\!\ast}\, r^{-\alpha}|x|^{\alpha-1}x\, , \hskip1cm \alpha=\frac{\textnormal{Mod}
{\mathbb A}^\ast}{\textnormal{Mod}{\mathbb A}}
\end{equation}
The $ \mathscr{F}$-energy of $h^\alpha$ is easily computed
\begin{equation}
{{\mathscr F}}_{h^\alpha} = \left(\alpha^2+ n-1
\right)^\frac{n}{2}\, \mbox{Mod}{\mathbb A}
\end{equation}
We now test the infimum on the right hand side of (\ref{ExE2}) with
the spherical mapping (which is not radial),
\begin{equation}
h(x)= r_{\!\ast}\, r^{-\alpha}|x|^{\alpha-1}\Phi^\lambda \left(\frac{x}{|x|}\right)
\end{equation}
where $\Phi^\lambda : \mathbb S^{n-1} \rightarrow \mathbb S^{n-1}$
is the spherical homothety defined at (\ref{E327}). If $\lambda$ is
sufficiently close to 1, but different from 1, then
\begin{equation}
\Aint_{\mathbb S^{n-1}} \left[\;\alpha^2 +(n-1) \normav{D\Phi^\lambda
}^2\right]^\frac{n}{2} < \left(\alpha^2+ n-1 \right)^\frac{n}{2}
\end{equation}
see (\ref{E349}). Now, the ${\mathscr F}$-energy of $h$ is computed
as follows
\begin{eqnarray}
{\mathscr F}_h&=& \int_{\mathbb A}
\left(\frac{\left|h_N\right|^2}{\left|h\right|^2}+(n-1)\frac{\left|h_T\right|^2}{\left|h\right|^2}
\right)^\frac{n}{2} \nonumber \\
&=& \int_{\mathbb A} \left(\alpha^2 +(n-1) \normav{D\Phi^\lambda
}^2\right)^\frac{n}{2} \frac{dx}{|x|^n} \nonumber \\
&=& \Aint_{\mathbb S^{n-1}} \left(\alpha^2 +(n-1) \normav{D\Phi^\lambda
}^2\right)^\frac{n}{2} \cdot \mbox{Mod}\,{\mathbb A}
\nonumber \\
&<& \left(\alpha^2+ n-1 \right)^\frac{n}{2}\,\mbox{Mod}\,{\mathbb
A}= {\mathscr F}_{h^\alpha}
\end{eqnarray}
as desired. \hskip0.5cm

The spherical homothety $\Phi^\lambda : \mathbb S^{n-1} \rightarrow
\mathbb S^{n-1}$ works as well for the $n$-harmonic energy, though a
computation is  more involved. It results in the proof of Theorem~\ref{n4no}.
\begin{theorem}
Let $n \geqslant 4$ and
\begin{equation}\label{E398}
\delta_n = \left(\frac{\sqrt{n-1}+ \sqrt{n-3}}{\sqrt{n-1}-\sqrt{n-3}}\right)^\frac{1}{2}
\textnormal{\large{exp}}\left[\frac{n-2}{n\sqrt{n-1}}\tan^{-1}\sqrt{n-3}
\right] \geqslant \sqrt{n}.
\end{equation}
Consider the annuli $\mathbb A = \mathbb A (r,R)$ and $\mathbb A^\ast = \mathbb A (r_\ast,R_\ast)$, such that
\begin{equation}\label{Exstd}
1 < \frac{R}{r}< \delta_n \quad \mbox{and} \quad \frac{R_{\!\ast}}{r_\ast}  > \frac{  H_-(\delta_n)}{ H_-\left(\delta_n \frac{r}{R}\right)}.
\end{equation}
Then
\begin{equation}\label{ExE3}
\inf_{h\in {\mathscr R} ({\mathbb A}\, , \, {\mathbb A}^\ast)
}\int_{\mathbb A}\norm Dh\norm^n > \inf_{h\in {\mathscr P} ({\mathbb
A}\, , \, {\mathbb A}^\ast) }\int_{\mathbb A} \norm Dh\norm^n
\end{equation}
\end{theorem}
\proof First we show that the infimum on the left hand side is attained at the radial
$n$-harmonic map
\begin{equation}
h(x)=H\big(|x|\big)\frac{x}{|x|}
\end{equation}
of the form $H(t)=\lambda H_-(kt)$, for suitable parameters $k>\frac{1}{r}$
and $\lambda >0$.  Indeed, the energy of radial mappings is given by
\begin{equation}
{\mathscr E}_h= \omega_{n-1} \int_r^R \left[t^2
\dot{H}^2+(n-1)H^2\right]^\frac{n}{2}\, \frac{\d t}{t}
\end{equation}
By standard convexity arguments the infimum is attained and the
minimizer satisfies the Lagrange-Euler equation. This equation
simply means that $h$ is $n$-harmonic. Since $\mbox{Mod}\, {\mathbb
A}\leqslant \mbox{Mod}\, {\mathbb A}^\ast$ we are in a situation in
which $h$ lies in the class generated by the principal solution
$H_-$; that is, $H(t)=\lambda\, H_-(kt)$, as claimed. These parameters are
uniquely determined by the size of the annuli ${\mathbb A}$ and
${\mathbb A}^{\!\ast}$. Precisely, we have
\begin{equation}
\lambda H_-(kr)=r_\ast \; \; \; \mbox{ and }\; \; \; \lambda
H_-(kR)=R_\ast
\end{equation}
We eliminate the factor $\lambda$ by dividing the above equations, and then find $k$ from the equation
\begin{equation}\label{Exst}
\frac{H_-(kR)}{H_-(kr)}=\frac{R_\ast}{r_\ast}> \frac{R}{r}.
\end{equation}
To find $k$ we observe that  the  function $Q(t) \bydef \frac{H_-(tr)}{H_-(tR)}$,
$\frac{1}{r}< t < \infty$, is increasing from $0$ to
$\frac{r}{R}$. Indeed
\[\frac{t \dot{Q}(t)}{Q(t)}= \frac{tr \dot{H}_- (tr)}{H_-(tr)}  - \frac{tR \dot{H}_- (tR)}{H_-(tR)} = \eta_- (tr)- \eta_- (tR)>0    \]
because the function $\eta_- = \eta_-(s)$ is decreasing, see Figure~\ref{Figuree8} in Chapter~\ref{secrad}. This latter statement is immediate from equation~\eqref{equ123}; that is,
\[\left(1+ \frac{\eta_-^2}{n-1}\right)^\frac{n-2}{2} (\eta_-^2-1)= \frac{-c}{\abs{H(s)}^n}\;, \;\;\;\textnormal{where}\;\;\; H(s) = H_-(s) \;,\;\textnormal{for simplicity} \]
Here $-c>0$ and $\abs{H(s)}$ is increasing. Also note that
 \[ \lim_{t\rightarrow \infty} Q(t)\;=\; \frac{r}{R} \cdot  \lim_{t \to \infty} \frac{ \frac{H(tr)}{tr}   }{ \frac{H(tR)}{tR}  } = \frac{r}{R} \cdot  \frac{\Theta}{\Theta}= \frac{r}{R}\]
see~\eqref{E118}. Now the  equation (\ref{Exst}) has exactly one solution.
Our hypothesis at (\ref{Exstd}) is equivalent to
\begin{equation}
Q(k_\circ) = \frac{r_\ast}{R_\ast} \le Q \left(\frac{\delta_n}{R}\right).
\end{equation}
Since $Q$ is increasing we infer that $k_\circ R \le \delta_n$. Next, since $\eta_-$ is decreasing then for $r<t<R$, we have
\[\eta_- (k_\circ t)> \eta_- (k_\circ R) \geqslant \eta_- (\delta_n). \]
Now we appeal to formula~\eqref{E398} and~\eqref{252}. Accordingly,
\[\Gamma_- \left(\frac{n-3}{n-1}\right)= \delta_n.\]
On the other hand, in view of~\eqref{ES108} and~\eqref{254}
\[\eta_- (\delta_n)= \frac{1}{u(\delta_n)}= \frac{1}{\Gamma^{-1} (\delta_n)}= \sqrt{\frac{n-1}{n-3}}.\]
Therefore, $\eta_- (k_\circ t) \geqslant \sqrt{\frac{n-1}{n-3}}$ for every $r<t<R$. This means that
\begin{equation}
\eta_{_H}(t)> \sqrt{\frac{n-1}{n-3}}\, , \; \; \; \; \mbox{ for all
}\; r \leqslant t \leqslant R
\end{equation}
Now, we return to the computation of the infimum in the left hand
side of (\ref{ExE3}) to obtain
\begin{eqnarray}
\inf_{h\in {\mathscr R} ({\mathbb A}\, , \, {\mathbb A}^\ast)
}\int_{\mathbb A}\norm Dh\norm^n &= &\omega_{n-1} \int_r^R \left[t^2
\dot{H}^2 +
(n-1)H^2\right]^\frac{n}{2}\frac{dt}{t} \nonumber \\
&=& \omega_{n-1} \int_r^R \left[H(t)\right]^n \left[\eta_{_H}^2(t) +
n-1\right]^\frac{n}{2}\frac{dt}{t}
\end{eqnarray}
Then we test the infimum in  the right hand side of (\ref{ExE3})
with the  mapping
\begin{equation}
h_\lambda (x)= H\big( |x|\big)\, \Phi^\lambda
\left(\frac{x}{|x|}\right)
\end{equation}
where, as in the previous case, $\Phi^\lambda : \mathbb S^{n-1}
\rightarrow \mathbb S^{n-1}$ is the spherical homothety. We find
that
\begin{eqnarray}
\inf_{h\in {\mathscr P} ({\mathbb A}\, , \, {\mathbb A}^\ast)
}\int_{\mathbb A}\norm Dh\norm^n &\leqslant &  \int_{\mathbb A}\norm
Dh_\lambda \norm^n
\nonumber \\
&=& \omega_{n-1}  \int_r^R \left[H(t)\right]^n\Aint_{\mathbb S^n}
\left[\eta_{_H}^2(t) + (n-1)| D\Phi^\lambda |^2
\right]^\frac{n}{2}\frac{dt}{t} \nonumber \\
&<& \omega_{n-1} \int_r^R \left[H(t)\right]^n \left[\eta_{_H}^2(t) +
n-1\right]^\frac{n}{2}\frac{dt}{t}  \\
&=& \inf_{h\in {\mathcal R} ({\mathbb A}\, , \, {\mathbb A}^\ast)
}\int_{\mathbb A}\norm Dh\norm^n \nonumber
\end{eqnarray}
The above strict inequality follows from (\ref{E349}) applied to
$\alpha = \eta_{_H}(t)>\sqrt{\frac{n-1}{n-3}}$, where we have chosen
$\lambda \not= 1$ sufficiently close to 1.

\chapter{Quasiconformal mappings between annuli}\label{SecQc}
In this final chapter we present an application of free Lagrangians to obtain
sharp estimates for quasiconformal homeomorphisms $h: {\mathbb A}
\rightarrow {\mathbb A}^\ast$ between annuli ${\mathbb A}={\mathbb
A}(r,R)$ and ${\mathbb A}^\ast={\mathbb A}(r_\ast , R_\ast)$. With
the aid of M\"obius transformations (reflections about the spheres
and $(n-1)$-dimensional hyperplanes) we may assume that $h$
preserves the orientation and the order of the boundary components.
We shall employ the operator norm of the differential matrix,
commonly used in the literature \cite{IMb, Reb, Rib}.
Accordingly, a homeomorphism $h
 : \; {\mathbb A}  \stackrel{\textrm{\tiny{onto}}}{\longrightarrow} {\mathbb
 A}^\ast$ of Sobolev class ${\mathscr W}^{1,1}_{\loc}({\mathbb A}, \, {\mathbb
 A}^\ast)$ is quasiconformal if
 \begin{equation}
 |Dh(x)|^n \leqslant K\, J(x,h)\; \; \; \; \; \; \; \mbox{ a.e. in
 }\; {\mathbb A}
 \end{equation}
The smallest such number $K\geqslant 1$, denoted by $K_O=K_O(h)$, is
called the outer dilatation of $h$. The  inner dilatation is defined
to be the smallest number $K_I=K_I(h)\geqslant 1$ such that
 \begin{equation}
 |D^\sharp h(x)|^n \leqslant K_I\, J(x,h)^{n-1}\; \; \; \; \; \; \; \mbox{ a.e. in
 }\; {\mathbb A}
 \end{equation}
Note the relations
\begin{equation}\label{E205}
K_O\leqslant K_I^{n-1}\; \; \; \mbox{ and }\; \; \; K_I\leqslant
K_O^{n-1}
\end{equation}
For the power stretching $h(x)=|x|^{\alpha -1}x$, $\alpha >0$ we
have
\begin{equation}
K_O=\max \left\{ {\alpha}^{-1}\, , \, \alpha^{n-1}\right\} \; \; \;
\mbox{ and }\; \; \; K_I=\max \left\{\alpha^{1-n}\, , \, \alpha
\right\}
\end{equation}
which shows that both estimates at (\ref{E205}) are sharp. Since
$J(x,h)$ is integrable and $K< \infty$ we see that
 $h\in {\mathscr W}^{1,n}({\mathbb A}, \, {\mathbb A}^\ast)$.
\begin{theorem}\label{QTh1}
Suppose  $h :  {\mathbb A}
\stackrel{\textrm{\tiny{onto}}}{\longrightarrow} {\mathbb A}^\ast$
is a quasiconformal map between annuli. Then
\begin{equation}\label{QuaModIne}
\frac{1}{K_I}\leqslant \left(\frac{\textnormal{Mod}\, {\mathbb
A}^\ast}{{\textnormal{Mod}}{\mathbb A}}\right)^{n-1} \leqslant K_O
\end{equation}
\end{theorem}
This estimate is classic in the theory of quasiconformal mappings, see the pioneering work by F.W. Gehring \cite{Ge}. In the proof below we shall not appeal to any advances in
Quasiconformal Theory or PDEs. In fact, our proof provides a method to tackle
the uniqueness problem in the borderline cases of (\ref{QuaModIne}),
which seems to be unknown in higher dimensions.
\begin{theorem}\label{QTh2}
If one of the two estimates at (\ref{QuaModIne}) becomes equality,
then it is attained only on the corresponding  extremal mappings $h : \; {\mathbb A}
\stackrel{\textrm{\tiny{onto}}}{\longrightarrow} {\mathbb A}^\ast$
of the form
\begin{equation}
h(x)=\sqrt{r_\ast \, R_\ast}\,
\left(\frac{\sqrt{rR}}{|x|}\right)^{\pm \alpha}
\Phi\left(\frac{x}{|x|}\right)\, , \; \; \; \; \; r \leqslant |x|
\leqslant R
\end{equation}
where
\begin{equation}
\alpha = \frac{\textnormal{Mod}\, {\mathbb
A}^\ast}{\textnormal{Mod}\, {\mathbb A}}= \left\{\begin{array}{lll}
K_O^\frac{1}{n-1} \; \; \; \; &  \textrm{ if } \; \textnormal{Mod}\,
{\mathbb A}^\ast \geqslant \textnormal{Mod}\,
{\mathbb A}\\
\;& \\
  K_I^\frac{1}{1-n} &  \textrm{ if } \; \textnormal{Mod}\,
{\mathbb A}^\ast \leqslant \textnormal{Mod}\, {\mathbb A}
\end{array} \right.
\end{equation}
Here the spherical part $\Phi : \; {\mathbb S}^{n-1}
\stackrel{\textrm{\tiny{onto}}}{\longrightarrow} {\mathbb S}^{n-1}$
can be any homeomorphism satisfying
\begin{itemize}
\item[(i)]{volume condition:\; $J(\omega , \Phi)\equiv \pm 1$\,, \, for a.e. $\omega \in {\mathbb S}^{n-1}$}
\item[(ii)]{$\alpha$-contraction condition: meaning that}
\item  {$\,|\,D\Phi(\omega)\,| \leqslant \alpha $ ,\, for the equality in the upper bound of (\ref{QuaModIne})}
\item  {$\,|\,(D\Phi(\omega))^{-1}\,| \leqslant 1/\alpha $ ,\, for the equality in the lower bound of (\ref{QuaModIne})}
\end{itemize}
\end{theorem}
Note,  that the only volume preserving homeomorphisms $\Phi : \;
{\mathbb S}^{1} \stackrel{\textrm{\tiny{onto}}}{\longrightarrow}
{\mathbb S}^{1}$ are isometries. Thus, in dimension $n=2$, the
extremal quasiconformal mappings take the form
\begin{equation}
h(z)=\lambda \, |z|^{\alpha -1} z
\end{equation}
where $\lambda$ is a complex number of modulus $r\,
r_\ast^{-\alpha}$. Example \ref{smper} shows that in higher
dimensions there exist non-isometric homeomorphisms $\Phi : \;
{\mathbb S}^{n-1} \stackrel{\textrm{\tiny{onto}}}{\longrightarrow}
{\mathbb S}^{n-1}$ which are volume preserving and satisfy the
$\alpha$-contraction condition, see also~\cite{Ra10}. 
The exception is the conformal case
of $\alpha =1$ for which the contraction condition $|D\Phi (\omega)|
\leqslant 1$ together with $J(\omega , \Phi) \equiv \pm 1$ imply
that $\Phi$ is an isometry. As a corollary to these observations we
obtain
 Schottky's theorem (1877), \cite{Sc}, in $\rn$.
\begin{theorem}\label{QTh3}
An annulus ${\mathbb A}= {\mathbb A}(r,R)$ can be mapped conformally
onto ${\mathbb A}^\ast= {\mathbb A}(r_\ast ,R_\ast)$, if and only if
$\frac{R}{r}=\frac{R_\ast}{r_\ast}$. Moreover, modulo isometry and
rescaling, every conformal mapping takes the form
\begin{equation}
h(x) \;=\; \left\{\begin{array}{lll} x \; \; \; \; &  \textrm{ the identity}\\
 \frac{x}{|x|^2} &  \textrm{ the inversion} \end{array} \right.
\end{equation}
\end{theorem}
For both Theorems it involves no loss of generality in assuming that
$h$ preserves orientation and the order of boundary components of
the annuli. And we do so from now on.
\subsection{Proof of Theorems \ref{QTh1}}
 Let us first prove the
inequality at the right hand side of (\ref{QuaModIne}). We begin
with the identity (\ref{Eq178}) and, after H\"older's inequality,
use the free Lagrangian (\ref{Jaco})
\begin{eqnarray}\label{KO}
\mbox{Mod}\, {\mathbb A}^\ast &=& \int_{\mathbb A}\frac{\d |h|\wedge
\, \ast \d |x|}{|h|\, |x|^{n-1}} = \int_{\mathbb A}
\left\langle[Dh]^\ast
\frac{h}{|h|}\, , \, \frac{x}{|x|} \right\rangle \, \frac{\d x}{|h|\, |x|^{n-1}}\nonumber \\
&\leqslant &  \int_{\mathbb A} \frac{|Dh|\, \d x}{|h|\, |x|^{n-1}}
\leqslant \left(\int_{\mathbb A}
\frac{|Dh|^n}{|h|^n}\right)^\frac{1}{n} \left(\int_{\mathbb
A}\frac{\d x}{|x|^n}\right)^\frac{n-1}{n} \nonumber \\ & \leqslant &
K_O^\frac{1}{n} \left(\int_{\mathbb A} \frac{J(x,h)\, \d
x}{|h(x)|^n}\right)^\frac{1}{n} \left(\mbox{Mod}\, {\mathbb
A}\right)^\frac{n-1}{n} \\
&=&  K_O^\frac{1}{n} \, \left(\mbox{Mod}\, {\mathbb A}^\ast
\right)^\frac{1}{n} \, \left(\mbox{Mod}\, {\mathbb
A}\right)^\frac{n-1}{n} \nonumber
\end{eqnarray}
as desired.

For the left hand side at (\ref{QuaModIne}) we begin with the
 identity  (\ref{331}) to compute in the
similar fashion that
\begin{eqnarray}
\mbox{Mod}\, {\mathbb A} &=& \int_{\mathbb A} \frac{\d |x|}{|x|}
\wedge h^\sharp \omega = \int_{\mathbb A} \left\langle [D^\sharp h]
\frac{x}{|x|}\, , \, \frac{h}{|h|} \right\rangle \, \frac{\d
x}{|x|\,
|h|^{n-1}} \nonumber \\
&\leqslant & \int_{\mathbb A} \frac{|D^\sharp h|\, \d x}{|x|\,
|h|^{n-1}} \leqslant \left( \int_{\mathbb A}  \frac{|D^\sharp
h|^\frac{n}{n-1}}{|h|^n} \right)^\frac{n-1}{n} \left(\int_{\mathbb
A}\frac{\d x}{|x|^n}\right)^\frac{1}{n} \nonumber \\
&\leqslant &K_I^\frac{1}{n}\, \left(\int_{\mathbb A}\frac{J(x,h)\,
\d x}{|h|^n}\right)^\frac{n-1}{n}\, \left(\mbox{Mod}\, {\mathbb
A}\right)^\frac{1}{n} \\
&=& K_I^\frac{1}{n}\, \left(\mbox{Mod}\, {\mathbb A}^\ast
\right)^\frac{n-1}{n} \, \left(\mbox{Mod}\, {\mathbb
A}\right)^\frac{1}{n}  \nonumber
\end{eqnarray}
as desired.

\subsection{Proof of Theorem \ref{QTh2}}
Concerning the borderline cases, let $h : \; {\mathbb A}
\stackrel{\textrm{\tiny{onto}}}{\longrightarrow} {\mathbb A}^\ast$
denote  the extremal mapping for the right hand side of
(\ref{QuaModIne}). We simplify the matters by assuming that
$r=r_\ast =1$. Note that the power mapping
\begin{equation}
h^\alpha (x)=|x|^{\alpha -1}x\, , \; \; \; \; \; \; \alpha =
\frac{\mbox{Mod}\, {\mathbb A}^\ast}{\mbox{Mod}\, {\mathbb A}}
\end{equation}
is among the extremals.  The borderline equation
\begin{equation}\label{Qdo}
\left(\mbox{Mod}\, {\mathbb A}^\ast\right)^{n-1} = K_O\,
\left(\mbox{Mod}\, {\mathbb A} \right)^{n-1}
\end{equation}
implies that we have equalities everywhere in (\ref{KO}). This
amounts to the following two conditions.
\begin{itemize}
\item[(i)]{$\left\langle D^\ast h \cdot  \frac{h}{|h|}\, , \, \frac{x}{|x|} \right\rangle \equiv \left|D^\ast h \cdot \frac{h}{|h|}\right| \equiv |Dh| \equiv \alpha \, \frac{|h|}{|x|}$\, , \hskip0.2cm $\alpha \; $ - a positive constant}
\end{itemize}
and
\begin{itemize}
\item[(ii)]{$|Dh|^n \equiv K_O\, J(x,h)$}
\end{itemize}
The first set of equations is fulfilled
  if and only if
\begin{equation}\label{QTri}
D^\ast h \cdot \frac{h}{|h|} \equiv |Dh(x)|\, \frac{x}{|x|} \equiv
\alpha \frac{|h|}{|x|} \frac{x}{|x|}\, , \; \; \; \; \alpha \mbox{ -
a positive constant}
\end{equation}
We identify  $\alpha$ from the following equations
\begin{equation}
\frac{\alpha^n}{|x|^n} \equiv \frac{|Dh|^n}{|h|^n} \equiv K_O
\frac{J(x,h)}{|h|^n}
\end{equation}
Upon integration over the annulus ${\mathbb A}$ we obtain
\begin{equation}
\alpha^n \, \mbox{Mod}\, {\mathbb A} = K_O\, \mbox{Mod}\, {\mathbb
A}^\ast
\end{equation}
which in view of (\ref{Qdo}) yields
\begin{equation}
\alpha = \frac{\mbox{Mod}\, {\mathbb A}}{\mbox{Mod}\, {\mathbb A}}=
\sqrt[n-1]{K_O} \geqslant 1
\end{equation}
Now, an interesting nonlinear PDE arises from (i) and, in a more direct manner, from (\ref{QTri})
\begin{equation}
D^\ast h \cdot \frac{h}{|h|^2}= \alpha \frac{x}{|x|^2}
\end{equation}
Rather unexpectedly we can easily solve this equation for $|h|$. Let us
express it as
\begin{equation}
\nabla \log |h|^2 \; = \; \alpha \, \nabla \log |x|^2
\end{equation}
Hence, in view of the normalization $|h(x)|=1$ for $|x|=1$, we find
that
\begin{equation}
|h(x)|=|x|^\alpha
\end{equation}
Now, $h$ takes the form
$$h(x)=|x|^\alpha \Psi (x)\, , \; \; \; \; \mbox{ with }\; \Psi : {\mathbb A} \rightarrow {\mathbb S}^{n-1}$$
The map $h$, being a homeomorphism in ${\mathscr W}^{1,n}({\mathbb
A}, {\mathbb A}^\ast)$, is differentiable almost everywhere, then so
is $\Psi$. From (\ref{QTri}) we also find the operator norm of the
differential
\begin{equation}
|Dh(x)|= \alpha \frac{|h|}{|x|}= \alpha |x|^{\alpha -1} = |Dh^\alpha
(x)|
\end{equation}
This,  in view of (ii), yields
\begin{equation}\label{QDoG}
J(x,h)=\alpha |x|^{n \alpha -n}=J(x,h^\alpha)
\end{equation}
Further examination of the mapping $\Psi : {\mathbb A} \rightarrow
{\mathbb S}^{n-1}$ will reveal that its normal derivative vanishes
almost everywhere. To this end, let $x_\circ \in {\mathbb A}$ be a
point of differentiability of $\Psi$. For small real numbers
$\epsilon$ we consider Taylor's expansions of order one,
\begin{equation}
\Psi(x_\circ + \epsilon x_\circ)= \Psi(x_\circ)+\epsilon
\Psi_N(x_\circ)+o(\epsilon )
\end{equation}
and
\begin{equation}
h(x_\circ + \epsilon x_\circ)-h(x_\circ)= \epsilon
\big[Dh(x_\circ)\big]x_\circ + o(\epsilon)
\end{equation}
Upon elementary computation, letting $\epsilon$ go to zero, we
arrive at the following estimate
\begin{equation}
\left|\alpha \Psi (x_\circ) + \Psi_N(x_\circ)\right| \leqslant
\alpha
\end{equation}
We square it to obtain
\begin{equation}
\alpha^2+2\alpha \, \langle \Psi \, , \, \Psi_N \rangle +
\left|\Psi_N \right|^2 \leqslant \alpha^2
\end{equation}
It is important to observe that $\Psi$ is orthogonal to its
directional derivatives, because $ \langle \Psi , \Psi \rangle =
\left| \Psi \right|^2 \equiv 1$. Thus, in particular, $\Psi_N
(x_\circ)=0$. This simply means that
\begin{equation}
\Psi(x)= \Phi \left( \frac{x}{|x|}\right) \, , \; \; \; \; \; \; \;
\mbox{ where } \; \Phi : {\mathbb S}^{n-1} \rightarrow {\mathbb
S}^{n-1}
\end{equation}
Obviously, $\Phi \in {\mathscr W}^{1,n}({\mathbb S}^{n-1}, {\mathbb
S}^{n-1})$ is a homeomorphism of ${\mathbb S}^{n-1}$ onto itself. It
induces the linear tangent map
\begin{equation}
D\Phi (\omega)\, : \; T_\omega {\mathbb S}^{n-1} \rightarrow
T_\sigma {\mathbb S}^{n-1}\, , \; \; \; \; \; \sigma = \Phi(\omega)
\end{equation}
The Jacobian determinant of $\Phi$, with respect to the standard
volume form on ${\mathbb S}^{n-1}$, will be denoted by
\begin{equation}
J(\omega , \Phi)= \mbox{det}\, D\Phi
\end{equation}
It relates to the Jacobian determinant of $h$ by the rule
\begin{eqnarray}
J(x,h)&=& \alpha |x|^{n\alpha -n}\langle \Phi, \Phi_{T_2} \times ...
\times \Phi_{T_n} \rangle \\
&=& \alpha |x|^{n\alpha -n} J(\omega , \Phi)\, , \; \; \; \omega=
\frac{x}{|x|}
\end{eqnarray}
Thus, in view of (\ref{QDoG}) we see that $\Phi : {\mathbb S}^{n-1}
\rightarrow {\mathbb S}^{n-1}$ is volume preserving. Here  we can
take the cross product of directional derivatives with respect to an
arbitrary positively oriented  orthonormal frame $T_2 , ..., T_n$ in
${\bf T}_\omega {\mathbb S}^{n-1}$. Thus, in view of (\ref{QDoG}), we
obtain
\begin{equation}
J(\omega , \Phi)= \left|\Phi_{T_2} \times ... \times
\Phi_{T_n}\right| \leqslant \left|\Phi_{T_2}\right| \cdots
\left|\Phi_{T_n}\right| \leqslant |\,D \Phi\,|^{n-1}
\end{equation}
Next, we look at the Cauchy-Green tensor of $h$ in terms of the
orthonormal frame $N, T_2 , ..., T_n$ at $T_x {\mathbb A}$
\begin{eqnarray}
D^\ast h \, Dh &=& \left[ \begin{array}{ccc} &- \; \; - \; \; h_N \; \; - \; \; -  & \\
 & - \; \; - \; \; h_{T_2} \; \; - \; \; - & \\
 &\vdots &  \\
 & - \; \; - \; \; h_{T_n} \; \; - \; \; - &
\end{array}\right]\;  \left[ \begin{array}{ccc}  &|& |  \hskip0.8cm | \\
&|& | \hskip0.8cm | \\
&h_N & h_{T_2}   \cdots    h_{T_n} \\
 &|& |  \hskip0.8cm | \\
&|& |  \hskip0.8cm |
\end{array}\right] \nonumber \\
&=& |x|^{2\alpha -2}  \left[ \begin{array}{ccc} &- \; \; - \; \; \alpha \Phi \; \; - \; \; -  & \\
 & - \; \; - \; \; \Phi_{T_2} \; \; - \; \; - & \\
 &\vdots &  \\
 & - \; \; - \; \; \Phi_{T_n} \; \; - \; \; - &
\end{array}\right]\;  \left[ \begin{array}{ccc}  &|& |  \hskip0.8cm | \\
&|& | \hskip0.8cm | \\
&\alpha \Phi & \Phi_{T_2}   \cdots    \Phi_{T_n} \\
 &|& |  \hskip0.8cm | \\
&|& |  \hskip0.8cm |
\end{array}\right]  \\
&=& |x|^{2\alpha -2}  \left[ \begin{array}{cccc} &\alpha^2 &  0 \hskip0.8cm & 0 \\
 & 0 &|\Phi_{T_2}|^2 & \langle \Phi_{T_2} , \Phi_{T_n} \rangle \\
 &\;  & \; &  \\ &\; & \; & \\
& 0  & \langle \Phi_{T_n} , \Phi_{T_2} \rangle & |\Phi_{T_n}|^2
 \end{array}\right]\nonumber
\end{eqnarray}
and
\begin{equation}
D^\ast \Phi \, D \Phi =  \left[ \begin{array}{cccc}
|\Phi_{T_2}|^2 & \ldots & \langle \Phi_{T_2} , \Phi_{T_n} \rangle \\
 &  \vdots &\\
\langle \Phi_{T_n} , \Phi_{T_2} \rangle & \ldots & |\Phi_{T_n}|^2
 \end{array}\right]
\end{equation}
Let $0 \leqslant \lambda_2(x) \leqslant ... \leqslant \lambda_n(x)$
denote the singular values of $D\Phi$, meaning that $\lambda_2^2(x)
\leqslant ... \leqslant \lambda_n^2 (x)$ are eigenvalues of $D^\ast
\Phi \, D \Phi$. It follows from (\ref{KO}) that the numbers
$$\lambda_2(x)|x|^{\alpha -1}\, , \ldots , \lambda_n(x)|x|^{\alpha -1} \; \mbox{ and }\; \alpha |x|^{\alpha -1}$$
are the singular values of $Dh$. Since $\alpha |x|^{\alpha
-1}=|Dh|$, this latter number is the largest singular value,
\begin{equation}
0 \leqslant \lambda_2(x) \leqslant ... \leqslant \lambda_n(x)
\leqslant \alpha
\end{equation}
In particular,
\begin{equation}
\normav{D \Phi (x)} \equiv \lambda_n (x) \leqslant \alpha
\end{equation}
as desired.

This computation also shows that the inequality $|D \Phi (x)|
\leqslant \alpha$ is both sufficient and necessary for the equation
$|Dh(x)|= \alpha |x|^{\alpha -1}$.

A backwards inspection of the above arguments reveals that  $h(x)=
|x|^\alpha \Phi \left(\frac{x}{|x|}\right)$ is an extremal map if
and only if $|D \Phi (x)| \leqslant \alpha$ and $J(\omega , \Phi)
\equiv 1$.

It is rewarding to look at the inverse map $f=h^{-1} :{\mathbb
A}^\ast \overset{\textnormal{\tiny{onto}}}{\longrightarrow} {\mathbb
A}$, which takes the form
\begin{equation}
f(y)=|y|^\frac{1}{\alpha} \, \Phi^{-1} \left( \frac{y}{|y|}\right)
\end{equation}
 We  observe that for
every quasiconformal map $h$ and its inverse $h^{-1}=f$ it holds
\begin{equation}
K_O(x, h)= K_I(y,f)\, , \; \;  \; \; \; \; y=h(x)
\end{equation}
Now in much the same way we find all extremals for the left hand
side of (\ref{QuaModIne}). Theorem \ref{QTh2} follows.

Let us finish this subsection with an example of the volume
preserving mappings of ${\mathbb S}^2 \subset {\mathbb R}^3$.
\begin{example}\label{smper}
For every $\epsilon >0$ there exists a  homeomorphism  $\Phi
\in {\mathscr W}^{1, \infty} ({\mathbb S}^2 , {\mathbb S}^2)$, not
rotation, such that $|\Phi (\omega) - \omega| \leqslant \epsilon$
and $J(\omega , \Phi) \equiv 1$ for almost every $\omega \in
{\mathbb S}^2$.
\end{example}
\subsection{Construction of $\Phi$} 
\begin{figure}[!h]
\begin{center}
\includegraphics*[height=3.0in]{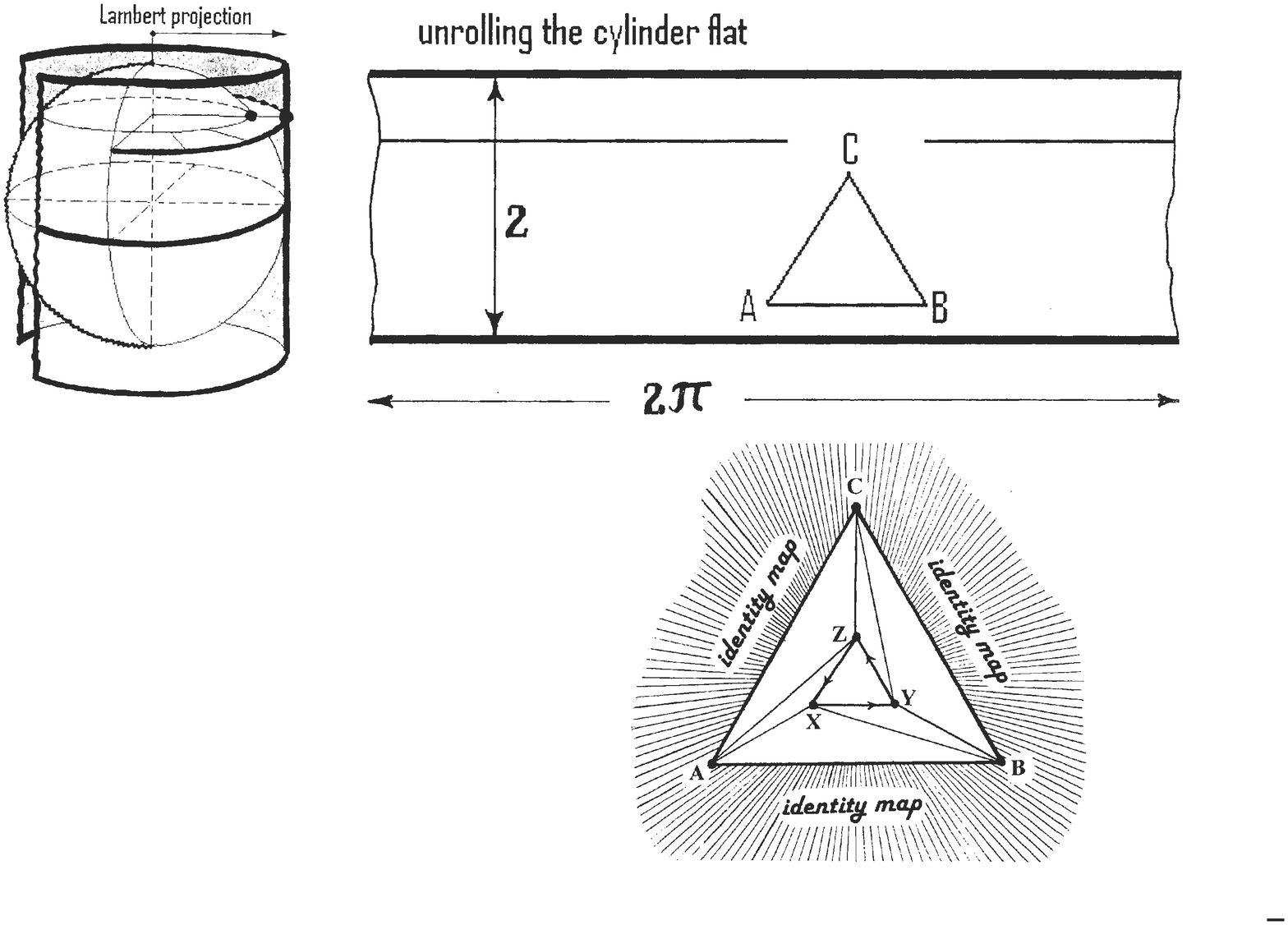}
\caption{Constructing an area-preserving Lipschitz
homeomorphism.}\label{smperex}
\end{center}
\end{figure}

We shall view $\Phi$ as a small
perturbation of $\mbox{id}\, : \; {\mathbb S}^2 \rightarrow {\mathbb
S}^2$. The construction of such perturbation is made in three steps.
In the first step we project the sphere $x^2+y^2+z^2=1$ with two
poles removed, $z \not= \pm 1$  onto the cylinder $x^2+y^2=1\, , \;
-1<z<1$. The horizontal rays from the axis of the cylinder project
the points of the sphere towards the surface of the cylinder. This
projection, sometimes attributed to Archimedes, is well known as
Lambert's Cylindrical Projection (Johann H. Lambert 1772). In the
second step we cut the cylinder along a path from the top to  the
bottom circles and unroll it flat. The third map is a piece-wise
linear perturbation of the triangle $ABC$ inside the flat region,
see Figure \ref{smperex}. It keeps $A, B, C$ fixed while permuting
$X,Y,Z$ in the cyclic way, $X\rightarrow Y\rightarrow Z \rightarrow
X$. It is geometrically clear that such piece-wise linear
deformation preserves the area. Moreover, choosing $X,Y,Z$ close to
the barycenter makes the deformation arbitrarily close to the
identity. Finally we roll it back onto the cylinder and use inverse
of Lambert's projection to end up with the desired area preserving
perturbation of the identity.
\begin{remark} As pointed out by the referee, a simpler example can be given : $(\theta, \varphi)\mapsto (\theta + g(\theta) , \varphi ),$  where $\,\theta\,,\varphi\,$ are the longitude and latitude, and $\,g\,$ is any Lipschitz function.
\end{remark}

\subsection*{Acknowledgements}
We thank Leonid Kovalev for  valuable discussions on the subject of this paper and  the referee for very careful reading of the paper and many useful suggestions.

\backmatter
\bibliographystyle{amsalpha}

\begin{thebibliography}{100}
\bibitem{An}
Antman, S. S. {\em Fundamental mathematical problems in the theory of nonlinear elasticity.} Univ.
Maryland, College Park, Md.,  (1975), pp. 35--54.

\bibitem{AIM}
K. Astala, T. Iwaniec, and G.  Martin,  \textit{Deformations of annuli with
smallest mean distortion}. Arch. Ration. Mech. Anal. {\bf 195}, no.~3 (2010), 899--921.

\bibitem{AIMb}
Astala, K., Iwaniec, T., and Martin, G. J. {\em Elliptic partial differential equations and quasiconformal mappings in the plane}, Princeton University Press, 2009.

\bibitem{AIMO}
Astala, K., Iwaniec, T., Martin, G. J., and Onninen, J. {\em
Extremal mappings of finite distortion.} Proc. London Math. Soc. (3)
{\bf 91} (2005), no. 3, 655--702.

\bibitem{Ba}
 Ball, J. {\em Convexity conditions and existence theorems in nonlinear
elasticity}.  Arch. Rational Mech. Anal. {\bf 63} (1976/77), no. 4,
337--403.
\bibitem{BCO}
Ball, J. M., Currie, J. C., and Olver, P. J. {\em Null Lagrangians,
weak continuity, and variational problems of arbitrary order.} J.
Funct. Anal. 41 (1981), no. 2, 135--174.


\bibitem{BI}
Bojarski, B., and Iwaniec, T. {\em Analytical foundations of the theory of quasiconformal mappings in $
R\sp{n}$.}  Ann. Acad. Sci. Fenn. Ser. A I Math.  {\bf 8}  (1983),  no. 2, 257--324.
\bibitem{CH}
Calabi, E. and Hartman, P. {\em On the smoothness of isometries.}  Duke Math. J.  {\bf 37}  (1970) 741--750.

\bibitem{Ch}
Chen, Y. W.   \textit{Discontinuity and representations of minimal surface solutions}, Proceedings of the conference on differential equations (dedicated to A. Weinstein), pp. 115--138. University of Maryland, College Park, MD (1956).

\bibitem{Ci}
Ciarlet, P. G.  {\em Mathematical elasticity.} Vol. I.
Three-dimensional elasticity. Studies in Mathematics and its
Applications, 20. North-Holland Publishing Co., Amsterdam (1988).

\bibitem{CIKO}
Cristina, J.,   Iwaniec, T, Kovalev, L. V.,  and  Onninen, J. {\em The Hopf-Laplace equation}, arXiv:1011.5934.

\bibitem{CHM}
Cs\"ornyei, M.,  Hencl, S., and Mal\'y, J. {\em Homeomorphisms in the Sobolev space $W^{1,n-1}$}, J. Reine Angew. Math. 644 (2010), 221--235.


\bibitem{DM}
F. Duzaar and G. Mingione, {\em The $p$-harmonic approximation and the regularity of $p$-harmonic maps},
Calc. Var. Partial Differential Equations 20 (2004), no. 3, 235--256.



\bibitem{Du}  Duren, P. {\em Harmonic mappings in the plane},  Cambridge University Press, Cambridge, 2004.
\bibitem{Ed}
Edelen, D. G. B. {\em The null set of the Euler-Lagrange operator.}
Arch. Rational Mech. Anal. 11 (1962) 117--121.

\bibitem{Ge}
Gehring, F. W. {\em Rings and quasiconformal mappings in space}, Trans. Amer. Math. Soc. \textbf{103} (1962), 353--393.


\bibitem{Ha}
Han, Z-C.  {\em Remarks on the geometric behavior of harmonic maps between surfaces}, Elliptic and parabolic methods in geometry (Minneapolis, MN, 1994), 57--66.


\bibitem{HL}
Hardt R. and  Lin. F-H. {\em Mappings minimizing the $L^p$ norm of the gradient},
Comm. Pure Appl. Math. 40 (1987), no. 5, 555--588.



\bibitem{HK}
Hencl, S., and  Koskela, P. {\em Regularity of the inverse of a
planar Sobolev homeomorphism}, Arch. Ration. Mech. Anal. {\bf 180} (2006), no. 1, 75--95.

\bibitem{HKO} Hencl, S.  Koskela, P., and  Onninen, J. {\em
A note on extremal mappings of finite distortion},  Math. Res. Lett. {\bf 12}
(2005), no. 2-3, 231--237.


\bibitem{HN}
Hildebrandt S. and   Nitsche, J. C. C.  \textit{A uniqueness theorem for surfaces of least area with partially free boundaries on obstacles}, Arch. Rational Mech. Anal. {\bf 79} (1982), no.~3, 189--218.


\bibitem{Iw}
Iwaniec, T. {\em On the concept of the weak Jacobian and Hessian.}
Rep. Univ. Jyv\"askyl\"a Dep. Math. Stat., 83, Univ. Jyv\"askyl\"a,
Jyv\"askyl\"a, (2001) 181--205.

\bibitem{IKO1}
Iwaniec, T.,  Koskela, P., and Onninen, J. {\em Mappings of finite
distortion: monotonicity and continuity.}  Invent. Math.   {\bf 144}
(2001),  no. 3, 507--531.

\bibitem{IKKO}
Iwaniec, T.,   Koh, N.-T., Kovalev, L. V.,   and  Onninen, J.
{\em Existence of energy-minimal diffeomorphisms between doubly connected domains}, arXiv:1008.0652.

\bibitem{IKOni1}
Iwaniec, T.,   Kovalev, L. V. and  Onninen, J. {\em Harmonic mappings of an annulus, Nitsche conjecture and its generalizations},  Amer. J. Math. {\bf 132} (2010), no. 5, 1397--1428.

\bibitem{IKOni}
Iwaniec, T.,   Kovalev, L. V. and  Onninen, J. {\em The Nitsche conjecture}, J. Amer. Math. Soc. {\bf 24} (2011), no. 2, 345--373.

\bibitem{IKOaf}
Iwaniec, T.,   Kovalev, L. V. and  Onninen, J.  {\em Harmonic mapping problem and affine capacity},  Proc. Roy. Soc. Edinburgh Sect. A, to appear.

\bibitem{IKOge}
Iwaniec, T.,   Kovalev, L. V. and  Onninen, J. {\em Doubly connected minimal surfaces and extremal harmonic mappings}, J. Geom. Anal., to appear.


\bibitem{IMb}
Iwaniec, T. and Martin, G. J.  {\em Geometric Function Theory and
Non-linear Analysis.} Oxford Mathematical Monographs. 2001.

\bibitem{IOt}
Iwaniec, T. and  Onninen, J. {\em Deformations of finite conformal energy: boundary behavior and limit theorems}, Trans. Amer. Math. Soc., to appear.

\bibitem{IOj}
Iwaniec, T. and Onninen, J. {\em Deformations of finite conformal energy: existence and removability of singularities}, Proc. Lond. Math. Soc. (3) {\bf 100} (2010), no. 1, 1--23.

\bibitem{IS}
Iwaniec, T. and \v Sver\'ak, V. {\em On mappings with integrable dilatation.}  Proc. Amer. Math. Soc.  {\bf 118}  (1993),  no. 1, 181--188.
\bibitem{JS}     Jost, J. and  Schoen, R. {\em  On the existence of harmonic diffeomorphisms},  Invent. Math., {\bf 66},  (1982), 353-359.


\bibitem{Ka}
Kalaj D.,  {\em On the Nitsche conjecture for harmonic mappings in $\R^2$ and $\R^3$.} Israel J. Math. {\bf 150} (2005), 241--251.


\bibitem{L}
Lyzzaik, A. {\em The modulus of the image annuli under univalent harmonic mappings and a conjecture of J.C.C. Nitsche},  J. London Math. Soc., {\bf 64}, (2001), 369--384.

\bibitem{Mib}
Mingione, G. {\em Regularity of minima: an invitation to the dark side of the calculus of variations},
Appl. Math. {\bf 51} (2006), no. 4, 355--426.


\bibitem{Mor}
Morrey, C.B. {\em The Topology of (Path) Surfaces.} Amer. J. Math.
{\bf 57} (1935), no. 1, 17--50.
\bibitem{Mo}
Morrey, C.B. {\em Quasi-convexity and the lower semicontinuity of multiple integrals.} Pacific J. Math.
{\bf 2}, (1952) 25--53.
\bibitem{Mo1}
Morrey, C.B. {\em A variational method in the theory of harmonic integrals. II}, Amer. J. Math. {\bf 78} (1956), 137--170.

\bibitem{Mob}
 Morrey, C.B.  {\em Multiple integrals in the calculus of variations},  Springer-Verlag (1966).


\bibitem{Nitsche}
Nitsche, J.C.C.  {\em On the modulus of doubly connected regions under harmonic mappings},  Amer. Math. Monthly,  {\bf 69}, (1962), 781--782.

\bibitem{Ra10}
Rajala, K. {\em Distortion of quasiconformal and quasiregular mappings at extremal points.} Michigan Math. J. {\bf 53} (2005), no. 2, 447--458.

\bibitem{Reb}
Reshetnyak, Yu. G. {\em Space mappings with bounded distortion.}
American Mathematical
Society, Providence, RI, 1989.
\bibitem{Rib}
Rickman, S. {\em Quasiregular mappings.}  Springer-Verlag, Berlin,
1993.
\bibitem{Sc}
Schottky, F.H. {\em \"Uber konforme Abbildung von mehrfach
zusammenh\"angenden Fl\"ache.}  J. f\"ur Math., {\bf 83}, (1877).


\bibitem{Tu}
Turowski, G. \textit{Behaviour of doubly connected minimal surfaces at the edges of the support surface}, Arch. Math. (Basel) {\bf 77} (2001), no.~3, 278--288.


\bibitem{Uh}
 Uhlenbeck, K. {\em Regularity for a class of non-linear elliptic systems},  Acta Math.  {\bf 138}  (1977), no. 3-4, 219--240.

\bibitem{Ur}
Uraltseva, N. {\em  Degenerate quasilinear elliptic systems},  Zapiski Nauch. Sem. LOMI, {\bf 7}, (1968), 184--222.


 \bibitem{Wa}
 Wang, C. {\em A compactness theorem of $n$-harmonic maps},
Ann. Inst. H. Poincar\'e Anal. Non Lin\'eaire {\bf 22} (2005), no. 4, 509--519.

 \bibitem{W}  Weitsman, A. {\em Univalent harmonic mappings of annuli and a conjecture of J.C.C. Nitsche}, Israel J. Math., {\bf 124},  (2001),  327--331.

\end{thebibliography}

\end{document}